\documentclass[11pt]{scrartcl}

\usepackage[utf8]{inputenc}

\usepackage{config}
\usepackage{titling}

\usepackage[a4paper, top=3cm, bottom=2cm, left=2.5cm, right=2.5cm, includefoot]{geometry}

\usepackage[inline]{enumitem}

\setlist[itemize]{topsep=0pt,partopsep=0pt,itemsep=0pt,parsep=0pt}
\setlist[itemize,1]{label={\small\textbullet}}
\setlist[itemize,2]{label={\tiny\textbullet}}
\setlist[itemize,3]{label=$\cdot$}
\setlist[enumerate]{topsep=0pt,partopsep=0pt,itemsep=0pt,parsep=0pt}
\setlist[enumerate,1]{label=\roman*)}
\setlist[enumerate,2]{label=\alph*)}
\setlist[enumerate,3]{label=\arabic*)}

\hypersetup{
	colorlinks=true,
	linkcolor=AO!65!black,
	citecolor=AO!65!black,
	urlcolor=AppleGreen!65!black,
	bookmarksopen=true,
	bookmarksnumbered,
	bookmarksopenlevel=2,
	bookmarksdepth=3
}

\title{A Flat Wall Theorem for Matching Minors in Bipartite Graphs}
\predate{}
\date{}
\postdate{}

\preauthor{}
\DeclareRobustCommand{\authorthing}{
	\begin{center}
		Archontia C.\@ Giannopoulou\thanks{Supported by the ERC consolidator grant DISTRUCT-648527.} -- National and Kapodistrian University of Athens\\
		\href{mailto:archontia.giannopoulou@gmail.com}{archontia.giannopoulou@gmail.com}\\
		Sebastian Wiederrecht\thanks{Supported by the ANR project ESIGMA (ANR-17-CE23-0010).} -- LIRMM, University of Montpellier\\
		\href{mailto:sebastian.wiederrecht@gmail.com}{sebastian.wiederrecht@gmail.com}
\end{center}}
\author{\authorthing}
\postauthor{}

\begin{document}
\maketitle

\begin{abstract}
A major step in the graph minors theory of Robertson and Seymour is the transition from the Grid Theorem which, in some sense uniquely, describes areas of large treewidth within a graph, to a notion of local flatness of these areas in form of the existence of a large flat wall within any huge grid of an $H$-minor free graph.
In this paper, we prove a matching theoretic analogue of the Flat Wall Theorem for bipartite graphs excluding a fixed matching minor.
Our result builds on a a tight relationship between structural digraph theory and matching theory and allows us to deduce a Flat Wall Theorem for digraphs which substantially differs from a previously established directed variant of this theorem.

\end{abstract}

\section{Introduction}\label{sec:introduction}

The general structure theorem for graphs excluding a fixed graph $H$ as a minor \cite{robertson2003graph} which sits at the heart of the graph minors series by Robertson and Seymour says, in a simplified way, that every graph that excludes $H$ as a minor can be obtained from a surprisingly small list of ingredients:
\begin{itemize}
	\item Graphs of bounded size,
	\item clique sums,
	\item graphs on surfaces,
	\item vortices, and
	\item apex vertices.
\end{itemize}

The proof of this general structure theorem can be summarised into the following three steps:

\begin{enumerate}
	\item Robertson and Seymour provide a structure theorem for $H$-minor free graphs where $H$ is planar.
	This is exactly the \textbf{Grid Theorem} \cite{robertson1986graph}.
	
	\item  Next, they provide a local (or weak) structure theorem for $H$-minor free graphs where $H$ is non-planar.
	This is the \textbf{Flat Wall Theorem} \cite{robertson1995graph}.
	
	\item Finally, they provide a global structure theorem for $H$-minor free graphs where $H$ is non-planar using all five ingredients from the list above \cite{robertson2003graph}.
\end{enumerate}

Each of these three steps gave rise to its own consequences and applications, both algorithmic and theoretical.
Replicating parts of these structural descriptions of proper minor closed families of graphs for other partial orderings of graphs, for example for the notion of `butterfly minors' in digraphs, has since become a fruitful quest which has given rise to many new discoveries and tools for the structural study of graphs.

A probable application of the ideas from the Graph Minors Project was proposed by Norin \cite{norine2005matching,thomas2006survey} in the form of a potential road to the solution of the \textsc{Pfaffian Recognition} problem.
Pfaffian orientations are of significance because they allow for an efficient way to compute the generating function of perfect matchings within a graph; \textsc{Pfaffian Recognition} asks to decide whether a given graph has a Pfaffian orientation.
See \cite{mccuaig2004polya,thomas2006survey} for more in depth discussions on the topic.
While this problem was solved for bipartite graphs by McCuaig et al.\@ \cite{robertson1999permanents,mccuaig2004polya}, the non-bipartite case remains open.
Norin proposed a variant of treewidth designed for the study of graphs with perfect matchings and conjectured the existence of a grid theorem based on so called `matching minors', a specialised version of minors fit for the study of graphs with perfect matchings.
For non-bipartite graphs this `perfect matching width' and a grid as a potential dual object could give rise to a polynomial time recognition algorithm for non-bipartite Pfaffian graphs.
While the non-bipartite case of Norin's conjecture remains open, the bipartite case was solved in \cite{hatzel2019cyclewidth} and has since let to additional insight into structural aspects of the study of `butterfly minors' in digraphs \cite{giannopoulou2021excluding}.

\paragraph{Our contribution}

This paper is part of an ongoing series in which we try to extend the findings of the Graph Minors Project to matching minors in bipartite graphs.
Moreover, we aim to deepen the understanding of the relation between bipartite graphs with perfect matchings and digraphs in order to obtain new tools and insights for the study of butterfly minors.

The first step towards the global structure theorem for $H$-minor free graphs has been replicated in the setting of matching minors in bipartite graphs in \cite{giannopoulou2021excluding} where it was shown that a class of bipartite graphs has bounded perfect matching width if and only if it excludes some bipartite and planar graph with a perfect matching as a matching minor.
The next step according to the roadmap above would be to describe the local structure of areas of large perfect matching width in $H$-matching minor free bipartite graphs with perfect matchings.
Following this plan, in this paper we establish a matching theoretic version of the \textbf{Flat Wall Theorem} for bipartite graphs. This is enabled by the previous paper in this series in which we have found a matching theoretic version of the \textbf{Two Paths Theorem} \cite{giannopoulou2021two}.

\begin{theorem}\label{thm:matchingflatwall}
	Let $r,t\in\N$ be positive integers.
	There exist functions $\ApexBound\colon\N\rightarrow\N$ and $\WallBound\colon\N\times\N\rightarrow\N$ such that for every bipartite graph $B$ with a perfect matching $M$ the following is true:
	If $W$ is an $M$-conformal \hyperref[def:matchingwall]{matching $\Fkt{\WallBound}{t,r}$-wall} in $B$ such that $M\cap\E{W}$ is the canonical matching of $W$, then either
	\begin{enumerate}
		\item $B$ has a $K_{t,t}$-\hyperref[def:matchingminor]{matching minor} grasped by $W$, or
		
		\item there exist an $M$-conformal set $A\subseteq\V{B}$ with $\Abs{A}\leq\Fkt{\ApexBound}{t}$ and an $M$-conformal matching $r$-wall $W'\subseteq W-A$ such that $W'$ is \hyperref[def:matchingflat]{$\Perimeter{W'}$-flat} in $B$ with respect to $A$.
	\end{enumerate}
\end{theorem}
We present the necessary definitions, in particular those of \emph{matching minors} and \emph{$\Perimeter{W'}$-flatness} in \Cref{sec:background,sec:flatwall}.

As a consequence of the flat wall theorem, we also provide a new directed version of the \textbf{Flat Wall Theorem} which is incomparable with the directed flat wall theorem of Giannopoulou et al.\@ \cite{giannopoulou2020directed} and we use these insights to obtain a new duality theorem for (undirected) treewidth based on matching minors.
All results in this paper are constrictive in the sense that they lead to efficient algorithms for the computation of the respective flat walls and their apex sets.

\subsection{Organisation}\label{sec:organisation}

In the remainder of this introduction we give a more extensive explanation and overview on the \textbf{Flat Wall Theorem} and highlight some of its main features.
We also introduce the directed flat wall theorem of Giannopoulou et al.\@ \cite{giannopoulou2020directed} and explain its merits and drawbacks to get a better point of reference for the comparison with our own results.

In \cref{sec:background} we give an in depth introduction to the relevant parts of matching theory to make this article as self contained as possible and provide further insight into the solution of the bipartite case of \textsc{Pfaffian Recognition} as these findings make up the corner stones of our notion of flatness.

\Cref{sec:flatwall} contains the additional definitions necessary for \Cref{thm:matchingflatwall} and the resulting weak structure theorem for bipartite graphs that do not contain $K_{t,t}$ as a matching minor.
These findings are then carried over into the setting of digraphs in \cref{sec:consequences}, Here we explain the relation between digraphs and bipartite graphs with perfect matchings and give statements and proofs of some immediate consequences of our main theorem.

The remaining sections are dedicated to the proof of the matching theoretic flat wall theorem for bipartite graphs.

\subsection{Preliminaries}\label{sec:preliminaries}

All graphs and digraphs in this article are considered simple, that is we do now allow multiple edges or loops and wherever such objects would arise from contraction, we identify multiple edges and remove loops.
For general notation not introduced in this paper we stick to the book on graph theory by Diestel \cite{diestel2012graph} while for digraph theory we recommend \cite{bang2018classes}.

This paper is relatively long and contains several lesser known concepts and highly technical definitions.
To increase readability we present key definitions in their own block environment and provide hyperlinks to the definitions wherever those are central to the understanding of statements and proofs.

Since the majority of our findings is focussed on bipartite graphs we fix the following convention.
Wherever possible we use $B$ as the standard name for a bipartite graph and $G$ for arbitrary graphs if not stated explicitly otherwise.
Moreover, we assume every bipartite graph to come with a bipartition into the \emph{colour classes} $V_1$ and $V_2$, where in our figures $V_1$ is represented by black vertices and the vertices in $V_2$ are depicted white.
In case ambiguity arises we either treat $V_i$ as the placeholder for all possible vertices of colour $i\in[1,2]$, or we write $\Vi{i}{B}$ to specify which graph we are talking about.
Since we also work a lot on digraphs, we will use $D$ as our standard name for digraphs to make these three different cases more distinguishable.

For integers $i,j\in\Z$ we use the notation $[i,j]$ for the set $\Set{i,i+1,\dots,j}$, where $[i,j]=\emptyset$ if $i>j$.

If $X$ and $Y$ are two finite sets, we denote the \emph{symmetric difference} by $X\Delta Y\coloneqq \Brace{X\setminus Y}\cup \Brace{Y\setminus X}$.

Given a graph $G$ we denote by $\Bidirected{G}$ its \emph{bidirection}, that is the digraph with vertex set $\V{G}$ and edge set $\CondSet{(u,v),(v,u)}{uv\in\E{G}}$.

If $\mathcal{F}$ is a family of (di)graphs we denote by $\V{\mathcal{F}}$ the set $\bigcup_{H\in\mathcal{F}}\V{F}$ and by $\E{\mathcal{F}}$ the set $\bigcup_{H\in\mathcal{F}}\E{F}$.

\begin{definition}[(Directed) Separation]\label{def:separation}
	Let $D$ be a digraph.
	A \emph{directed separation} is a tuple $(X,Y)$ such that $X\cup Y=\V{D}$ and there is no directed edge with tail in $Y\setminus X$ and head in $X\setminus Y$.
	The set $X\cap Y$ is called the \emph{separator}and the \emph{order} of $(X,Y)$ is $\Abs{X\cap Y}$.
	A tuple $(X,Y)$ is called a \emph{separation}, or \emph{undirected separation} if we want to emphasise this fact, if $(X,Y)$ and $(Y,X)$ both are directed separations.
	If $G$ is a graph, then $(X,Y)$ is a separation if it is a separation in $\Bidirected{G}$.
\end{definition}

\subsection{Graph Minor Theory and the Flat Wall}

The combination of clique sums and graphs of bounded size captures exactly the classes of $H$-minor free graphs which are of bounded treewidth.
From the \textbf{Grid Theorem} \cite{robertson1986graph} it follows that this class already contains all $H$-minor free graph classes for which $H$ is planar.
The structure of $H$-minor free graphs where $H$ is non-planar is substantially more complicated.
To better illustrate some of the problems one encounters when excluding a non-planar graph $H$ let us digress a bit.

All graphs and digraphs in this paper are considered simple, so they will not contain loops or parallel edges if not explicitly stated to do so.

\paragraph{The Two Paths Theorem}

Let $C$ be a cycle in a graph $G$.
A \emph{$C$-cross} in $G$ is a pair of disjoint paths $P_1$ and $P_2$ with endpoints $s_1,t_1$, $s_2,t_2$ respectively such that $s_1,s_2,t_1,t_2$ appear on $C$ in the order listed and $P_1$, $P_2$ are otherwise disjoint from $C$.

Let $s_1,s_2,t_1,t_2\in\V{G}$ be any four distinct vertices in a graph $G$.
The \textsc{Two Disjoint Paths Problem} asks whether there exist two disjoint paths $P_1$ and $P_2$ in $G$ such that for each $i\in[1,2]$ $P_i$ has endpoints $s_i$ and $t_i$.
We can reduce the \textsc{Two Disjoint Paths Problem} to the question whether there exists a $C$-cross over the $4$-cycle $C\coloneqq (s_1,s_2,t_1,t_2)$.
Please note that we may add potentially missing edges of $C$ to $G$ without any impact on on the existence of the two disjoint paths since no solution to the \textsc{Two Disjoint Paths Problem} may use any of the edges of $C$.
It turns out that the existence of a $C$-cross in general is linked to two special cases of the list of ingredients from above: \emph{clique sums} (of order at most three) and \emph{graphs on surfaces} (in particular planar graphs).

Let $X\subseteq\V{G}$ be a set of vertices of a graph $G$ and let $(A,B)$ be a \hyperref[def:separation]{separation} of order at most three in $G$ with $X\subseteq A$ such that there exists some $v\in B$ for which $A\cap B$ is a minimum separator between $v$ and $X$.
Let $H'$ be the graph obtained from $\InducedSubgraph{G}{A}$ by joining all pairs of vertices of $A\cap B$ with an edge.
Then $H'$ is called an \emph{elementary $X$-reduction} of $G$.
A graph $H$ is an \emph{X-reduction} of $G$ if it can be obtained by a sequence $(H_1,\dots,H_{\ell})$, $\ell\geq 1$, such that $H_1=G$, $H_{\ell}=H$ and for all $i\in[2,\ell]$ we have that $H_i$ is an elementary $X$-reduction of $H_{i-1}$.

It follows that taking $C$-reductions for some cycle $C$ in $G$ preserves the existence of a $C$-cross.
The key observation now is, in case no further $C$-reduction is possible in $H$, there does not exist a $C$-cross in $H$ if $H$ can be drawn into the plane without crossing edges such that $C$ bounds a face.
Indeed, this is the only obstruction for the existence of a $C$-cross in a $C$-reduction minimal graph as captured by the so called \textbf{Two Paths Theorem}.
The theorem has been obtained in many different forms with various techniques by a plethora of authors over the years.

\begin{theorem}[Two Paths Theorem, \cite{jung1970verallgemeinerung,seymour1980disjoint,shiloach1980polynomial,thomassen19802,robertson1990graph2}]\label{thm:twopaths}
	A cycle $C$ in a graph $G$ has \textbf{no} $C$-cross in $G$ if and only if there exists a $C$-reduction $H$ of $G$ which is planar and in which $C$ bounds a face.
\end{theorem}

\paragraph{The Flat Wall Theorem}

It is possible to slightly generalise the property of cycles $C$ without a $C$-cross to connected planar subgraphs $J$ as follows.
Let $G$ be a graph and $A\subseteq\V{G}$ be a set of vertices, called the \emph{apex set}.
Moreover, let $J\subseteq G-A$ be a subdivision of a $3$-connected planar graph and let $C_J$ be the outer face of $J$.
We say that $J$ is \emph{$A$-flat} in $G$ if there exists a \hyperref[def:separation]{separation} $(X,Y)$ in $G-A$ with $\V{C_J}=X\cap Y$ such that $\V{J}\subseteq X$ and $\InducedSubgraph{G}{X}$ has a $J$-reduction $H$ which is planar and in which $C_J$ bounds a face.

Let $k,t\in\N$ be positive integers.
The \emph{$k\times t$-grid} is the graph with vertex set $\CondSet{v_{i,j}}{i\in[1,k],j\in[1,t]}$, and edge set
$$\CondSet{v_{i_1,j_1}v_{i_2,j_2}}{i_1=i_2\text{ and }\Abs{j_1-j_2}=1\text{, or }j_1=j_2\text{ and }\Abs{i_1-i_2}=1}.$$
An \emph{elementary $k$-wall} is obtained from the $k\times 2k$-grid by deleting all edges with endpoints $v_{2i-1,2j-1}$ and $v_{2i-1,2j}$ for all $i\in[1,\Floor{\frac{k}{2}}]$ and $j\in[1,k]$, and all edges with endpoints $v_{2i,2j}$ and $v_{2i,2j+1}$ for all $i\in[1,\Floor{\frac{k}{2}}]$ and $j\in[1,k]$, and then deleting the two resulting vertices of degree one.
A subdivision of an elementary $k$-wall is called a \emph{$k$-wall}.

The following is another form of the \textbf{Grid Theorem} which uses the fact that for any subcubic graph $H$ a graph $G$ contains $H$ as a minor if and only if it contains a subdivision of $H$ as a subgraph.

\begin{theorem}[Wall Theorem, \cite{robertson1986graph}]\label{cor:undirwall}
	There exists a function $\UndirectedGrid\colon\N\rightarrow\N$ such that for every $k\in\N$ and every graph $G$ we have $\tw{G}\leq\Fkt{\UndirectedGrid}{2k}$, or $G$ contains a $k$-wall as a subgraph.
\end{theorem}

Notice that, if $G$ has a $K_t$-minor, there exist pairwise disjoint sets $X_1,\dots,X_t\subseteq\V{G}$ such that $\InducedSubgraph{G}{X_t}$ is connected and for every pair $X_i$, $X_j$, $i\neq j$, there is an edge $u_iu_j\in\E{G}$ with $u_i\in X_i$ and $u_j\in X_j$.
We say that the $X_1,\dots, X_t$ form a \emph{model} of $K_t$ in $G$.
Let $W$ be some $k$-wall for $k\geq t$ with horizontal paths $P_1,P_2,\dots,P_k$, and vertical paths $Q_1,Q_2,\dots,Q_k$.
We say that a model of a $K_t$-minor is \emph{grasped} by $W$ if for every $X_h$ there are distinct indices $i_1,\dots,i_t\in[1,k]$, and $j_1,\dots,j_t\in[1,k]$ such that $\V{P_{i_{\ell}}\cap Q_{j_{\ell}}}\subseteq X_h$ for all $\ell\in[1,t]$.

The following theorem was first obtained in \cite{robertson1995graph}, but there it was stated in a slightly weaker form and the bounds were not given explicitly. 
So we credit Kawarabayashi et al.\@ \cite{kawarabayashi2018new} for simplifying the proof and providing explicit bounds.

\begin{theorem}[Flat Wall Theorem, \cite{robertson1995graph,kawarabayashi2018new}]\label{thm:undirectedflatwall}
	Let $r,t\in\N$ be positive integers, $R=49152t^{24}\Brace{40t^2+r}$, and $G$ be a graph.
	Then the following is true:
	If $W$ is an $R$-wall in $G$, then either
	\begin{enumerate}
		\item $G$ has a $K_t$-minor grasped by $W$, or
		
		\item there exists a set $A\subseteq\V{G}$ with $\Abs{A}\leq 12288t^{24}$ and an $r$-wall $W'\subseteq W-A$ such that $W'$ is $A$-flat in $G$.
	\end{enumerate}
\end{theorem}

Indeed, one can take the Flat Wall Theorem and obtain an almost-characterisation of all graphs $G$ that exclude $K_t$ as a minor.

\begin{theorem}[Weak Structure Theorem\footnote{Traditionally the Flat Wall Theorem and the Weak Structure Theorem are used synonymous, but since there are two slightly different statements it might make sense to differentiate between the two.}, \cite{kawarabayashi2018new}]\label{thm:undirectedwaekstructure}
	Let $r,t\in\N$ be positive integers, $R=49152t^{24}\Brace{40t^2+r}$, and $G$ be a graph.
	If $G$ has no $K_t$-minor, then for every $R$-wall $W$ in $G$ there exists a set $A\subseteq\V{G}$ with $\Abs{A}\leq 12288t^{24}$ and an $r$-wall $W'\subseteq W-A$ such that $W'$ is $A$-flat in $G$.
	Conversely, if $t\geq 2$, $r\geq 80t^{12}$, and for every $R$-wall $W$ in $G$ there is a set $A\subseteq\V{G}$ of size at most $12288t^{24}$ and an $r$-wall $W'\subseteq W-A$ which is $A$-flat in $G$, then $G$ has no $K_{t'}$-minor, where $t'=2R^2$.
\end{theorem}

Recall the list of ingredients from the start of the introduction.
In the \textbf{Flat Wall Theorem} clique sums and graphs of bounded size are present implicitly since it only treats the case of $H$-minor free graphs with large enough treewidth.
Moreover, clique sums also appear in the form of reductions in the definition of $A$-flatness where we use small ordered \hyperref[def:separation]{separation} to cut away non-planar parts which still exist after deleting $A$, but cannot contribute to any substantial connectivity over the wall.
Additionally, we are given a first glimpse of the `graphs on surfaces' part by the flatness of the wall and apex vertices make a first appearance.

\subsection{The Strange Case of Digraphs}\label{sec:digraphintro}

A key feature of the \textbf{Flat Wall Theorem} is the definition of flatness itself.
The fact that the existence of two disjoint but crossing paths is directly linked to a topological property plays a central role in this theorem.
In the case of digraphs one encounters several problems when trying to follow along the three steps from above in order to create a structural theory of minors in digraphs.
In this subsection we introduce some of the main concepts currently used in structural digraph theory and discuss why they appear to be unavoidable while still being somewhat flawed when compared to their undirected counterparts.

To this end let us start by introducing some basic concepts.
The Graph Minors Project was started to investigate the structure of $H$-minor free graphs.
If we want to achieve something similar for digraphs we first need a notion of minors that works well with directed graphs.
However, simply allowing to contract any directed edge might lead to the introduction of new directed cycles where no directed cycles might have existed before.
In particular, it is possible to give rise to non-trivial strong components in digraphs without any directed cycles.
While this is not necessarily a bad thing depending on the application, it is possible that the introduction of new directed cycles might change the reachability of vertices by also introducing directed paths that did not exist before the contraction.
Such a behaviour makes it hart to control the structure using the means of digraphs and does not lead to a structural theory that is different from the structure of the underlying undirected graph.
Thus more restricted notions for directed minors were pursued.

The first notion that enabled structural characterisations of classes of digraphs, although not under the name they would eventually become known by, was the idea of butterfly minors \cite{seymour1987characterization}.

Let $D$ be a digraph and $\Brace{u,v}\in\E{D}$.
The edge $\Brace{u,v}$ is \emph{butterfly contractible} if $\OutNeighbours{D}{u}=\Set{v}$, or $\InNeighbours{D}{v}=\Set{u}$.
As in undirected graphs we identify parallel edges and remove loops after contracting a butterfly contractible edge.
A digraph $H$ is a \emph{butterfly minor} of $D$ if it can be obtained from $D$ by a sequence of edge-deletions, vertex-deletions, and butterfly contractions.

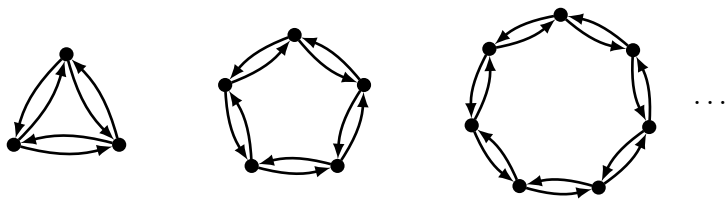
\begin{figure}[h!]
	\begin{center}
		\begin{tikzpicture}
			
			\node (D1) [v:ghost] {};
			\node (D2) [v:ghost,position=0:30mm from D1] {};
			\node (D3) [v:ghost,position=0:35mm from D2] {};
			\node (D4) [v:ghost,position=0:20mm from D3] {$\dots$};

			\node (C1) [v:ghost,position=0:0mm from D1] {
				
				\begin{tikzpicture}[scale=0.8]
					
					\pgfdeclarelayer{background}
					\pgfdeclarelayer{foreground}
					
					\pgfsetlayers{background,main,foreground}
					
					\begin{pgfonlayer}{main}
						
						\node (C) [] {};

						

						
						
						\node (v1) [v:main,position=90:10mm from C] {};
						\node (v2) [v:main,position=210:10mm from C] {};
						\node (v3) [v:main,position=330:10mm from C] {};
						
						
						
						

						

						
						
						\draw (v1) [e:main,->,bend right=15] to (v2);
						\draw (v2) [e:main,->,bend right=15] to (v3);
						\draw (v3) [e:main,->,bend right=15] to (v1);
						
						\draw (v1) [e:main,->,bend right=15] to (v3);
						\draw (v2) [e:main,->,bend right=15] to (v1);
						\draw (v3) [e:main,->,bend right=15] to (v2);

						

						
						
					\end{pgfonlayer}
					

					\begin{pgfonlayer}{background}
						
					\end{pgfonlayer}	
					
					\begin{pgfonlayer}{foreground}

					\end{pgfonlayer}
				\end{tikzpicture}
				
			};
			
			\node (C2) [v:ghost,position=0:0mm from D2] {
				
				\begin{tikzpicture}[scale=0.8]
					
					\pgfdeclarelayer{background}
					\pgfdeclarelayer{foreground}
					
					\pgfsetlayers{background,main,foreground}
					
					\begin{pgfonlayer}{main}
						
						\node (C) [] {};

						

						
						
						\node (v1) [v:main,position=90:12mm from C] {};
						\node (v2) [v:main,position=162:12mm from C] {};
						\node (v3) [v:main,position=234:12mm from C] {};
						\node (v4) [v:main,position=306:12mm from C] {};
						\node (v5) [v:main,position=18:12mm from C] {};
						
						
						
						

						

						
						
						\draw (v1) [e:main,->,bend right=15] to (v2);
						\draw (v2) [e:main,->,bend right=15] to (v3);
						\draw (v3) [e:main,->,bend right=15] to (v4);
						\draw (v4) [e:main,->,bend right=15] to (v5);
						\draw (v5) [e:main,->,bend right=15] to (v1);
						
						\draw (v1) [e:main,->,bend right=15] to (v5);
						\draw (v2) [e:main,->,bend right=15] to (v1);
						\draw (v3) [e:main,->,bend right=15] to (v2);
						\draw (v4) [e:main,->,bend right=15] to (v3);
						\draw (v5) [e:main,->,bend right=15] to (v4);
						
						

						
						
					\end{pgfonlayer}
					

					\begin{pgfonlayer}{background}
						
					\end{pgfonlayer}	
					
					\begin{pgfonlayer}{foreground}

					\end{pgfonlayer}
				\end{tikzpicture}
				
			};
			
			\node (C3) [v:ghost,position=0:0mm from D3] {
				
				\begin{tikzpicture}[scale=0.8]
					
					\pgfdeclarelayer{background}
					\pgfdeclarelayer{foreground}
					
					\pgfsetlayers{background,main,foreground}
					
					\begin{pgfonlayer}{main}
						
						\node (C) [] {};

						

						
						
						\node (v1) [v:main,position=90:15mm from C] {};
						\node (v2) [v:main,position=141.4:15mm from C] {};
						\node (v3) [v:main,position=192.8:15mm from C] {};
						\node (v4) [v:main,position=243.2:15mm from C] {};
						\node (v5) [v:main,position=294.6:15mm from C] {};
						\node (v6) [v:main,position=346.1:15mm from C] {};
						\node (v7) [v:main,position=37.5:15mm from C] {};
						
						
						
						

						

						
						
						\draw (v1) [e:main,->,bend right=15] to (v2);
						\draw (v2) [e:main,->,bend right=15] to (v3);
						\draw (v3) [e:main,->,bend right=15] to (v4);
						\draw (v4) [e:main,->,bend right=15] to (v5);
						\draw (v5) [e:main,->,bend right=15] to (v6);
						\draw (v6) [e:main,->,bend right=15] to (v7);
						\draw (v7) [e:main,->,bend right=15] to (v1);
						
						\draw (v2) [e:main,->,bend right=15] to (v1);
						\draw (v3) [e:main,->,bend right=15] to (v2);
						\draw (v4) [e:main,->,bend right=15] to (v3);
						\draw (v5) [e:main,->,bend right=15] to (v4);
						\draw (v6) [e:main,->,bend right=15] to (v5);
						\draw (v7) [e:main,->,bend right=15] to (v6);
						\draw (v1) [e:main,->,bend right=15] to (v7);

						

						
						
					\end{pgfonlayer}
					

					\begin{pgfonlayer}{background}
						
					\end{pgfonlayer}	
					
					\begin{pgfonlayer}{foreground}

					\end{pgfonlayer}
				\end{tikzpicture}
				
			};
			
		\end{tikzpicture}
	\end{center}
	\caption{The anti-chain $\Antichain{\Bidirected{C_3}}$ of odd bicycles.}
	\label{fig:oddbicycles}
\end{figure}

With this restricted setting of legal contractions it was ensured that no new directed cycles would be created.
However, there exist many infinite anti-chains for the butterfly minor relation such as the family of \emph{odd bicycles} from \cite{seymour1987characterization} which can be found in \cref{fig:oddbicycles}.

To circumvent this issue different generalisations of butterfly minors were introduced such as `directed minors' \cite{kintali2017forbidden} which in addition to the butterfly minors operations also allow for the contraction of whole directed cycles at once, and `strong minors' which replace the directed-cycle rule of directed minors by allowing to contract whole strongly connected subdigraphs at once\footnote{Observe that there is actually no difference between directed minors and strong minors.} \cite{kim2015tournament}.
Strong minors were shown to yield a well-quasi ordering for tournaments \cite{kim2015tournament}, but beyond some limited results these two notions did not receive a lot of attention.
A possible reason for this can be found in the introduction of directed treewidth and the results around it.

Let $D$ be a digraph and $X,Y\subseteq\V{D}$.
A directed walk $W$ is a \emph{directed $X$-walk} if it starts and ends in $X$, and contains a vertex of $\V{D-X}$.
We say that $Y$ \emph{strongly guards} $X$ if every directed $X$-walk in $D$ contains a vertex of $Y$.
The set $Y$ \emph{weakly guards} $X$ if every directed $X$-$\V{D-X}$-path contains a vertex of $Y$, 

An \emph{arborescence} is a digraph $\vec{T}$ obtained from a tree $T$ by selecting a \emph{root} $r\in\V{T}$ and orienting all edges of $T$ away from $r$.
If $e$ is a directed edge and $v$ is an endpoint of $e$ we write $v\sim e$.

\begin{definition}[Directed Treewidth]
	Let $D$ be a digraph.
	A \emph{directed tree decomposition} for $D$ is a tuple $\Brace{T,\beta,\gamma}$ where $T$ is an arborescence, $\beta\colon\Fkt{V}{T}\rightarrow 2^{V(D)}$ is a function that partitions $\Fkt{V}{D}$, into sets called the \emph{bags}\footnote{This means $\CondSet{\Fkt{\beta}{t}}{t\in\V{T}}$ is a partition of $\V{D}$ into non-empty sets.}, and $\gamma\colon\Fkt{E}{T}\rightarrow 2^{V(D)}$ is a function, giving us sets called the \emph{guards}, satisfying the following requirement:
	\begin{enumerate}
		\item[] For every $\Brace{d,t}\in\Fkt{E}{T}$, $\Fkt{\gamma}{d,t}$ strongly guards $\Fkt{\beta}{T_t}\coloneqq\bigcup_{t'\in V(T_t)}\Fkt{\beta}{t'}$.
	\end{enumerate}
	Here $T_t$ denotes the subarboresence of $T$ with root $t$.
	For every $t\in\Fkt{V}{T}$ let $\Fkt{\Gamma}{t}\coloneqq\Fkt{\beta}{t}\cup\bigcup_{t\sim e}\Fkt{\gamma}{e}$.
	The \emph{width} of $\Brace{T,\beta,\gamma}$ is defined as
	\begin{align*}
		\Width{T,\beta,\gamma}\coloneqq\max_{t\in V(T)}\Abs{\Fkt{\Gamma}{t}}-1.
	\end{align*}
	The \emph{directed treewidth} of $D$, denoted by $\dtw{D}$, is the minimum width over all directed tree decompositions for $D$.
\end{definition}

In the paper by Johnson et al.\@ \cite{johnson2001directed} it was conjectured that, similar to the undirected case, there would exist a grid like structure that acted as a universal obstruction to small directed treewidth.
Moreover, they conjectured this grid to exist as a butterfly minor, so no further additions to the minor relation would be necessary.
The directed grid theorem was proven almost $15$ years later by Kawarabayashi and Kreutzer.

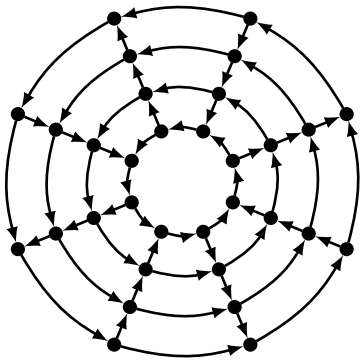
\begin{figure}[!h]
	\centering
	\begin{tikzpicture}[scale=0.9]
		\pgfdeclarelayer{background}
		\pgfdeclarelayer{foreground}
		\pgfsetlayers{background,main,foreground}
		
		\node (mid) [v:ghost] {};
		
		\node (u11) [v:main,position=22.5:8mm from mid] {};
		\node (u21) [v:main,position=67.5:8mm from mid] {};
		\node (u31) [v:main,position=112.5:8mm from mid] {};
		\node (u41) [v:main,position=157.5:8mm from mid] {};
		\node (u51) [v:main,position=202.5:8mm from mid] {};
		\node (u61) [v:main,position=247.5:8mm from mid] {};
		\node (u71) [v:main,position=292.5:8mm from mid] {};
		\node (u81) [v:main,position=337.5:8mm from mid] {};
		
		\node (u12) [v:main,position=22.5:14mm from mid] {};
		\node (u22) [v:main,position=67.5:14mm from mid] {};
		\node (u32) [v:main,position=112.5:14mm from mid] {};
		\node (u42) [v:main,position=157.5:14mm from mid] {};
		\node (u52) [v:main,position=202.5:14mm from mid] {};
		\node (u62) [v:main,position=247.5:14mm from mid] {};
		\node (u72) [v:main,position=292.5:14mm from mid] {};
		\node (u82) [v:main,position=337.5:14mm from mid] {};
		
		\node (u13) [v:main,position=22.5:20mm from mid] {};
		\node (u23) [v:main,position=67.5:20mm from mid] {};
		\node (u33) [v:main,position=112.5:20mm from mid] {};
		\node (u43) [v:main,position=157.5:20mm from mid] {};
		\node (u53) [v:main,position=202.5:20mm from mid] {};
		\node (u63) [v:main,position=247.5:20mm from mid] {};
		\node (u73) [v:main,position=292.5:20mm from mid] {};
		\node (u83) [v:main,position=337.5:20mm from mid] {};
		
		\node (u14) [v:main,position=22.5:26mm from mid] {};
		\node (u24) [v:main,position=67.5:26mm from mid] {};
		\node (u34) [v:main,position=112.5:26mm from mid] {};
		\node (u44) [v:main,position=157.5:26mm from mid] {};
		\node (u54) [v:main,position=202.5:26mm from mid] {};
		\node (u64) [v:main,position=247.5:26mm from mid] {};
		\node (u74) [v:main,position=292.5:26mm from mid] {};
		\node (u84) [v:main,position=337.5:26mm from mid] {};
		
		\begin{pgfonlayer}{background}
			
			\draw[e:main,bend right=15,->] (u11) to (u21);
			\draw[e:main,bend right=15,->] (u21) to (u31);
			\draw[e:main,bend right=15,->] (u31) to (u41);
			\draw[e:main,bend right=15,->] (u41) to (u51);
			\draw[e:main,bend right=15,->] (u51) to (u61);
			\draw[e:main,bend right=15,->] (u61) to (u71);
			\draw[e:main,bend right=15,->] (u71) to (u81);
			\draw[e:main,bend right=15,->] (u81) to (u11);
			
			\draw[e:main,bend right=15,->] (u12) to (u22);
			\draw[e:main,bend right=15,->] (u22) to (u32);
			\draw[e:main,bend right=15,->] (u32) to (u42);
			\draw[e:main,bend right=15,->] (u42) to (u52);
			\draw[e:main,bend right=15,->] (u52) to (u62);
			\draw[e:main,bend right=15,->] (u62) to (u72);
			\draw[e:main,bend right=15,->] (u72) to (u82);
			\draw[e:main,bend right=15,->] (u82) to (u12);
			
			\draw[e:main,bend right=15,->] (u13) to (u23);
			\draw[e:main,bend right=15,->] (u23) to (u33);
			\draw[e:main,bend right=15,->] (u33) to (u43);
			\draw[e:main,bend right=15,->] (u43) to (u53);
			\draw[e:main,bend right=15,->] (u53) to (u63);
			\draw[e:main,bend right=15,->] (u63) to (u73);
			\draw[e:main,bend right=15,->] (u73) to (u83);
			\draw[e:main,bend right=15,->] (u83) to (u13);
			
			\draw[e:main,bend right=15,->] (u14) to (u24);
			\draw[e:main,bend right=15,->] (u24) to (u34);
			\draw[e:main,bend right=15,->] (u34) to (u44);
			\draw[e:main,bend right=15,->] (u44) to (u54);
			\draw[e:main,bend right=15,->] (u54) to (u64);
			\draw[e:main,bend right=15,->] (u64) to (u74);
			\draw[e:main,bend right=15,->] (u74) to (u84);
			\draw[e:main,bend right=15,->] (u84) to (u14);
			
			\draw[e:main,->] (u11) to (u12);
			\draw[e:main,->] (u12) to (u13);
			\draw[e:main,->] (u13) to (u14);
			
			\draw[e:main,->] (u24) to (u23);
			\draw[e:main,->] (u23) to (u22);
			\draw[e:main,->] (u22) to (u21);
			
			\draw[e:main,->] (u31) to (u32);
			\draw[e:main,->] (u32) to (u33);
			\draw[e:main,->] (u33) to (u34);
			
			\draw[e:main,->] (u44) to (u43);
			\draw[e:main,->] (u43) to (u42);
			\draw[e:main,->] (u42) to (u41);
			
			\draw[e:main,->] (u51) to (u52);
			\draw[e:main,->] (u52) to (u53);
			\draw[e:main,->] (u53) to (u54);
			
			\draw[e:main,->] (u64) to (u63);
			\draw[e:main,->] (u63) to (u62);
			\draw[e:main,->] (u62) to (u61);
			
			\draw[e:main,->] (u71) to (u72);
			\draw[e:main,->] (u72) to (u73);
			\draw[e:main,->] (u73) to (u74);
			
			\draw[e:main,->] (u84) to (u83);
			\draw[e:main,->] (u83) to (u82);
			\draw[e:main,->] (u82) to (u81);
			
		\end{pgfonlayer}
	\end{tikzpicture}
	\caption{The cylindrical grid of order $4$.}
	\label{fig:directedgrid1}
\end{figure}

Let $k\in\N$ be a positive integer.
The \emph{cylindrical grid of order $k$} is the digraph obtained from the cycles $C_1,\dots C_k$, with $C_i=\Brace{v_0^i,e^i_0,v_1^i,e_1^i,\dots,e_{2k-3}^i,v_{2k-2}^i,e_{ek-2}^i,v_{2k-1}^i,e_{2k-1}^i,v_0^i}$
for each $i\in[1,k]$, by adding the directed paths $P_i = v_{2i}^1v_{2i}^2\dots v_{2i}^{k-1}v_{2i}^k$ and $Q_i = v_{2i+1}^1v_{2i+1}^2\dots v_{2i+1}^{k-1}v_{2i+1}^k$ for every $i\in[0,k-1]$.
See \cref{fig:directedgrid1} for an illustration.

\begin{theorem}[Directed Grid Theorem, \cite{kawarabayashi2015directed}]\label{thm:directedgrid}
	There exists a function $\DirectedGrid\colon\N\rightarrow \N$ such that for every $k\in\N$ and every digraph $D$ we have $\dtw{D}\leq\Fkt{\DirectedGrid}{k}$, or $D$ contains the cylindrical grid of order $k$ as a butterfly minor.
\end{theorem}

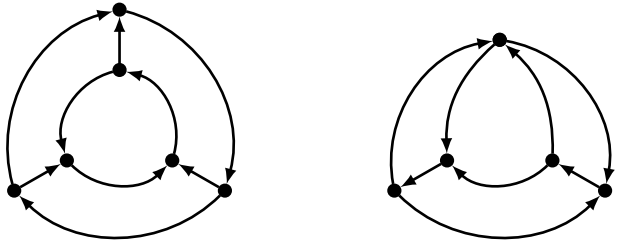
\begin{figure}[!h]
	\centering
	\begin{tikzpicture}
		\pgfdeclarelayer{background}
		\pgfdeclarelayer{foreground}
		\pgfsetlayers{background,main,foreground}
		
		\node (mid) [v:ghost] {};
		
		\node (v1) [v:main,position=90:8mm from mid] {};
		\node (v2) [v:main,position=210:8mm from mid] {};
		\node (v3) [v:main,position=330:8mm from mid] {};
		
		\node (v4) [v:main,position=90:16mm from mid] {};
		\node (v5) [v:main,position=210:16mm from mid] {};
		\node (v6) [v:main,position=330:16mm from mid] {};
		
		\draw [e:main,->,bend right=45] (v1) to (v2);
		\draw [e:main,->,bend right=45] (v2) to (v3);
		\draw [e:main,->,bend right=45] (v3) to (v1);
		
		\draw [e:main,->,bend left=45] (v4) to (v6);
		\draw [e:main,->,bend left=45] (v6) to (v5);
		\draw [e:main,->,bend left=45] (v5) to (v4);
		
		\draw [e:main,->] (v1) to (v4);
		\draw [e:main,->] (v5) to (v2);
		\draw [e:main,->] (v6) to (v3);
		
		\node (Bmid) [v:ghost,position=0:50mm from mid] {};
		
		\node (Bv1) [v:main,position=90:12mm from Bmid] {};
		\node (Bv2) [v:main,position=210:8mm from Bmid] {};
		\node (Bv3) [v:main,position=330:8mm from Bmid] {};
		
		\node (Bv4) [v:main,position=90:12mm from Bmid] {};
		\node (Bv5) [v:main,position=210:16mm from Bmid] {};
		\node (Bv6) [v:main,position=330:16mm from Bmid] {};
		
		\draw [e:main,->,bend right=25] (Bv1) to (Bv2);
		\draw [e:main,<-,bend right=45] (Bv2) to (Bv3);
		\draw [e:main,->,bend right=25] (Bv3) to (Bv1);
		
		\draw [e:main,->,bend left=45] (Bv4) to (Bv6);
		\draw [e:main,<-,bend left=45] (Bv6) to (Bv5);
		\draw [e:main,->,bend left=45] (Bv5) to (Bv4);
		
		\draw [e:main,<-] (Bv5) to (Bv2);
		\draw [e:main,->] (Bv6) to (Bv3);
		
	\end{tikzpicture}
	\caption{Two strongly connected and planar digraphs that are \textbf{not} a butterfly minor of the cylindrical grid.
	The digraph on the left is strongly planar while the digraph on the right is not.}
	\label{fig:nongridminorplanardigraph}
\end{figure}

An interesting observation about the cylindrical grid is, that of course it itself hast large directed treewidth.
That means any statement for all digraphs of large directed treewidth must also be true for the cylindrical grid itself.
In \cref{fig:nongridminorplanardigraph} we present two strongly connected and planar digraphs both of which are not butterfly minors of any cylindrical grid.

The fact that such digraphs exist can be seen as an indication that obtaining any structural results beyond the directed grid theorem for digraphs that exclude a fixed butterfly minor would be a difficult task.
To date no structural description of $H$-butterfly minor free digraphs where $H$ is strongly connected and planar is known.
One could interpret this as a problem inherited from butterfly minors, but indeed one can observe that the only strong minor contained in the cylindrical grid which has large directed treewidth is the cylindrical grid itself.
So by changing the notion of minors we use, we cannot expect much more.
Another way to resolve this problem would be to search for an alternative for directed treewidth.
But while many such attempts have been made (see \cite{bang2018classes} for an overview), no structural parameter has been found that could fill this role.
Moreover, there are several strong indications, in form of other dualities like brambles and tangles \cite{giannopoulou2020canonical}, that directed treewidth is indeed the best we can do.

\paragraph{The Directed Flat Wall Theorem}

Let $H$ be a strongly connected digraph.
Suppose we are content with the fact that the class of $H$-butterfly minor free digraphs has bounded directed treewidth if and only if $H$ is a butterfly minor of a cylindrical grid.
Then we would still be left with the question, what if $H$ is not a butterfly minor of some cylindrical grid?
We know that $H$ might still be planar, but it might as well not be.
Note that \cref{thm:undirectedflatwall} only considers the complete graph rather than general graphs $H'$.
Still, since $H'$ is a subgraph, and hence a minor, of $K_{\Abs{\V{H'}}}$, this is enough to roughly capture the structure for all $H'$-minor free graphs.
Let us now take a similar route.

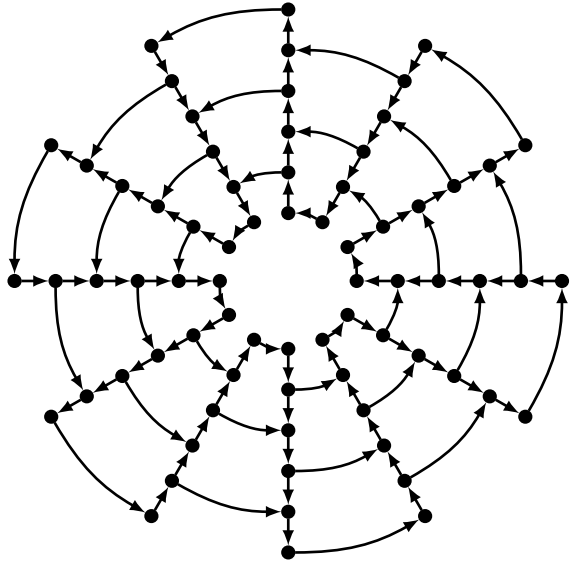
\begin{figure}[!h]
	\centering
	\begin{tikzpicture}[scale=0.9]
		\pgfdeclarelayer{background}
		\pgfdeclarelayer{foreground}
		\pgfsetlayers{background,main,foreground}
		
		\node (mid) [v:ghost] {};
		
		\node (u11) [v:main,position=30:10mm from mid] {};
		\node (u21) [v:main,position=60:10mm from mid] {};
		\node (u31) [v:main,position=90:10mm from mid] {};
		\node (u41) [v:main,position=120:10mm from mid] {};
		\node (u51) [v:main,position=150:10mm from mid] {};
		\node (u61) [v:main,position=180:10mm from mid] {};
		\node (u71) [v:main,position=210:10mm from mid] {};
		\node (u81) [v:main,position=240:10mm from mid] {};
		\node (u91) [v:main,position=270:10mm from mid] {};
		\node (u101) [v:main,position=300:10mm from mid] {};
		\node (u111) [v:main,position=330:10mm from mid] {};
		\node (u121) [v:main,position=0:10mm from mid] {};
		
		\node (u12) [v:main,position=30:16mm from mid] {};
		\node (u22) [v:main,position=60:16mm from mid] {};
		\node (u32) [v:main,position=90:16mm from mid] {};
		\node (u42) [v:main,position=120:16mm from mid] {};
		\node (u52) [v:main,position=150:16mm from mid] {};
		\node (u62) [v:main,position=180:16mm from mid] {};
		\node (u72) [v:main,position=210:16mm from mid] {};
		\node (u82) [v:main,position=240:16mm from mid] {};
		\node (u92) [v:main,position=270:16mm from mid] {};
		\node (u102) [v:main,position=300:16mm from mid] {};
		\node (u112) [v:main,position=330:16mm from mid] {};
		\node (u122) [v:main,position=0:16mm from mid] {};
		
		\node (u13) [v:main,position=30:22mm from mid] {};
		\node (u23) [v:main,position=60:22mm from mid] {};
		\node (u33) [v:main,position=90:22mm from mid] {};
		\node (u43) [v:main,position=120:22mm from mid] {};
		\node (u53) [v:main,position=150:22mm from mid] {};
		\node (u63) [v:main,position=180:22mm from mid] {};
		\node (u73) [v:main,position=210:22mm from mid] {};
		\node (u83) [v:main,position=240:22mm from mid] {};
		\node (u93) [v:main,position=270:22mm from mid] {};
		\node (u103) [v:main,position=300:22mm from mid] {};
		\node (u113) [v:main,position=330:22mm from mid] {};
		\node (u123) [v:main,position=0:22mm from mid] {};
		
		\node (u14) [v:main,position=30:28mm from mid] {};
		\node (u24) [v:main,position=60:28mm from mid] {};
		\node (u34) [v:main,position=90:28mm from mid] {};
		\node (u44) [v:main,position=120:28mm from mid] {};
		\node (u54) [v:main,position=150:28mm from mid] {};
		\node (u64) [v:main,position=180:28mm from mid] {};
		\node (u74) [v:main,position=210:28mm from mid] {};
		\node (u84) [v:main,position=240:28mm from mid] {};
		\node (u94) [v:main,position=270:28mm from mid] {};
		\node (u104) [v:main,position=300:28mm from mid] {};
		\node (u114) [v:main,position=330:28mm from mid] {};
		\node (u124) [v:main,position=0:28mm from mid] {};
		
		\node (u15) [v:main,position=30:34mm from mid] {};
		\node (u25) [v:main,position=60:34mm from mid] {};
		\node (u35) [v:main,position=90:34mm from mid] {};
		\node (u45) [v:main,position=120:34mm from mid] {};
		\node (u55) [v:main,position=150:34mm from mid] {};
		\node (u65) [v:main,position=180:34mm from mid] {};
		\node (u75) [v:main,position=210:34mm from mid] {};
		\node (u85) [v:main,position=240:34mm from mid] {};
		\node (u95) [v:main,position=270:34mm from mid] {};
		\node (u105) [v:main,position=300:34mm from mid] {};
		\node (u115) [v:main,position=330:34mm from mid] {};
		\node (u125) [v:main,position=0:34mm from mid] {};
		
		\node (u16) [v:main,position=30:40mm from mid] {};
		\node (u26) [v:main,position=60:40mm from mid] {};
		\node (u36) [v:main,position=90:40mm from mid] {};
		\node (u46) [v:main,position=120:40mm from mid] {};
		\node (u56) [v:main,position=150:40mm from mid] {};
		\node (u66) [v:main,position=180:40mm from mid] {};
		\node (u76) [v:main,position=210:40mm from mid] {};
		\node (u86) [v:main,position=240:40mm from mid] {};
		\node (u96) [v:main,position=270:40mm from mid] {};
		\node (u106) [v:main,position=300:40mm from mid] {};
		\node (u116) [v:main,position=330:40mm from mid] {};
		\node (u126) [v:main,position=0:40mm from mid] {};
		
		\begin{pgfonlayer}{background}
			
			\draw[e:main,bend right=15,->] (u21) to (u31);
			\draw[e:main,bend right=15,->] (u41) to (u51);
			\draw[e:main,bend right=15,->] (u61) to (u71);
			\draw[e:main,bend right=15,->] (u81) to (u91);
			\draw[e:main,bend right=15,->] (u101) to (u111);
			\draw[e:main,bend right=15,->] (u121) to (u11);
			
			\draw[e:main,bend right=15,->] (u12) to (u22);
			\draw[e:main,bend right=15,->] (u32) to (u42);
			\draw[e:main,bend right=15,->] (u52) to (u62);
			\draw[e:main,bend right=15,->] (u72) to (u82);
			\draw[e:main,bend right=15,->] (u92) to (u102);
			\draw[e:main,bend right=15,->] (u112) to (u122);
			
			\draw[e:main,bend right=15,->] (u23) to (u33);
			\draw[e:main,bend right=15,->] (u43) to (u53);
			\draw[e:main,bend right=15,->] (u63) to (u73);
			\draw[e:main,bend right=15,->] (u83) to (u93);
			\draw[e:main,bend right=15,->] (u103) to (u113);
			\draw[e:main,bend right=15,->] (u123) to (u13);
			
			\draw[e:main,bend right=15,->] (u14) to (u24);
			\draw[e:main,bend right=15,->] (u34) to (u44);
			\draw[e:main,bend right=15,->] (u54) to (u64);
			\draw[e:main,bend right=15,->] (u74) to (u84);
			\draw[e:main,bend right=15,->] (u94) to (u104);
			\draw[e:main,bend right=15,->] (u114) to (u124);
			
			\draw[e:main,bend right=15,->] (u25) to (u35);
			\draw[e:main,bend right=15,->] (u45) to (u55);
			\draw[e:main,bend right=15,->] (u65) to (u75);
			\draw[e:main,bend right=15,->] (u85) to (u95);
			\draw[e:main,bend right=15,->] (u105) to (u115);
			\draw[e:main,bend right=15,->] (u125) to (u15);
			
			\draw[e:main,bend right=15,->] (u16) to (u26);
			\draw[e:main,bend right=15,->] (u36) to (u46);
			\draw[e:main,bend right=15,->] (u56) to (u66);
			\draw[e:main,bend right=15,->] (u76) to (u86);
			\draw[e:main,bend right=15,->] (u96) to (u106);
			\draw[e:main,bend right=15,->] (u116) to (u126);
			
			\draw[e:main,->] (u11) to (u12);
			\draw[e:main,->] (u12) to (u13);
			\draw[e:main,->] (u13) to (u14);
			\draw[e:main,->] (u14) to (u15);
			\draw[e:main,->] (u15) to (u16);
			
			\draw[e:main,->] (u26) to (u25);
			\draw[e:main,->] (u25) to (u24);
			\draw[e:main,->] (u24) to (u23);
			\draw[e:main,->] (u23) to (u22);
			\draw[e:main,->] (u22) to (u21);
			
			\draw[e:main,->] (u31) to (u32);
			\draw[e:main,->] (u32) to (u33);
			\draw[e:main,->] (u33) to (u34);
			\draw[e:main,->] (u34) to (u35);
			\draw[e:main,->] (u35) to (u36);
			
			\draw[e:main,->] (u46) to (u45);
			\draw[e:main,->] (u45) to (u44);
			\draw[e:main,->] (u44) to (u43);
			\draw[e:main,->] (u43) to (u42);
			\draw[e:main,->] (u42) to (u41);
			
			\draw[e:main,->] (u51) to (u52);
			\draw[e:main,->] (u52) to (u53);
			\draw[e:main,->] (u53) to (u54);
			\draw[e:main,->] (u54) to (u55);
			\draw[e:main,->] (u55) to (u56);
			
			\draw[e:main,->] (u66) to (u65);
			\draw[e:main,->] (u65) to (u64);
			\draw[e:main,->] (u64) to (u63);
			\draw[e:main,->] (u63) to (u62);
			\draw[e:main,->] (u62) to (u61);
			
			\draw[e:main,->] (u71) to (u72);
			\draw[e:main,->] (u72) to (u73);
			\draw[e:main,->] (u73) to (u74);
			\draw[e:main,->] (u74) to (u75);
			\draw[e:main,->] (u75) to (u76);
			
			\draw[e:main,->] (u86) to (u85);
			\draw[e:main,->] (u85) to (u84);
			\draw[e:main,->] (u84) to (u83);
			\draw[e:main,->] (u83) to (u82);
			\draw[e:main,->] (u82) to (u81);
			
			\draw[e:main,->] (u91) to (u92);
			\draw[e:main,->] (u92) to (u93);
			\draw[e:main,->] (u93) to (u94);
			\draw[e:main,->] (u94) to (u95);
			\draw[e:main,->] (u95) to (u96);
			
			\draw[e:main,->] (u106) to (u105);
			\draw[e:main,->] (u105) to (u104);
			\draw[e:main,->] (u104) to (u103);
			\draw[e:main,->] (u103) to (u102);
			\draw[e:main,->] (u102) to (u101);
			
			\draw[e:main,->] (u111) to (u112);
			\draw[e:main,->] (u112) to (u113);
			\draw[e:main,->] (u113) to (u114);
			\draw[e:main,->] (u114) to (u115);
			\draw[e:main,->] (u115) to (u116);
			
			\draw[e:main,->] (u126) to (u125);
			\draw[e:main,->] (u125) to (u124);
			\draw[e:main,->] (u124) to (u123);
			\draw[e:main,->] (u123) to (u122);
			\draw[e:main,->] (u122) to (u121);
			
		\end{pgfonlayer}
	\end{tikzpicture}
	\caption{The elementary directed $3$-wall.}
	\label{fig:directedwall1}
\end{figure}

Let $k\in\N$ be a positive integer.
An \emph{elementary cylindrical $k$-wall} $W$ is the digraph obtained from the cylindrical grid $G$ of order $2k$ by deleting the edges $\Brace{v_{2i}^{2j},v_{2i+1}^{2j}}$ and $\Brace{v_{2i+1}^{2j+1},v_{2i+2}^{2j+1}}$ for every $i\in[0,2k-1]$ and every $j\in[0,k-1]$.
A \emph{cylindrical $k$-wall} is a subdivision of $W$.
See \cref{fig:directedwall1} for an illustration.

As in the undirected case, we can translate the Directed Grid Theorem into a Directed Wall Theorem by simply doubling the necessary quantities.
Since this is just a linear factor we may use the same function as the directed grid theorem and simply double the argument.

\begin{theorem}[Directed Wall Theorem]\label{cor:directedwall}
	There exists a function $\DirectedGrid\colon\N\rightarrow \N$ such that for every $k\in\N$ and every digraph $D$ we have $\dtw{D}\leq\Fkt{\DirectedGrid}{2k}$, or $D$ contains a cylindrical $k$-wall.
\end{theorem}

Another problems that must be addressed when attempting to generalise the \textbf{Flat Wall Theorem} into the setting of digraphs is the lack of a directed version of the \textbf{Two Paths Theorem}.
It is a well known fact that the \textsc{Directed $t$-Disjoint Paths Problem} is $\NP$-complete already for $t=2$ \cite{fortune1980directed}.
And thus we cannot hope for any nice characterisation for the existence for two crossing directed paths over a cycle that can be used to build a large clique as a butterfly minor.
Even worse, since we are dealing with digraphs there might even exist disjoint crossing directed paths, but they could in fact be useless for our purposes.
Moreover, we are somewhat bound to the direction of the concentric cycles within a cylindrical wall, so there might exist many crosses which might even be useful, but if they all appear in a single `row' of the wall it is not possible for us to route through them as often as it would be necessary to create a large clique as a butterfly minor.
Because of all these problems Giannopoulou et al.\@ had to settle for a rather relaxed version of `flatness' for their directed flat wall theorem.

In what follows we introduce all necessary definitions to properly state the directed flat wall theorem of Giannopoulou et al.
We make use of many of these definitions in our proofs as well.
Especially the parametrisation of the wall is a tool that provides some convenience later on.

Suppose $\mathcal{P}=\Set{P_1,\dots,P_k}$ is a family of pairwise disjoint directed paths, and $Q$ is a directed path that meets all paths in $\mathcal{P}$ such that $P_i\cap Q$ is a directed path.
\begin{itemize}
	\item We say that the paths $P_1,\dots,P_k$ \emph{appear in this order on $Q$} if for all $i\in[1,k-1]$, $P_i\cap Q$ occurs on $Q$ strictly before $P_{i+1}\cap Q$ with respect to the orientation of $Q$.
	
	\item In this case, for $i,j\in[1,k]$ with $i<j$, we denote by $\InducedSubgraph{Q}{P_i,\dots,P_j}$ the minimal directed subpath of $Q$ containing all vertices of $Q\cap P_{\ell}$ for $\ell\in[i,j]$.
\end{itemize}

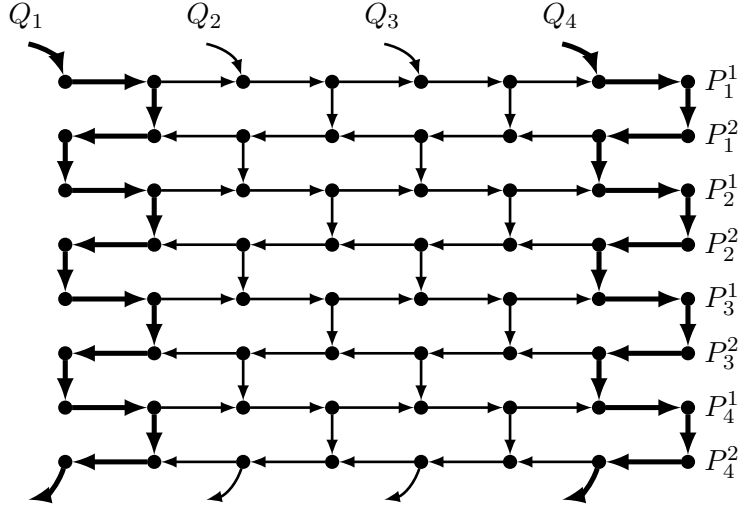
\begin{figure}[!h]
	\centering
	\begin{tikzpicture}[scale=0.9]
		\pgfdeclarelayer{background}
		\pgfdeclarelayer{foreground}
		\pgfsetlayers{background,main,foreground}
		
		\node(v11) [v:main] {};
		\node(v12) [v:main,position=0:13mm from v11] {};
		\node(v13) [v:main,position=0:13mm from v12] {};
		\node(v14) [v:main,position=0:13mm from v13] {};
		\node(v15) [v:main,position=0:13mm from v14] {};
		\node(v16) [v:main,position=0:13mm from v15] {};
		\node(v17) [v:main,position=0:13mm from v16] {};
		\node(v18) [v:main,position=0:13mm from v17] {};
		
		\node(v01) [v:ghost,position=135:8mm from v11] {};
		\node(v03) [v:ghost,position=135:8mm from v13] {};
		\node(v05) [v:ghost,position=135:8mm from v15] {};
		\node(v07) [v:ghost,position=135:8mm from v17] {};
		
		\node(v21) [v:main,position=270:8mm from v11] {};
		\node(v22) [v:main,position=0:13mm from v21] {};
		\node(v23) [v:main,position=0:13mm from v22] {};
		\node(v24) [v:main,position=0:13mm from v23] {};
		\node(v25) [v:main,position=0:13mm from v24] {};
		\node(v26) [v:main,position=0:13mm from v25] {};
		\node(v27) [v:main,position=0:13mm from v26] {};
		\node(v28) [v:main,position=0:13mm from v27] {};
		
		\node(v31) [v:main,position=270:8mm from v21] {};
		\node(v32) [v:main,position=0:13mm from v31] {};
		\node(v33) [v:main,position=0:13mm from v32] {};
		\node(v34) [v:main,position=0:13mm from v33] {};
		\node(v35) [v:main,position=0:13mm from v34] {};
		\node(v36) [v:main,position=0:13mm from v35] {};
		\node(v37) [v:main,position=0:13mm from v36] {};
		\node(v38) [v:main,position=0:13mm from v37] {};
		
		\node(v41) [v:main,position=270:8mm from v31] {};
		\node(v42) [v:main,position=0:13mm from v41] {};
		\node(v43) [v:main,position=0:13mm from v42] {};
		\node(v44) [v:main,position=0:13mm from v43] {};
		\node(v45) [v:main,position=0:13mm from v44] {};
		\node(v46) [v:main,position=0:13mm from v45] {};
		\node(v47) [v:main,position=0:13mm from v46] {};
		\node(v48) [v:main,position=0:13mm from v47] {};
		
		\node(v51) [v:main,position=270:8mm from v41] {};
		\node(v52) [v:main,position=0:13mm from v51] {};
		\node(v53) [v:main,position=0:13mm from v52] {};
		\node(v54) [v:main,position=0:13mm from v53] {};
		\node(v55) [v:main,position=0:13mm from v54] {};
		\node(v56) [v:main,position=0:13mm from v55] {};
		\node(v57) [v:main,position=0:13mm from v56] {};
		\node(v58) [v:main,position=0:13mm from v57] {};
		
		\node(v61) [v:main,position=270:8mm from v51] {};
		\node(v62) [v:main,position=0:13mm from v61] {};
		\node(v63) [v:main,position=0:13mm from v62] {};
		\node(v64) [v:main,position=0:13mm from v63] {};
		\node(v65) [v:main,position=0:13mm from v64] {};
		\node(v66) [v:main,position=0:13mm from v65] {};
		\node(v67) [v:main,position=0:13mm from v66] {};
		\node(v68) [v:main,position=0:13mm from v67] {};
		
		\node(v71) [v:main,position=270:8mm from v61] {};
		\node(v72) [v:main,position=0:13mm from v71] {};
		\node(v73) [v:main,position=0:13mm from v72] {};
		\node(v74) [v:main,position=0:13mm from v73] {};
		\node(v75) [v:main,position=0:13mm from v74] {};
		\node(v76) [v:main,position=0:13mm from v75] {};
		\node(v77) [v:main,position=0:13mm from v76] {};
		\node(v78) [v:main,position=0:13mm from v77] {};
		
		\node(v81) [v:main,position=270:8mm from v71] {};
		\node(v82) [v:main,position=0:13mm from v81] {};
		\node(v83) [v:main,position=0:13mm from v82] {};
		\node(v84) [v:main,position=0:13mm from v83] {};
		\node(v85) [v:main,position=0:13mm from v84] {};
		\node(v86) [v:main,position=0:13mm from v85] {};
		\node(v87) [v:main,position=0:13mm from v86] {};
		\node(v88) [v:main,position=0:13mm from v87] {};
		
		\node(v91) [v:ghost,position=225:8mm from v81] {};
		\node(v93) [v:ghost,position=225:8mm from v83] {};
		\node(v95) [v:ghost,position=225:8mm from v85] {};
		\node(v97) [v:ghost,position=225:8mm from v87] {};
		
		\node(Q1) [v:ghost,position=90:4mm from v01] {$Q_1$};
		\node(Q2) [v:ghost,position=90:4mm from v03] {$Q_2$};
		\node(Q3) [v:ghost,position=90:4mm from v05] {$Q_3$};
		\node(Q4) [v:ghost,position=90:4mm from v07] {$Q_4$};
		
		\node(P1_1) [v:ghost,position=0:5mm from v18] {$P_1^1$};
		\node(P1_2) [v:ghost,position=0:5mm from v28] {$P_1^2$};
		
		\node(P2_1) [v:ghost,position=0:5mm from v38] {$P_2^1$};
		\node(P2_2) [v:ghost,position=0:5mm from v48] {$P_2^2$};
		
		\node(P3_1) [v:ghost,position=0:5mm from v58] {$P_3^1$};
		\node(P3_2) [v:ghost,position=0:5mm from v68] {$P_3^2$};
		
		\node(P4_1) [v:ghost,position=0:5mm from v78] {$P_4^1$};
		\node(P4_2) [v:ghost,position=0:5mm from v88] {$P_4^2$};
		
		\begin{pgfonlayer}{background}
			
			\draw[e:main,->,line width=1.9pt] (v11) to (v12);
			\draw[e:main,->] (v12) to (v13);
			\draw[e:main,->] (v13) to (v14);
			\draw[e:main,->] (v14) to (v15);
			\draw[e:main,->] (v15) to (v16);
			\draw[e:main,->] (v16) to (v17);
			\draw[e:main,->,line width=1.9pt] (v17) to (v18);
			
			\draw[e:main,<-,line width=1.9pt] (v21) to (v22);
			\draw[e:main,<-] (v22) to (v23);
			\draw[e:main,<-] (v23) to (v24);
			\draw[e:main,<-] (v24) to (v25);
			\draw[e:main,<-] (v25) to (v26);
			\draw[e:main,<-] (v26) to (v27);
			\draw[e:main,<-,line width=1.9pt] (v27) to (v28);
			
			\draw[e:main,->,line width=1.9pt] (v31) to (v32);
			\draw[e:main,->] (v32) to (v33);
			\draw[e:main,->] (v33) to (v34);
			\draw[e:main,->] (v34) to (v35);
			\draw[e:main,->] (v35) to (v36);
			\draw[e:main,->] (v36) to (v37);
			\draw[e:main,->,line width=1.9pt] (v37) to (v38);
			
			\draw[e:main,<-,line width=1.9pt] (v41) to (v42);
			\draw[e:main,<-] (v42) to (v43);
			\draw[e:main,<-] (v43) to (v44);
			\draw[e:main,<-] (v44) to (v45);
			\draw[e:main,<-] (v45) to (v46);
			\draw[e:main,<-] (v46) to (v47);
			\draw[e:main,<-,line width=1.9pt] (v47) to (v48);
			
			\draw[e:main,->,line width=1.9pt] (v51) to (v52);
			\draw[e:main,->] (v52) to (v53);
			\draw[e:main,->] (v53) to (v54);
			\draw[e:main,->] (v54) to (v55);
			\draw[e:main,->] (v55) to (v56);
			\draw[e:main,->] (v56) to (v57);
			\draw[e:main,->,line width=1.9pt] (v57) to (v58);
			
			\draw[e:main,<-,line width=1.9pt] (v61) to (v62);
			\draw[e:main,<-] (v62) to (v63);
			\draw[e:main,<-] (v63) to (v64);
			\draw[e:main,<-] (v64) to (v65);
			\draw[e:main,<-] (v65) to (v66);
			\draw[e:main,<-] (v66) to (v67);
			\draw[e:main,<-,line width=1.9pt] (v67) to (v68);
			
			\draw[e:main,->,line width=1.9pt] (v71) to (v72);
			\draw[e:main,->] (v72) to (v73);
			\draw[e:main,->] (v73) to (v74);
			\draw[e:main,->] (v74) to (v75);
			\draw[e:main,->] (v75) to (v76);
			\draw[e:main,->] (v76) to (v77);
			\draw[e:main,->,line width=1.9pt] (v77) to (v78);
			
			\draw[e:main,<-,line width=1.9pt] (v81) to (v82);
			\draw[e:main,<-] (v82) to (v83);
			\draw[e:main,<-] (v83) to (v84);
			\draw[e:main,<-] (v84) to (v85);
			\draw[e:main,<-] (v85) to (v86);
			\draw[e:main,<-] (v86) to (v87);
			\draw[e:main,<-,line width=1.9pt] (v87) to (v88);
			
			\draw[e:main,->,line width=1.9pt] (v12) to (v22);
			\draw[e:main,->] (v14) to (v24);
			\draw[e:main,->] (v16) to (v26);
			\draw[e:main,->,line width=1.9pt] (v18) to (v28);
			
			\draw[e:main,->,line width=1.9pt] (v21) to (v31);
			\draw[e:main,->] (v23) to (v33);
			\draw[e:main,->] (v25) to (v35);
			\draw[e:main,->,line width=1.9pt] (v27) to (v37);
			
			\draw[e:main,->,line width=1.9pt] (v32) to (v42);
			\draw[e:main,->] (v34) to (v44);
			\draw[e:main,->] (v36) to (v46);
			\draw[e:main,->,line width=1.9pt] (v38) to (v48);
			
			\draw[e:main,->,line width=1.9pt] (v41) to (v51);
			\draw[e:main,->] (v43) to (v53);
			\draw[e:main,->] (v45) to (v55);
			\draw[e:main,->,line width=1.9pt] (v47) to (v57);
			
			\draw[e:main,->,line width=1.9pt] (v52) to (v62);
			\draw[e:main,->] (v54) to (v64);
			\draw[e:main,->] (v56) to (v66);
			\draw[e:main,->,line width=1.9pt] (v58) to (v68);
			
			\draw[e:main,->,line width=1.9pt] (v61) to (v71);
			\draw[e:main,->] (v63) to (v73);
			\draw[e:main,->] (v65) to (v75);
			\draw[e:main,->,line width=1.9pt] (v67) to (v77);
			
			\draw[e:main,->,line width=1.9pt] (v72) to (v82);
			\draw[e:main,->] (v74) to (v84);
			\draw[e:main,->] (v76) to (v86);
			\draw[e:main,->,line width=1.9pt] (v78) to (v88);
			
			\draw[e:main,bend left=30,->,line width=1.9pt] (v01) to (v11);
			\draw[e:main,bend left=30,->] (v03) to (v13);
			\draw[e:main,bend left=30,->] (v05) to (v15);
			\draw[e:main,bend left=30,->,line width=1.9pt] (v07) to (v17);
			
			\draw[e:main,bend left=30,->,line width=1.9pt] (v81) to (v91);
			\draw[e:main,bend left=30,->] (v83) to (v93);
			\draw[e:main,bend left=30,->] (v85) to (v95);
			\draw[e:main,bend left=30,->,line width=1.9pt] (v87) to (v97);
			
		\end{pgfonlayer}
	\end{tikzpicture}
	\caption{The elementary cylindrical $4$-wall. The thick edges of the cycles $Q_1$ and $Q_4$ mark its perimeter.}
	\label{fig:cylindrical4wall}
\end{figure}

It is convenient for us to imagine cylindrical walls, and similarly, as we will see later, their matching theoretic analogues, as illustrated in \cref{fig:cylindrical4wall}.
So we think of the concentric cycles as paths from the top to the bottom with an additional edge from the lowest row back to the top.
The directed paths that alternately go in and out of the wall are then seen as the horizontal paths.

Since the underlying undirected graph is planar and $3$-connected Whitney's Theorem \cite{whitney1992congruent} ensures a unique planar embedding.
In particular, we will refer to the faces of the embedding as the faces of the graph itself.
A face that is bounded by a cycle that is not the perimeter is called a \emph{cell} of the wall.

\begin{definition}[Vertical and Horizontal Paths]
	Let $k\in\N$ be a positive integer and $W$ be a cylindrical $k$-wall.
	
	We denote the vertical paths of $W$ by $Q_1,\dots,Q_k$, ordered from left to right.
	Let $\CondSet{P_j^i}{i\in[1,2],~j\in[1,k]}$ be the horizontal directed paths such that the paths $P^1_j$, $j\in[1,k]$, are oriented from left to right and the paths $P^2_j$, $j\in[1,k]$, are oriented from right to left such that $P^i_j$ is above $P^{i'}_{j'}$ whenever $j<j'$ and $P^1_j$ is above $P^2_j$ for all $j\in[1,k]$.
	The top line is $P_1^1$.
	
	By $\hat{P}_j$ we denote the disjoint union of $P_j^1$ and $P_j^2$ for all $j\in[1,k]$.
	
	Two horizontal paths $P_j^i$ and $P_{j'}^{i'}$ are \emph{consecutive} if $i\neq i'$, and $j'\in[j-1,j+1]$ or if $P_j^i+P_{j'}^{i'}=\hat{P}_j$.
	A family $\mathcal{P}\subseteq \CondSet{P_j^i}{i\in[1,2],~j\in[1,k]}$ is said to be \emph{consecutive} if there do not exist paths $P_1,P_2\in\mathcal{P}$, and $P_3\in\CondSet{P_j^i}{i\in[1,2],~j\in[1,k]}\setminus\mathcal{P}$ such that there is no directed path from $P_1$ to $P_2$ in $W-P_3$.
	We extend our notation for $\InducedSubgraph{P_j^i}{Q_p,\dots,Q_q}$ for $p<q$ in the natural way for $\InducedSubgraph{\mathcal{P}}{Q_p,\dots,Q_q}$ and, in a slight abuse of notation, identify $\hat{P}_i$ and $\Set{P^1_i,P^2_i}$.
	
	For more convenience we write ``Let $W=\Brace{Q_1,\dots,Q_k,\hat{P}_1,\dots,\hat{P}_k}$ be a cylindrical $k$-wall.'' to fix the embedding and naming of the vertical cycles and horizontal paths as explained above and depicted in \cref{fig:cylindrical4wall} for the case $k=4$.
\end{definition}

\begin{definition}[Perimeter]\label{def:perimeter}
	Let $k\in\N$ be a positive integer and $W=\Brace{Q_1,\dots,Q_k,\hat{P}_1,\dots,\hat{P}_k}$ be a cylindrical $k$-wall.
	The vertical cycles $Q_1$ and $Q_k$ each bounds a face of $W$, we define their union to be the \emph{perimeter} of $W$ an write $\Perimeter{W}=Q_1\cup Q_k$.
\end{definition}

\begin{definition}[$W$-Distance]\label{def:Wdistance}
	Let $k\in\N$ be a positive integer and $W=\Brace{Q_1,\dots,Q_k,\hat{P}_1,\dots,\hat{P}_k}$ be a cylindrical $k$-wall.
	Given two vertices $u,v\in\V{W}$, we say that they have \emph{$W$-distance} at least $i$ if there exist $i$ distinct vertical or $i$ distinct horizontal paths whose removal separates $u$ and $v$ in $W$.
\end{definition}

\begin{definition}[Slice]\label{def:slice}
	Let $k\in\N$ be a positive integer and $W=\Brace{Q_1,\dots,Q_k,\hat{P}_1,\dots,\hat{P}_k}$ be a cylindrical $k$-wall.
	A \emph{slice} $W'$ of $W$ is a cylindrical wall containing the vertical paths $Q_i,\dots,Q_{i+\ell}$ for all $i\in[1,k]$ and some $\ell\in[1,k-i]$, and the horizontal paths $\InducedSubgraph{P_1^1}{Q_i,\dots,Q_{i+\ell}},\dots,\InducedSubgraph{P_k^2}{Q_i,\dots,Q_{i+\ell}}$.
	We say that $W'$ is the \emph{slice of $W$ between $Q_i$ and $Q_{i+\ell}$} and that it is of \emph{width $\ell+1$}.
\end{definition}

\begin{definition}[Strip]\label{def:strip}
	Let $k\in\N$ be a positive integer and $W=\Brace{Q_1,\dots,Q_k,\hat{P}_1,\dots,\hat{P}_k}$ be a cylindrical $k$-wall.
	A \emph{strip of height $j-i+1$ between $i$ and $j$} of $W$ is the subgraph of $W$ induced by the horizontal paths $\hat{P}_i,\dots,\hat{P}_j$ for some $i<j\in[1,k]$ and the subpaths $\InducedSubgraph{Q_{\ell}}{\hat{P}_i,\dots,\hat{P}_j}$ for $\ell\in[1,k]$.
\end{definition}

\begin{figure}[!h]
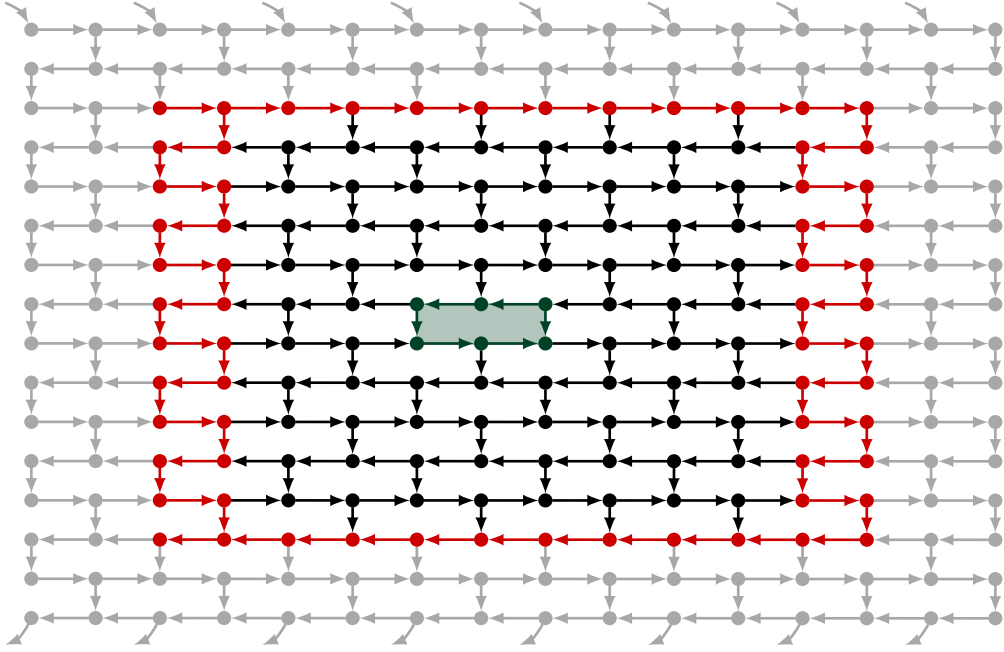

	\centering

	\caption{The tile $T$ at $\Brace{2,2}$ of width $2$ in a cylindrical $8$-wall.
		The green cell is the centre, the red paths are the perimeter of $T$, and the vertices that belong to $T$, but neither to the red paths, nor the green cell, are the internal vertices of $T$.}
	\label{fig:atile}
\end{figure}

\begin{definition}[Tiles]\label{def:tile}
	Let $k\in\N$ be a positive integer and $W=\Brace{Q_1,\dots,Q_k,\hat{P}_1,\dots,\hat{P}_k}$ be a cylindrical $k$-wall.
	Let $i,j\in[1,k]$ and $d\in\N$ be positive integers.
	The \emph{tile $T$ of $W$ at $\Brace{i,j}$ of width $d$} is defined as the subgraph of $W$ induced by
	\begin{align*}
		\bigcup_{\ell\in[i,i+2d+1]}\InducedSubgraph{Q_{\ell}}{\hat{P}_j,\dots,\hat{P}_{j+2d+1}}\cup\bigcup_{\ell\in[j,j+2d+1]}\InducedSubgraph{\hat{P}_{\ell}}{Q_i,\dots,Q_{i+2d+1}}.
	\end{align*}
	We call $i$ the \emph{column index} of the tile, $j$ the \emph{row index} of the tile, and say that the $j$-th row has a tile $T$ if the row index of $T$ is $j$.
	To make the notation a bit more compact we write $T_{i,j,d}$ for the tile of $W$ at $\Brace{i,j}$ of width $w$.
	
	Since $W$ is a cylindrical wall, there exist subgraphs of $W$ that technically also form tiles, but that do not necessarily fit into our parametrisation of $W$.
	To overcome this, we agree for $\ell>k$ to set $P_{\ell}\coloneqq P_{\Brace{\Brace{\ell-1}\mod k}+1}$.
	This means that tiles that start near the bottom are allowed to continue at the top.
	Indeed, the notions of top and bottom are only present because of the way we parametrised the wall, and thus even those tiles are well defined.
	
	The \emph{perimeter} of the tile $T$ is $T\cap\Brace{Q_i\cup Q_{i+2d+1}\cup P_j^1\cup P_{j+2d+1}^2}$.
	We call $Q_i$ the \emph{left path} of the perimeter, $Q_{i+2d+1}$ its \emph{right path}, $P_j^1$ the \emph{upper path} of the perimeter, and finally $P_{j+2d+1}^2$ is its \emph{lower path}.
	
	The \emph{corners} of a tile are the vertices $a,b,c,d\in\V{T}$ where $a$, the \emph{upper left corner}, is the common starting point of $T\cap Q_i$ and $T\cap P_j^1$, $b$, the \emph{upper right corner}, is the end of $T\cap P_j^1$ and the starting point of $Q_{i+2d+1}$, $c$, the \emph{lower left corner}, is the common end of $T\cap Q_i$ and $T\cap P_{j+2d+1}^2$, finally $d$, the \emph{lower right corner}, is the end of $T\cap Q_{i+2d+1}$ and the starting point of $T\cap P_{j+2d+1}^2$.
	
	The \emph{centre} of $T$ is the boundary of the unique cell $C_T$ of $W$ whose boundary consists of vertices from $Q_{i+d+1}$, $Q_{i+d+2}$, $P^2_{j+d+1}$, and $P^1_{j+d+2}$.
	All vertices of $T$ which are not in the centre and not on the perimeter of $T$ are called \emph{internal}.
	See \cref{fig:atile} for an illustration of a tile.
\end{definition}

Please note that by this definition, only cells that lie between $P_i^1$ and $P_i^2$ for some $i\in[1,k]$ can be centre of a tile.
However, if we were to take the mirror image of our currently fixed embedding along a straight vertical line between $Q_{\Floor{\frac{k}{2}}}$ and $Q_{\Floor{\frac{k}{2}}+1}$, we obtain a new embedding for which we then can reapply our parametrisation.
By doing so, every path $P_i^1$ now becomes a path $P_{i'}^2$, and $P_i^2$ becomes $P_{i''}^1$ for $i,i',i''\in[1,k]$.
This means that we can define for every cell $F$ of $W$ a tile $T_F$ such that $F$ is the centre of $T_F$.

We are now finally ready to state the definition of `flatness' used in the directed flat wall theorem.

Let $D$ be a digraph and $A\subseteq\V{D}$ be a set of vertices.
Let $d,t\geq 1$ be integers and let $W\subseteq D-A$ be a cylindrical wall.
We say that $W$ is \emph{almost $A$-flat} in $D-A$ \emph{with directed treewidth bounded by $d$} if the following holds.
\begin{enumerate}
	\item There is a separation\footnote{Note that $(X,Y)$ is indeed a \hyperref[def:separation]{separation}, that is no edge of $D$ has one endpoint in $X\setminus Y$ and the other in $Y\setminus X$. We emphasise this since we are in a digraph and $(X,Y)$ could also denote a directed separation - which it does \textbf{not} here.} $(X,Y)$ of $D$ such that $X\cap Y= A\cup\Perimeter{W}$, $W-\Perimeter{W}\subseteq Y$, and every vertex in $Y$ reaches a vertex of $W-\Perimeter{W}$ or is reachable from it.
	
	\item For every path $Q$ in $D-A$ that is internally disjoint from $W-\Perimeter{W}$, such that $Q$ has end endpoint in $W-\Perimeter{W}$ and the other in $W$ there is a cell $C$ of $W$ such that the boundary of $C$ contains both endpoints of $Q$.
	Furthermore, for every cell $C$ of $W$, if $Z$ is the set of vertices of a path $P$ in $D-A$ with both endpoints on $C$ and internally disjoint from $W$, then the components of $\InducedSubgraph{D}{Z}$, which are called \emph{extensions}, have directed treewidth at most $d$.
	
	\item If $T$ is a tile $W$ of width five and $c_1, c_2$ are its two upper corners from left to right and $d_1, d_2$ are its two lower corners from left to right, then there are no two disjoint paths $P_1$, $P_2$ in $D-A-(W-T)$ connecting $c_1$ to $d_2$ and $c_2$ to $d_2$.
\end{enumerate}
If $W$ has at most $t$ rows whose tiles do not satisfy property \emph{iii)}, but $W$ satisfies properties \emph{i)} and \emph{ii)} then we say that $W$ is \emph{$t$-barely $A$-flat}.

Let us investigate this definitions and the flaws in the resulting notion of flatness more closely.
Property \emph{i)} states that, after removing the apex set $A$, the perimeter of $W$ separates $D$ into two parts, an `outside' part which has no connection to $W$, and an `inside' part which contains all of $W$.
Moreover, by property \emph{ii)}, the cells of $W$ already kind of represent some notion of planarity in the sense that every path connecting two distinct cells within the `inside' part must contain vertices at least from the boundaries of the two cells which are connected.
Additionally, the requirement that all attachments must have bounded directed treewidth is also something one can achieve in the undirected case.
In the directed version however it can be seen as an attempt to impose some further control over the internal connectivities via disjoint paths to at least approximate the implications from the \textbf{Flat Wall Theorem}.
Finally, and this circles back our discussion of `usable' and `unusable' crosses from before, if we consider a tile whose width is at least some constant, then we can be sure that there is no cross over this tile that goes in the same direction as the concentric cycles of the wall.
Except in case $W$ is only $t$-barely $A$-flat.
Here some of these crosses might exist, but only in a restricted area of the wall.
Moreover, if $T$ is a tile with corners $c_1,c_2$, $d_1,d_2$ as in property \emph{iii)} then a cross connecting $d_1$ to $c_2$ and $d_2$ to $c_1$ can still exist without any problem.
So while the notions of almost flatness and barely flatness manage to capture some parts of their undirected counterpart, they are inherently flawed.
Indeed, it seems like there is no way to avoid these problems as they arise from the highly restrictive nature of the cylindrical grid itself as well as the fact that \textsc{Directed $2$-Disjoint Paths} is $\NP$-complete.

\begin{theorem}[Directed Flat Wall Theorem, \cite{giannopoulou2020directed}]\label{thm:directedflatwall}
	There exist functions $a\colon\N\rightarrow\N$, $b\colon\N\rightarrow\N$, and $d\colon\N\times\N\rightarrow\N$ such that for all integers $t,r\in\N$ and every digraph $D$ one of the following is true:
	\begin{enumerate}
		\item $\dtw{D}\leq\Fkt{d}{t,r}$,
		
		\item $D$ contains $\Bidirected{K_t}$ as a butterfly minor, or
		
		\item there exist a set $A\subseteq\V{D}$ with $\Abs{A}\leq\Fkt{a}{t}$ and a cylindrical $r$-wall $W\subseteq D-A$ which is $\Fkt{b}{t}$-barely $A$-flat in $G-D$ with directed treewidth bounded by $\Fkt{d}{r,t}$.
	\end{enumerate}
\end{theorem}

Please note the difference in the way \cref{thm:directedflatwall} is stated in comparison to \cref{thm:undirectedflatwall}.
The later guarantees a flat subwall for every wall of sufficient size, while the directed version can only guarantee the existence of such a wall.
This discrepancy is due to the definition of barely flatness which requires the attachments to have bounded directed treewidth.
The way to prove the existence of such a wall is simply to start with any wall, make sure all other requirements are met and then, if there is some attachment of larger directed treewidth, move into this attachment and reiterate the arguments.
It is certainly possible to find a version of \cref{thm:directedflatwall} which allows every large enough wall to contain a somewhat flat subwall, but this would mean to sacrifice the bounded directed treewidth of the attachments.

\section{Matching Theoretic Background}\label{sec:background}

In this section we aim to introduce the matching theoretic concepts necessary for this paper.
For a deeper introduction to Matching Theory the reader may consult \cite{lovasz2009matching}.

Let $G$ be a graph, a set $F\subseteq\E{G}$ of pairwise disjoint edges is called a \emph{matching}, and a vertex $v\in\V{G}$ is said to be \emph{covered} by $F$ if $F$ contains an edge that has $v$ as an endpoint.
The set of all vertices covered by $F$ is denoted by $\V{F}$.
A matching $M$ is \emph{perfect} if $\V{M}=\V{G}$ and an edge $e\in\E{G}$ is \emph{admissible} if there exists a perfect matching $M'$ of $G$ with $e\in M'$.
A graph $G$ is \emph{matching covered} if it is connected and all of its edges are admissible.
A set $X\subseteq\V{G}$ is \emph{conformal} if $G-X$ has a perfect matching.
If $M$ is a perfect matching and $\E{G-X}\cap M$ is a perfect matching of $G-X$, we call $X$ \emph{$M$-conformal}.
A subgraph $H$ of $G$ for which $\V{H}$ is ($M$-)conformal is called \emph{($M$-)conformal}.

\begin{definition}[Matching Minors]\label{def:matchingminor}
Let $G$ be a graph with a perfect matching.
If $v\in\V{G}$ is a vertex of degree exactly two, we call the process of contracting both edges incident with $v$ and removing all loops and parallel edges afterwards the \emph{bicontraction} of $v$.
A graph $H$ that can be obtained from a conformal subgraph of $G$ by repeatedly applying bicontractions is called a \emph{matching minor}.
We say that $H$ is an \emph{$M$-minor} of $G$ if $H$ can be obtained from an $M$-conformal subgraph $H'$ of $G$ by repeated bicontractions where $M$ contains a perfect matching of $H'$.
\end{definition}

In many cases we are interested in a topological variant of matching minors.
Let $G$ be a graph.
We say that a graph $H$ is a \emph{bisubdivision} of $G$ if $H$ can be obtained from $G$ by replacing every edge of $G$ with a path of odd length such that any two of these paths are internally disjoint.

\begin{lemma}[\cite{lucchesi2018two}]\label{lemma:confmathingminors}
	Let $G$ and $H$ be matching covered graphs such that $\MaximumDegree{H}=3$.
	Then $G$ contains a conformal bisubdivision of $H$ if and only if it contains $H$ as a matching minor.
\end{lemma}

The idea of matching minors stems from a more general concept used for the investigation of matching covered graphs.

\begin{definition}[Tight Cuts and Braces]\label{def:braces}
Let $G$ be a graph and $X\subseteq\V{G}$.
We denote by $\Cut{X}$ the set of edges in $G$ with exactly one endpoint in $X$ and call $\Cut{X}$ the \emph{edge cut} around $X$ in $G$.
Let us denote by $\Perf{G}$ the set of all perfect matchings in $G$.
An edge cut $\Cut{X}$ is \emph{tight} if $\Abs{\Cut{X}\cap M}=1$ for all perfect matchings $M\in\Perf{G}$.
If $\Cut{X}$ is a tight cut and $\Abs{X},\Abs{\V{G}\setminus X}\geq 2$, it is \emph{non-trivial}.
Identifying the \emph{shore} $X$ of a non-trivial tight cut $\Cut{X}$ into a single vertex is called a \emph{tight cut contraction} and the resulting graph $G'$ can easily be seen to be matching covered again.
A bipartite matching covered graph without non-trivial tight cuts is called a \emph{brace}.	
\end{definition}

The following is an observation first made by Lov\'asz (see the proof of Lemma 1.4 in \cite{lovasz1987matching}).

\begin{observation}\label{obs:tightcutminoritymajority}
Let $B$ be a bipartite and matching covered graph and let $X\subseteq\V{B}$ be a set of vertices with $\Abs{X},\Abs{\V{B}\setminus X}\geq 3$.
Then $\CutG{}{X}$ is a tight cut if and only if there exists $i\in[1,2]$ such that $\Abs{X\cap V_i}=\Abs{X\setminus V_i}-1$, and $\Fkt{N_B}{X\cap V_i}\subseteq X$.
We say that $V_i$ is the \emph{minority} of $X$ and $V_{3-i}$ is the \emph{majority} of $X$.
\end{observation}

It follows from a famous result of Lov\'asz \cite{lovasz1987matching} that the braces of a bipartite matching covered graph $B$ are uniquely determined.
Similar to how every $2$-connected minor of a graph $G$ must be a minor of one of its blocks, every brace that is a matching minor of some bipartite matching covered graph $B$ is a matching minor of some brace of $B$ \cite{lucchesi2015thin}.
Moreover, the braces of a bipartite matching covered graph play a similar role for matching theory as blocks, i.\@e.\@ maximal $2$-connected subgraphs, do in general graph theory.
This is even more emphasised by the following classical result of Plummer.

\begin{definition}[Extendibility]\label{def:extendibility}
Let $k$ be a positive integer.
A graph $G$ is called \emph{$k$-extendible} if it has at least $2k+2$ vertices and for every matching $F\subseteq\E{G}$ of size $k$ there exists a perfect matching $M$ of $G$ with $F\subseteq M$.
\end{definition}

\begin{theorem}[\cite{lovasz2009matching}]\label{thm:braces}
	A bipartite graph $B$ is a brace if and only if it is either isomorphic to $C_4$, or it is $2$-extendible.
\end{theorem}

The following theorem is a collection of several different characterisations of $k$-extendibility in bipartite graphs.

\begin{theorem}[\cite{plummer1986matching,aldred2003m}]\label{thm:bipartiteextendibility}
	Let $B$ be a bipartite graph and $k\in\N$ a positive integer.
	The following statements are equivalent.
	\begin{enumerate}
		\item $B$ is $k$-extendible.
		\item $\Abs{V_1}=\Abs{V_2}$, and for all non-empty $S\subseteq V_1$, $\Abs{\NeighboursG{B}{S}}\geq \Abs{S}+k$.
		\item For all sets $S_1\subseteq V_1$ and $S_2\subseteq V_2$ with $\Abs{S_1}=\Abs{S_2}\leq k$ the graph $B-S_1-S_2$ has a perfect matching.
		\item $B$ is matching $k$-connected.
		\item There is a perfect matching $M\in\Perf{B}$ such that for every $v_1\in V_1$, every $v_2\in V_2$ there are $k$ pairwise internally disjoint internally $M$-conformal paths with endpoints $v_1$ and $v_2$.
		\item For every perfect matching $M\in\Perf{B}$, every $v_1\in V_1$, every $v_2\in V_2$ there are $k$ pairwise internally disjoint internally $M$-conformal paths with endpoints $v_1$ and $v_2$.
	\end{enumerate}
\end{theorem}

While ii) can be seen as a generalisation of Hall's Theorem, iii) is, in a way, an even stronger version of matching $k$-connectivity.
The statements iv) and v) can be seen as matching theoretic versions of Menger's Theorem.

\begin{theorem}[Menger's Theorem \cite{menger1927allgemeinen}]\label{thm:directedlocalmenger}
	Let $D$ be a digraph and $X,Y\subseteq\V{D}$ be two sets of vertices, then the maximum number of pairwise disjoint directed $X$-$Y$-paths in $D$ equals the minimum size of a set $S\subseteq\V{G}$ such that every directed $X$-$Y$-path in $D$ contains a vertex of $S$.
\end{theorem}

\paragraph{The \textsc{Pfaffian Recognition} Problem}

To obtain a matching theoretic version of the \textbf{Flat Wall Theorem}, that is a version that interacts with \hyperref[def:matchingminor]{matching minors} in bipartite graphs, we need a matching theoretic notion of flatness.
To understand where we derive our flatness from we need to know a bit about the structure of matching covered bipartite graphs that exclude $K_{3,3}$ as a matching minor.

The graph $K_{3,3}$ plays a key role in the theory of matching minors.
It was found to be the singular obstruction, in the sense of matching minors, for a bipartite graph to have a Pfaffian Orientation relatively early \cite{kasteleyn1967graph}.
At the time however this did not yield a solution for the \textsc{Pfaffian Recognition} problem as no algorithm was known to check for the presence of a specific matching minor.
In the following $30$ years many equivalent problems would be discovered by various authors (see \cite{mccuaig2004polya} for a good overview) but it took a complete structural characterisation of bipartite graphs without a $K_{3,3}$-matching minor that resembles similar results from (regular) minor theory such as, for example, Wagner's description of $K_5$-minor free graphs \cite{wagner1937eigenschaft}.

Given a bipartite graph $B$ with a perfect matching we say that $B$ \emph{contains} $K_{3,3}$ if it has a matching minor isomorphic to $K_{3,3}$, in this case we often say that $B$ is \emph{non-Pfaffian}.
If $B$ does not contain $K_{3,3}$ we say that $B$ is \emph{$K_{3,3}$-free} or \emph{Pfaffian}.

\begin{definition}[$4$-Cycle Sum]\label{def:cyclesum}
	For every $i\in\Set{1,2,3}$ let $B_i$ be a bipartite graph with a perfect matching and $C_i$ be a conformal cycle of length four in $B_i$.
	A \emph{$4$-cycle-sum} of $B_1$ and $B_2$ at $C_1$ and $C_2$ is a graph $B'$ obtained by identifying $C_1$ and $C_2$ into the cycle $C'$ and possibly forgetting some of its edges.
	If a bipartite graph $B''$ is a $4$-cycle-sum of $B'$ and some bipartite and matching covered graph $B_3$ at $C'$ and $C_3$, then $B''$ is called a \emph{trisum} of $B_1$, $B_2$ and $B_3$.
\end{definition}

\begin{figure}[!h]
	\centering
	\begin{tikzpicture}[scale=0.9]
		\pgfdeclarelayer{background}
		\pgfdeclarelayer{foreground}
		\pgfsetlayers{background,main,foreground}
		
		\foreach \x in {2,4,6,8,10,12,14}
		{
			\node[v:main] () at (\x*25.71:15mm){};
		}
		
		\foreach \x in {1,3,5,7,9,11,13}
		{
			\node[v:mainempty] () at (\x*25.71:15mm){};
		}
		
		\begin{pgfonlayer}{background}
			\foreach \x in {2,4,6,8,10,12,14}
			{
				\draw[e:mainthin] (\x*25.71:15mm) to (128.55+\x*25.71:15mm);
			}
			
			\foreach \x in {1,...,14}
			{
				\draw[e:main] (\x*25.71:15mm) to (25.71+\x*25.71:15mm);
			}
		\end{pgfonlayer}
	\end{tikzpicture}
	\caption{The Heawood graph $H_{14}$.}
	\label{fig:heawood}
\end{figure}
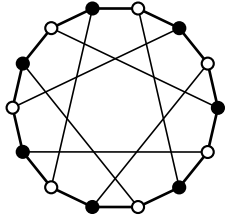

The \emph{Heawood graph} is the bipartite graph associated with the incidence matrix of the Fano plane, see \cref{fig:heawood} for an illustration.
It is the singular exceptional graph in the structure theorem for $K_{3,3}$-free braces, similar to how Wagner's graph is the only non-planar graph necessary to describe the structure of al $K_5$-minor free graphs.

\begin{theorem}\cite{mccuaig2004polya,robertson1999permanents}\label{thm:trisums}
	A brace is $K_{3,3}$-free if and only if it either is isomorphic to the Heawood graph, or it can be obtained from planar braces by repeated application of the \hyperref[def:cyclesum]{trisum} operation.
\end{theorem}

\begin{corollary}\cite{mccuaig2004polya,robertson1999permanents}\label{cor:pfaffianalg}
	There exists an algorithm that decides, given a brace $B$ as input, whether $B$ contains $K_{3,3}$ as a matching minor in time $\Fkt{\mathcal{O}}{\Abs{\V{B}}^3}$.
\end{corollary}

\paragraph{A Matching Theoretic Two Paths Theorem}

In $K_{3,3}$-free braces that are distinct from the Heawood graph we may replace any occurrences of non-planarity by $4$-cycles.
Moreover, if we are given a highly connected \hyperref[def:matchingminor]{matching minor} within a $K_{3,3}$-free brace, say for example a large grid, we know that this grid cannot be contained in many different parts which are separated by $4$-cycles and thus we must be able to find a planar brace within the tree-structure provided by \cref{thm:trisums} that contains our large grid.
This gives a solid idea of what flatness might be like in bipartite matching covered graphs.
So next we need a tool that allows us to argue that, in the absence of a large complete bipartite graph as a matching minor, we will find an area within every grid which is, essentially, $K_{3,3}$-free.
In \cite{giannopoulou2021two} we found an analogue of the \textbf{Two Paths Theorem} for crosses formed by alternating paths over conformal cycles in bipartite graphs with perfect matchings based on \hyperref[def:cyclesum]{$4$-cycle sums}.

\begin{definition}[Conformal Cross]\label{def:matchingcross}
Let $G$ be a graph with a perfect matching and let $C$ be a conformal cycle in $G$.
Two paths $P_1$ and $P_2$ where $P_i$ has endpoints $s_i$ and $t_i$ for each $i\in[1,2]$ are said to form a \emph{conformal cross} over $C$ if they are disjoint, internally disjoint from $C$, their endpoints $s_1$, $s_2$, $t_1$, and $t_2$ occur on $C$ in the order listed, and the graph $C+P_1+P_2$ is a conformal subgraph of $G$.
\end{definition}

\begin{definition}[Alternating Path]\label{def:alternatingpath}
	Let $G$ be a graph with a perfect matching $M$.
	A path $P$ in $G$ is \emph{$M$-alternating} if there exists a set $S\subseteq\V{P}$ of endpoints of $P$ such that $P-S$ is a conformal subgraph of $G$ and we say that $P$ is \emph{alternating} if there exists a perfect matching $M$ of $G$ such that $P$ is $M$-alternating.
	$P$ is $M$-conformal if $S=\emptyset$ and $P$ is \emph{internally $M$-conformal} if $S$ contains both endpoints of $P$.
\end{definition}

Please note that, if $P_1$ and $P_2$ form a \hyperref[def:matchingcross]{conformal cross} over some conformal cycle $C$ in a graph $G$, then there exists a perfect matching $M$ of $G$ such that both $P_1$ and $P_2$ are internally $M$-conformal and $C$ is an $M$-conformal cycle.
The following two results are the key to the main theorem of \cite{giannopoulou2021two} and will play the role of the \textbf{Two Paths Theorem} in this paper.

\begin{lemma}[\cite{giannopoulou2021two}]\label{lemma:goodcrossesmeanK33}
	Let $B$ be a brace and $C$ a $4$-cycle in $B$, then there is a conformal cross over $C$ in $B$ if and only if $C$ is contained in a conformal bisubdivision of $K_{3,3}$.	
\end{lemma}

\begin{theorem}[\cite{giannopoulou2021two}]\label{thm:4cycleK33}
	Let $B$ be a brace containing $K_{3,3}$ and $C$ a $4$-cycle in $B$, then there exists a conformal bisubdivision of $K_{3,3}$ with $C$ as a subgraph.
\end{theorem}

\subsection{A Relation Between Digraphs and Bipartite Graphs with Perfect Matchings}\label{subsec:digraphsandmatchings}

In many cases it suffices to fix a single perfect matching $M$ of a bipartite graph $B$ to analyse its matching theoretic properties.
This approach can lead to a simplification of the matter at hand or at least simplify some notation.
When working with bipartite graphs, fixing a single perfect matching yields an interesting and deep connection to the setting of digraphs.
Let $B$ be a bipartite graph with a perfect matching $M$ and consider the edges of $\E{B}\setminus M$.
Each such edge connects a black endpoint of some edge $e\in M$ to the white endpoint of another edge $e'\in M\setminus\Set{e}$.
Instead of encoding the bipartition of $B$ through colours we could also assign a direction to every edge of $\E{B}\setminus M$ fixing the convention that the tail of an oriented edge should always be black, while its head should always be white.
If we now contract each edge of $M$ individually into a vertex, the resulting graph is simply a digraph.

\begin{definition}[$M$-Direction]\label{def:Mdirection}
	Let $B=\Brace{V_1\cup V_2, E}$ be a bipartite graph and let $M\in\Perf{G}$ be a perfect matching of $B$. 
	The \emph{$M$-direction} $\DirM{B}{M}$ of $B$ is defined as follows.
	\begin{enumerate}
		\item $\Fkt{V}{\DirM{G}{M}}\coloneqq M$ and
		
		\item $\Fkt{E}{\DirM{G}{M}}\coloneqq\left\{\Brace{e,f}\in\Choose{M}{2}\text{there is $g\in\E{B}$ such that $\emptyset\neq e\cap g\subseteq V_1$ and}\right.$\\
		\phantom{x} ~~~~~~~~~~~~~~ $\left.\phantom{\Choose{M}{2}}\text{$\emptyset\neq f\cap g\subseteq V_2$}\right\}$.	
	\end{enumerate}
\end{definition}

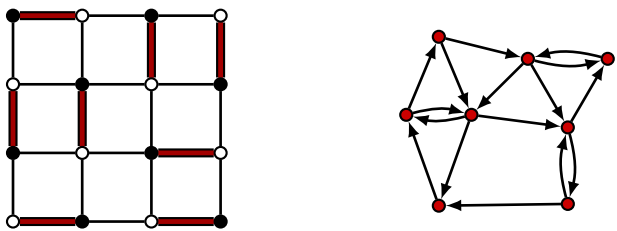
\begin{figure}[h!]
	\begin{center}
		\begin{tikzpicture}[scale=0.7]
			
			\pgfdeclarelayer{background}
			\pgfdeclarelayer{foreground}
			
			\pgfsetlayers{background,main,foreground}

			\begin{pgfonlayer}{main}
				
				\node (C) [] {};
				
				\node (C1) [v:ghost, position=180:40mm from C] {};
				
				\node (C2) [v:ghost, position=0:0mm from C] {};
				
				\node (C3) [v:ghost, position=0:40mm from C] {};

				
				
				\node (a) [v:main,position=0:0mm from C1] {};
				\node (ain) [v:ghost,position=0:2.25mm from a] {};
				
				\node (b) [v:main,position=0:13mm from a,fill=white] {};
				\node (e) [v:main,position=270:13mm from a,fill=white] {};
				
				\node (c) [v:main,position=0:13mm from b] {};
				\node (cin) [v:ghost,position=270:2.25mm from c] {};
				\node (f) [v:main,position=270:13mm from b] {};
				\node (fin) [v:ghost,position=270:2.25mm from f] {};
				\node (i) [v:main,position=270:13mm from e] {};
				\node (iin) [v:ghost,position=90:2.25mm from i] {};
				
				\node (d) [v:main,position=0:13mm from c,fill=white] {};
				\node (g) [v:main,position=270:13mm from c,fill=white] {};
				\node (j) [v:main,position=0:13mm from i,fill=white] {};
				\node (m) [v:main,position=270:13mm from i,fill=white] {};
				
				\node (h) [v:main,position=0:13mm from g] {};
				\node (hin) [v:ghost,position=90:2.25mm from h] {};
				\node (k) [v:main,position=270:13mm from g] {};
				\node (kin) [v:ghost,position=0:2.25mm from k] {};
				\node (n) [v:main,position=270:13mm from j] {};
				\node (nin) [v:ghost,position=180:2.25mm from n] {};
				
				\node (l) [v:main,position=0:13mm from k,fill=white] {};
				\node (o) [v:main,position=270:13mm from k,fill=white] {};
				
				\node (p) [v:main,position=0:13mm from o] {};
				\node (pin) [v:ghost,position=180:2.25mm from p] {};
				
				

				
				
				\node (ab) [v:main,fill=BostonUniversityRed,position=270:4mm from C3] {};
				\node (ei) [v:main,fill=BostonUniversityRed,position=247.5:16mm from ab] {};
				\node (fj) [v:main,fill=BostonUniversityRed,position=292.5:16mm from ab] {};
				\node (mn) [v:main,fill=BostonUniversityRed,position=270:32mm from ab] {};
				\node (cg) [v:main,fill=BostonUniversityRed,position=45:15mm from fj] {};
				\node (dh) [v:main,fill=BostonUniversityRed,position=0:15mm from cg] {};
				\node (kl) [v:main,fill=BostonUniversityRed,position=300:15mm from cg] {};
				\node (op) [v:main,fill=BostonUniversityRed,position=270:14.5mm from kl] {};
				
				

				
				
				\draw (b) [e:main] to (c);
				\draw (b) [e:main] to (f);
				
				\draw (e) [e:main] to (a);
				\draw (e) [e:main] to (f);
				
				\draw (d) [e:main] to (c);
				
				\draw (g) [e:main] to (f);
				\draw (g) [e:main] to (h);
				\draw (g) [e:main] to (k);
				
				\draw (j) [e:main] to (i);
				\draw (j) [e:main] to (k);
				\draw (j) [e:main] to (n);
				
				\draw (m) [e:main] to (i);
				
				\draw (l) [e:main] to (h);
				\draw (l) [e:main] to (p);
				
				\draw (o) [e:main] to (k);
				\draw (o) [e:main] to (n);

				
				
				
				
				\draw (ab) [e:main,->] to (fj);
				\draw (ab) [e:main,->] to (cg);
				
				\draw (ei) [e:main,->,bend left=15] to (fj);
				\draw (ei) [e:main,->] to (ab);
				
				\draw (fj) [e:main,->,bend left=15] to (ei);
				\draw (fj) [e:main,->] to (mn);
				\draw (fj) [e:main,->] to (kl);
				
				\draw (mn) [e:main,->] to (ei);
				
				\draw (cg) [e:main,->,bend right=15] to (dh);
				\draw (cg) [e:main,->] to (fj);
				\draw (cg) [e:main,->] to (kl);
				
				\draw (dh) [e:main,->,bend right=15] to (cg);
				
				\draw (kl) [e:main,->] to (dh);
				\draw (kl) [e:main,->,bend left=15] to (op);
				
				\draw (op) [e:main,->] to (mn);
				\draw (op) [e:main,->,bend left=15] to (kl);
				
				
				
			\end{pgfonlayer}
			

			\begin{pgfonlayer}{background}
				
				\draw (b) [e:coloredborder] to (a);
				\draw (e) [e:coloredborder] to (i);
				\draw (d) [e:coloredborder] to (h);
				\draw (g) [e:coloredborder] to (c);
				\draw (j) [e:coloredborder] to (f);
				\draw (m) [e:coloredborder] to (n);
				\draw (l) [e:coloredborder] to (k);
				\draw (o) [e:coloredborder] to (p);
				
				\draw (b) [e:colored,color=BostonUniversityRed] to (a);
				\draw (e) [e:colored,color=BostonUniversityRed] to (i);
				\draw (d) [e:colored,color=BostonUniversityRed] to (h);
				\draw (g) [e:colored,color=BostonUniversityRed] to (c);
				\draw (j) [e:colored,color=BostonUniversityRed] to (f);
				\draw (m) [e:colored,color=BostonUniversityRed] to (n);
				\draw (l) [e:colored,color=BostonUniversityRed] to (k);
				\draw (o) [e:colored,color=BostonUniversityRed] to (p);
				
			\end{pgfonlayer}	
			
			\begin{pgfonlayer}{foreground}

			\end{pgfonlayer}
		\end{tikzpicture}
	\end{center}
	\caption{Left: A bipartite graph $B$ with a perfect matching $M$. Right: The arising $M$-direction $\DirM{B}{M}$.}
	\label{fig:Mdirection}
\end{figure}

Several properties of matching covered bipartite graphs naturally correspond to properties of digraphs.
In particular this is the case for strong connectivity, as one can easily observe that the $M$-alternating cycles of a bipartite graph $B$ with a perfect matching $M$ are in bijection with the directed cycles of its $M$-direction.
The following statement is folklore (a proof can be found in \cite{zhang2010bipartite}, but the result was already known by \cite{robertson1999permanents}).

\begin{theorem}\label{thm:exttoconn}
	Let $B$ be a bipartite graph with a perfect matching $M$ and $k\in\N$ be a positive integer.
	Then $B$ is \hyperref[def:extendibility]{$k$-extendable} if and only if $\DirM{B}{M}$ is strongly $k$-connected.
\end{theorem}

Even more important, and part of the main motivation behind this project is the following relation between butterfly minors and \hyperref[def:matchingminor]{matching minors}.

\begin{lemma}[\cite{mccuaig2000even}]\label{lemma:mcguigmatminors}
	Let $B$ and $H$ be bipartite matching covered graphs.
	Then $H$ is a \hyperref[def:matchingminor]{matching minor} of $B$ if and only if there exist perfect matchings $M\in\Perf{B}$ and $M'\in\Perf{H}$ such that $\DirM{H}{M'}$ is a \hyperref[def:butterflyminor]{butterfly minor} of $\DirM{G}{M}$.
\end{lemma}

The $M$-direction and the relation between matching minors and butterfly minors allows us to state an interesting set of equivalences regarding bipartite Pfaffian graphs.

\begin{theorem}[\cite{little1975characterization,mccuaig2004polya}]\label{thm:pfaffian}
	Let $B$ be a bipartite graph with a perfect matching $M$.
	The following statements are equivalent.
	\begin{enumerate}
		\item $B$ is Pfaffian.
		
		\item $B$ does not contain $K_{3,3}$ as a matching minor.
		
		\item $\DirM{B}{M}$ does not contain a $\Bidirected{C}_{2k+1}$ for any positive $k\in\N$ as a butterfly minor.
	\end{enumerate}
\end{theorem}

It is worth to note here that in the setting of digraphs an infinitely large family of obstructions must be excluded to achieve the same result as simply excluding $K_{3,3}$ as a matching minor in the bipartite case.
We will revisit this briefly in \cref{sec:consequences}.

\section{The Matching Theoretic Flat Wall}\label{sec:flatwall}

In \cite{hatzel2019cyclewidth} it was proven that any bipartite graph with large enough perfect matching width contains a huge cylindrical wall as a \hyperref[def:matchingminor]{matching minor}.
Since the notion of perfect matching width is not needed for the results in this paper we will not introduce it, just note that the main theorem of this article represents the case of bipartite graphs without $K_{t,t}$ as a matching minor whose perfect matching width exceeds a certain threshold which only depends on $t$.
Still we need a formal introduction of the matching grid itself.

\begin{figure}[!h]
	\centering
	\begin{tikzpicture}

		\pgfdeclarelayer{background}
		\pgfdeclarelayer{foreground}
		\pgfsetlayers{background,main,foreground}

		\draw[e:main] (0,0) circle (11mm);
		\draw[e:main] (0,0) circle (16mm);
		\draw[e:main] (0,0) circle (21mm);
		\draw[e:main] (0,0) circle (26mm);
		
		\foreach \x in {1,...,4}
		{
			\draw[e:main] (\x*90:16mm) -- (\x*90+22.5:11mm);
			\draw[e:main] (\x*90:21mm) -- (\x*90+22.5:16mm);
			\draw[e:main] (\x*90:26mm) -- (\x*90+22.5:21mm);
		}
		
		\foreach \x in {1,...,4}
		{
			\draw[e:main] (\x*90-22.5:16mm) -- (\x*90-45:11mm);
			\draw[e:main] (\x*90-22.5:21mm) -- (\x*90-45:16mm);
			\draw[e:main] (\x*90-22.5:26mm) -- (\x*90-45:21mm);
			
		}
		
		\foreach \x in {1,...,8}
		{
			\draw[e:coloredthin,color=BostonUniversityRed,bend right=13] (\x*45:11mm) to (\x*45+22.5:11mm);
			\draw[e:coloredthin,color=BostonUniversityRed,bend right=13] (\x*45:16mm) to (\x*45+22.5:16mm);
			\draw[e:coloredthin,color=BostonUniversityRed,bend right=13] (\x*45:21mm) to (\x*45+22.5:21mm);
			\draw[e:coloredthin,color=BostonUniversityRed,bend right=13] (\x*45:26mm) to (\x*45+22.5:26mm);
		}
		
		\foreach \x in {1,...,8}
		{
			\node[v:main] () at (\x*45:11mm){};
			\node[v:main] () at (\x*45:16mm){};
			\node[v:main] () at (\x*45:21mm){};
			\node[v:main] () at (\x*45:26mm){};
			\node[v:mainempty] () at (\x*45+22.5:11mm){};
			\node[v:mainempty] () at (\x*45+22.5:16mm){};
			\node[v:mainempty] () at (\x*45+22.5:21mm){};
			\node[v:mainempty] () at (\x*45+22.5:26mm){};
		}

		\begin{pgfonlayer}{background}
			\foreach \x in {1,...,8}
			{
				\draw[e:coloredborder,bend right=13] (\x*45:11mm) to (\x*45+22.5:11mm);
				\draw[e:coloredborder,bend right=13] (\x*45:16mm) to (\x*45+22.5:16mm);
				\draw[e:coloredborder,bend right=13] (\x*45:21mm) to (\x*45+22.5:21mm);
				\draw[e:coloredborder,bend right=13] (\x*45:26mm) to (\x*45+22.5:26mm);
			}
		\end{pgfonlayer}
		
	\end{tikzpicture}
	\caption{The \hyperref[def:matchinggrid]{cylindrical matching grid} of order $4$ with the canonical matching.}
	\label{fig:cylindricalgrid}
\end{figure}
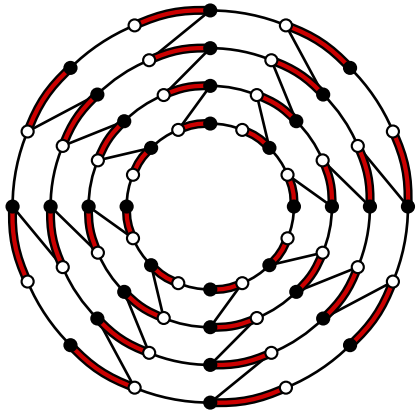

\begin{definition}[Cylindrical Matching Grid]\label{def:matchinggrid}
	The \emph{cylindrical matching grid} $CG_k$ of order $k$ is defined as follows.
	Let $C_1,\dots,C_k$ be $k$ vertex disjoint cycles of length $4k$.
	For every $i\in[1,k]$ let $C_i=\Brace{v_1^i,v_2^i,\dots,v_{4k}^i}$, $V_1^i\coloneqq\CondSet{v_j^i}{j\in\Set{1,3,5,\dots,4k-1}}$, $V_2^i\coloneqq\Fkt{V}{C_i}\setminus V_1^i$, and $M_i\coloneqq\CondSet{v_j^iv_{j+1}^i}{v_j^i\in V_1^i}$.
	Then $CG_k$ is the graph obtained from the union of the $C_i$ by adding
	\begin{align*}
		\CondSet{v_j^iv_{j+1}^{i+1}}{i\in[1,k-1]~\text{and}~j\in\Set{1,5,9,\dots,4k-3}}&\text{, and}\\
		\CondSet{v_j^iv_{j+1}^{i-1}}{i\in[2,k]~\text{and}~j\in\Set{3,7,11,\dots,4k-1}}&
	\end{align*}
	to the edge set.
	We call $M\coloneqq\bigcup_{i=1}^kM_i$ the \emph{canonical matching} of $CG_k$.
	See \cref{fig:cylindricalgrid} for an illustration.
\end{definition}

Please note that the cylindrical matching grid is indeed a subcubic graph and thus, by \cref{lemma:confmathingminors}, one can always find a conformal bisubdivision of $CG_k$ within a graph $G$ if it contains $CG_k$ as a matching minor.
To allows os for an easier transition between the setting of bipartite graphs with perfect matchings and digraphs we do not use the cylindrical matching grid itself as our wall, but a slight modification.

\begin{definition}[Matching Wall]\label{def:matchingwall}
	Let $k\in\N$ be a positive integer.
	The \emph{elementary matching $k$-wall} with its \emph{canonical matching} $M$ is the graph obtained from the \hyperref[def:matchinggrid]{cylindrical matching grid} of order $2k$ with canonical matching $M$ by deleting the non-$M$-edges whose index in congruent to $1$ modulo $4$ on each of the concentric cycles with even index and the non-$M$-edges whose index in congruent to $3$ modulo $4$ on each of the concentric cycles with odd index.
	A \emph{matching $k$-wall} $W'$ is a bisubdivision of the elementary matching $k$-wall, and a perfect matching $M'$ is its \emph{canonical matching} if $\DirM{W'}{M'}$ is a cylindrical $k$-wall.
	
	The \emph{perimeter} of $W'$, denoted by $\Perimeter{W'}$, is the union of the outer-most and the inner-most $M'$-conformal cycle of $W'$.
\end{definition}

Finally, to define what it means for a matching wall to `grasp' a \hyperref[def:matchingminor]{matching minor} we need to introduce the notion of `models' for matching minors as well.

Models, or embeddings, for matching minors have already been used and discussed in \cite{robertson1999permanents} and \cite{norine2007generating} and the definitions we give here are similar to those of Norine and Thomas.
Some parts of these definitions, however, have been changed to better suit our needs in the sections to come and therefore we provide the necessary proofs.

Let $T'$ be a tree and let $T$ be obtained from $T'$ by subdividing every edge an odd number of times.
Then $\Fkt{V}{T'} \subseteq \Fkt{V}{T}$.
The vertices of $T$ that belong to $T'$ are called \emph{old}, and the vertices in $\Fkt{V}{T} \setminus \Fkt{V}{T'}$ are called \emph{new}.
We say that $T$ is a \emph{barycentric tree}.

Let $G$ and $H$ be graphs with perfect matchings.
An \emph{embedding} or \emph{matching minor model} of $H$ in $G$ is a mapping $\mu \colon \Fkt{V}{H} \cup \Fkt{E}{H} \to \CondSet{F}{F\subseteq G}$,	such that the following requirements are met for all $v,v' \in \Fkt{V}{H}$ and $e,e' \in \Fkt{E}{H}$:
\begin{enumerate}
	\item $\Fkt{\mu}{v}$ is a barycentric subtree in $G$,
	
	\item if $v \neq v'$, then $\Fkt{\mu}{v}$ and $\Fkt{\mu}{v'}$ are vertex disjoint,
	
	\item $\Fkt{\mu}{e}$ is an odd path with no internal vertex in any $\Fkt{\mu}{v}$, and if $e' \neq e$, then $\Fkt{\mu}{e}$ and $\Fkt{\mu}{e'}$ are internally vertex disjoint,
	
	\item if $e=u_1u_2$, then the ends of $\Fkt{\mu}{e}$ can be labelled by $x_1,x_2$ such that $x_i$ is an old vertex of $\Fkt{\mu}{u_i}$,
	
	\item if $v$ has degree one, then $\Fkt{\mu}{v}$ is exactly one vertex, and
	
	\item $G-\Fkt{\mu}{H}$ has a perfect matching, where $\Fkt{\mu}{H'} \coloneqq \bigcup_{x\in \Fkt{V}{H'} \cup \Fkt{E}{H'}}\Fkt{\mu}{x}$ for every subgraph $H'$ of $H$.
\end{enumerate}	
If $\mu$ is a matching minor model of $H$ in $G$ we write $\mu\colon H\rightarrow G$.
Let $M$ be a perfect matching such that $\Fkt{\mu}{H}$ is $M$ conformal.
We call a matching minor model $\mu\colon H\rightarrow G$ an \emph{$M$-model} of $H$ in $G$ if $\Fkt{\mu}{H}$ is $M$-conformal.

\begin{lemma}[\cite{giannopoulou2021excluding}]\label{cor:Mmodels}
	Let $G$ and $H$ be graphs with perfect matchings and $M$ a perfect matching of $G$.
	Then $H$ is isomorphic to an \hyperref[def:matchingminor]{$M$-minor} of $G$ if and only if there exists an $M$-model $\mu_M\colon H\rightarrow G$ in $G$.
\end{lemma}

Let $k,t\in\N$ be positive integers and $B$ be a bipartite graph with a conformal matching $k$-wall $W$, and let $H$ be some bipartite matching covered graph.
Let $M$ be a perfect matching of $B$.
We say that $W$ \emph{$M$-grasps} an $H$-matching minor if there exists a matching minor model $\mu\colon H\rightarrow B$ such that $\Fkt{\mu}{H}$ is $M$-conformal, and for every $e\in\E{H}\cap\Remainder{M}{H}$ the $M$-conformal path $\Fkt{\mu}{e}$ is completely contained in $W$.
We say that $W$ \emph{grasps} an $H$-matching minor if there exists a perfect matching $M$ of $B$ such that $W$ $M$-grasps $H$.

et $B$ be a Pfaffian brace and $H$ be a planar brace.
We say that $H$ is a \emph{summand} of $B$ if there exist planar braces $H_1,\dots,H_{\ell}$ such that $B$ can be constructed from the $H_i$ by repeated applications of the trisum operation, and $H=H_1$.

What is left to do before we can formally state our main theorem is a definition of flatness in the matching theoretic sense.
In contrast to the undirected setting, where it makes sense to speak about subgraphs as a means of reductions, in the setting of graphs with perfect matchings, we sometimes will have to perform tight cut contractions.

Let $B$ and $H$ be bipartite graphs with a perfect matching such that $H$ has a single brace $J$ that is not isomorphic to $C_4$.
We say that $H$ is a \emph{$J$-expansion}.
A brace $G$ of $B$ is said to be a \emph{host of $H$} if $G$ contains a conformal subgraph $H'$ that is a $J$-expansion and can be obtained from $H$ by repeated applications of tight cut contractions.
The graph $H'$ is called the \emph{remnant of $H$}.

It makes sense for us to work with a cylindrical grid/wall rather than a square one.
However, the only problem this brings is that any cylindrical wall has two faces which might be considered the natural `outer face'.
Indeed, we would like the inner-most and the outer-most cycle of the cylindrical wall to both bound faces in an appropriate reduction.

Let $B$, $H$, and $J$ be bipartite graphs with perfect matchings such that $H$ and $J$ are conformal subgraphs of $B$.
We say that $H$ is \emph{$J$-bound} if there exists a subgraph $K$ of $B-J$ that is the union of elementary components of $B-J$ such that $K\cup J$ is matching covered, and $H$ is a conformal subgraph of $K\cup J$.
The graph $K\cup J$ is called a \emph{$J$-base of $H$} in $B$.

\begin{definition}[$P$-Flatness]\label{def:matchingflat}
	Let $B$ be a bipartite graph with a perfect matching, and let $H$ be a planar matching covered graph that is a $J$-expansion of some planar brace $J$.
	Moreover, let $P$ be a collection of pairwise vertex disjoint faces of $H$ such that $P$ is a conformal subgraph of $H$.
	At last, let $A\subseteq\V{B}$ be a conformal set.
	Then $H$ is \emph{$P$-flat} in $B$ \emph{with respect to $A$} if 
	\begin{enumerate}
		\item $H$ is a conformal subgraph of $B'\coloneqq B-A$,
		\item some $P$-base of $H$ in $B'$ has a Pfaffian brace $B''$ that is a host of $H$, and
		\item $B''$ has a summand $G$ that contains a remnant $H'$ of $H$ such that every remnant of a face from $P$ within $H'$ bounds a face of $G$.
	\end{enumerate}
\end{definition} 

The set $A$ in the definition above is the \emph{apex set}, similar to the one which occurs in the version of flatness used for the original Flat Wall Theorem.
The set $P$ takes on two roles at once, it mimics the separator of the separation $\Brace{X,Y}$ in the original definition and also allows us to essentially, prescribe which faces of $H$ should take the role of the outer face.
Since we do not require $B-A$ to be a brace or even matching covered, we need to remove all non-admissible edges and just take the matching covered subgraph that contains $H$.
This is modelled by selecting a certain $P$-base of $H$ in $B'$.
Now that we have reduced $B-A$ to a matching covered graph $B'''$, we must go one step further and get rid of the non-trivial tight cuts of $B'''$.
This is done by selecting $B''$ to be a host of $H$ in $B'''$.
By requiring $B''$ to be Pfaffian, we ready ourselves for the final reduction.
Indeed, since we insist $G$ to be a summand of $B''$ this means $B''$ cannot be isomorphic to the Heawood graph by \cref{thm:trisums}.
Since $G$ is a summand of $B''$, it must be a planar brace.
To get to this point, some tight cut contractions could have been necessary, and thus we are only able to talk about a remnant $H'$ of $H$, but since $H$ was chosen to be a $J$-expansion of some planar brace $J$, this remnant is well defined.
Similarly, the tight cut contractions could have shrunken some of the faces that were selected to form $P$, but we can still make out their remnants and thus (iii) resembles the third requirement of the original definition.
With this all necessary definitions for \cref{thm:matchingflatwall} are in place.
Please note that the matching minor from \cref{thm:matchingflatwall} is not necessarily conformal for the canonical perfect matching of the respective wall.

In the following, we sometimes say that a conformal matching $k$-wall $W$ is \emph{flat} in $B$, if $k$ is large enough and the second part of the theorem above holds true for $W$.

With \cref{thm:matchingflatwall} at hand, we can give an approximate characterisation of all bipartite graphs with perfect matchings that exclude $K_{t,t}$ as a matching minor for some $t\in\N$.
This weak structure theorem is similar to \cref{thm:undirectedflatwall} and in some sense can be seen as a generalisation of \cref{thm:trisums}.

\begin{theorem}\label{thm:weakmatchingstructure}
	Let $r,t\in\N$ be positive integers, $\ApexBound$ and $\WallBound$ be the two functions from \cref{thm:matchingflatwall}, and $B$ be a bipartite graph with a perfect matching.
	\begin{itemize}
		\item If $B$ has no $K_{t,t}$-\hyperref[def:matchingminor]{matching minor}, then for every conformal \hyperref[def:matchingwall]{matching $\Fkt{\WallBound}{t,r}$-wall} $W$ in $B$ and every perfect matching $M$ of $B$ such that $M\cap\E{W}$ is the canonical matching of $W$, there exist an $M$-conformal set $A\subseteq\V{B}$ with $\Abs{A}\leq\Fkt{\ApexBound}{t}$ and an $M$-conformal matching $r$-wall $W'\subseteq W-A$ such that $W'$ is \hyperref[def:matchingflat]{$\Perimeter{W'}$-flat} in $B$ with respect to $A$.
		
		\item Conversely, if $t\geq 2$ and $r\geq\sqrt{\Fkt{2\ApexBound}{t}}$, and for every conformal \hyperref[def:matchingwall]{matching $\Fkt{\WallBound}{r,t}$-wall} $W$ in $B$ and every perfect matching $M$ of $B$ such that $M\cap\E{W}$ is the canonical perfect matching of $W$, there exist an $M$-conformal set $A\subseteq\V{B}$ with $\Abs{A}\leq\Fkt{\ApexBound}{t}$ and an $M$-conformal matching $r$-wall $W'\subseteq W-A$ such that $W'$ is \hyperref[def:matchingflat]{$\Perimeter{W'}$-flat} in $B$ with respect to $A$, then $B$ has no \hyperref[def:matchingminor]{matching minor} isomorphic to $K_{t',t'}$, where $t'=16\Fkt{\WallBound}{t,r}^2$. 
	\end{itemize}
\end{theorem}

\begin{proof}
	The first part of the theorem follows immediately from \cref{thm:matchingflatwall}, since in case $B$ does not have $K_{t,t}$ as a matching minor, the first part of \cref{thm:matchingflatwall} can never be true and thus every matching $\Fkt{\WallBound}{t,r}$-wall must be flat in $B$.
	
	For the reverse, note that an elementary matching $\Fkt{\WallBound}{t,r}$-wall has exactly $16\Fkt{\WallBound}{t,r}^2$ vertices.
	Now suppose $B$ has a matching minor model $\mu\colon K_{t',t'}\rightarrow B$.
	Then there exists a perfect matching $M$ such that $\mu$ is $M$-conformal.
	Indeed, $K_{t',t'}$ contains an $\Remainder{M}{K_{t',t'}}$-conformal elementary matching $\Fkt{\WallBound}{t,r}$-wall, and thus $\Fkt{\mu}{K_{t',t'}}$ contains an $M$-conformal matching $\Fkt{\WallBound}{t,r}$-wall $W$.
	For every vertex $w$ of degree three in $W$ there exists a unique vertex $u_w\in\V{K_{t',t'}}$ such that $w\in\V{\Fkt{\mu}{u_w}}$, and in case $w\neq w'$ are both vertices of degree three in $W$, then $u_w\neq u_{w'}$.
	Moreover, if $P$ is a path in $W$ whose endpoints $w$ and $w'$ have degree three in $W$ and all internal vertices are vertices of degree two in $W$, then $\V{P}\subseteq \V{\Fkt{\mu}{u_w}}\cup\V{\Fkt{\mu}{u_{w'}}}$.
	
	By assumption there exist an $M$-conformal set $A\subseteq\V{B}$ and an $M$-conformal matching $r$-wall $W'\subseteq W$ such that $W'$ is \hyperref[def:matchingflat]{$\Perimeter{W'}$-flat} in $B$ with respect to $A$.
	Now $W'$ has $16r^2$ many vertices of degree three in $W'$, $16r$ of which lie on $\Perimeter{W'}$.
	Since $r\geq \sqrt{2\Fkt{\ApexBound}{t}}$, we have at least $32\Fkt{\ApexBound}{t}$ many such degree three vertices.
	Thus, with $\Abs{A}\leq \Fkt{\ApexBound}{t}$ and $t\geq 2$, there exist $w_1,\dots,w_6\in\V{W'-\Perimeter{W'}}$ such that $\V{\Fkt{\mu}{u_{w_i}}}\cap A=\emptyset$ for all $i\in[1,6]$.
	This, however, means that for every $\Perimeter{W'}$-base $H$ of $W'$, every brace $J$ of $H$ that is a host of $W'$ must contain $K_{3,3}$ as a \hyperref[def:matchingminor]{matching minor} and therefore no such $J$ can be Pfaffian by \cref{thm:pfaffian}.
	Hence $W'$ cannot be $\Perimeter{W'}$-flat in $B$ with respect to $A$ and we have reached a contradiction.
\end{proof}

\section{Consequences for (Di)Graphs}\label{sec:consequences}

In \cref{subsec:digraphsandmatchings} we already started to discuss the deeper connections between the structural theories of bipartite graphs with perfect matchings and digraphs.
In this section we utilise \cref{thm:matchingflatwall} to derive a flat wall theorem for digraphs which is different from \cref{thm:directedflatwall} in the sense that we exclude an entire infinite anti-chain rather than just $\Bidirected{K_t}$ as a butterfly minor but in return we obtain a much more structured result and a more comprehensible version of flatness.
This alternate directed flat wall theorem has a lot in common with \cref{thm:pfaffian}.
In \cref{sec:thomassen} we then apply the new directed flat wall theorem to obtain a completely new duality theorem for treewidth in general graphs.

For this we first need to introduce some additional definitions regarding infinite anti-chains of butterfly minors.
The groundwork for the theory necessary has been established in \cite{giannopoulou2021excluding}.
Recall the definition of \hyperref[def:Mdirection]{$M$-directions} of bipartite graphs and that we argued that the operation can be reversed to yield a unique bipartite graph with a perfect matching when given a digraph $D$ as input.
Let us formalise this inverse operation.

\begin{definition}[Split]\label{def:split}
	Let $D$ be a directed graph.
	We define $\Split{D}$ to be the bipartite graph $B$ for which a perfect matching $M$ exists such that $\DirM{B}{M}=D$.
\end{definition}

The digraph $J$ is a \emph{proper butterfly minor} of the digraph $D$ if $J$ is a butterfly minor of $D$ and $J\not\cong D$.
We say that $D$ is \emph{$J$-minimal} if $\Split{D}$ contains $\Split{J}$ as a \hyperref[def:matchingminor]{matching minor}, but for every proper butterfly minor $D'$ of $D$, $\Split{D'}$ is $\Split{J}$-free.

\begin{definition}[Fundamental Anti-Chain]
	Let $D$ be a digraph.
	The family
	\begin{align*}
		\Antichain{D}\coloneqq\CondSet{D'}{\text{$D'$ is a $D$-minimal digraph}}
	\end{align*}
	is called the \emph{fundamental anti-chain based on $D$}.
\end{definition}

\begin{lemma}[\cite{giannopoulou2021excluding}]\label{lemma:excludingantichains}
	Let $H$ and $D$ be digraphs.
	Then $D$ contains a butterfly minor from $\Antichain{H}$ if and only if $\Split{D}$ contains $\Split{H}$ as a matching minor.
\end{lemma}

Please note that the set $\CondSet{\Bidirected{C}_{2k+1}}{k\geq1,~k\in\N}$ is exactly the anti-chain $\Antichain{\Bidirected{K_3}}$.

\subsection{Another Directed Flat Wall Theorem}\label{sec:directedflatwall}

An advantage of relying heavily on the digraphic setting for the proof of \cref{thm:matchingflatwall} is that we may first select a perfect matching $M$ of our choosing and then just have to consider those $M$-conformal \hyperref[def:matchingwall]{matching walls} which have $M$ as their canonical matching, and only $M$-conformal sets as our apex set.
This means, we can directly translate our matching version of the flat wall theorem into the digraphic setting.
For this to work however, we need to translate our matching theoretic notion of flatness into the world of digraphs.
As this notion is tied to the structural description of bipartite $K_{3,3}$-\hyperref[def:matchingminor]{matching minor} free graphs, we may use a structural characterisation of $\Antichain{\Bidirected{K_3}}$-butterfly minor free digraphs from \cite{robertson1999permanents}.

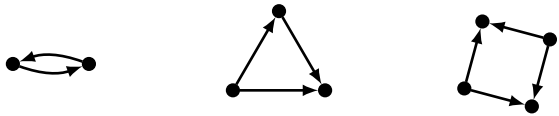
\begin{figure}[!ht]
	\centering
	\begin{tikzpicture}
		\pgfdeclarelayer{background}
		\pgfdeclarelayer{foreground}
		\pgfsetlayers{background,main,foreground}
		
		\node (mid) [v:ghost] {};
		\node (left) [v:ghost,position=180:30mm from mid] {};
		\node (right) [v:ghost,position=0:30mm from mid] {};

		\node(u1) [v:main,position=180:5mm from left] {};
		\node(u2) [v:main,position=0:5mm from left] {};
		
		\draw [e:main,->,bend right=20] (u1) to (u2);
		\draw [e:main,->,bend right=20] (u2) to (u1);
		
		\node (v1) [v:main,position=90:7mm from mid] {};
		\node (v2) [v:main,position=210:7mm from mid] {};
		\node (v3) [v:main,position=330:7mm from mid] {};

		\draw [e:main,<-] (v1) to (v2);
		\draw [e:main,->] (v2) to (v3);
		\draw [e:main,<-] (v3) to (v1);
		
		\node (w1) [v:main,position=30:6.5mm from right] {};
		\node (w2) [v:main,position=120:6.5mm from right] {};
		\node (w3) [v:main,position=210:6.5mm from right] {};
		\node (w4) [v:main,position=300:6.5mm from right] {};
		
		\draw [e:main,->] (w1) to (w2);
		\draw [e:main,->] (w1) to (w4);
		\draw [e:main,->] (w3) to (w2);
		\draw [e:main,->] (w3) to (w4);
		
	\end{tikzpicture}
	\caption{The subgraphs necessary for the small-cycle-sum operation.}
	\label{fig:smallcyclesum}
\end{figure}

\begin{definition}[Small-Cycle-Sum]\label{def:smallcyclesum}
	Let $D_0$ be a digraph, let $u,v\in\V{D_0}$, and let $(u,v),(v,u)\in\E{D_0}$.
	Let $D_1$ and $D_2$ be such that $D_1\cup D_2=D_0$, $\V{D_1}\cap\V{D_2}=\Set{u,v}$, $\V{D_1}\setminus\V{D_2}\neq\emptyset$, $\V{D_2}\setminus\V{D_1}\neq\emptyset$, and $\E{D_1}\cap\E{D_2}$.
	Let $D$ be obtained from $D_0$ by deleting some (possibly neither) of the edges $(u,v)$, $(v,u)$.
	We say that $D$ is a \emph{$2$-sum} of $D_1$ and $D_2$.
	
	Let $D_0$ be a digraph, let $u,v,w\in\V{D_0}$ and $(u,v),(w,v),(w,u)\in\E{D_0}$, and assume that $D_0$ has a directed cycle containing the edge $(w,v)$, but not the vertex $u$.
	Let $D_1$ and $D_2'$ be such that $D_1\cup D_2'\neq D_0$, $\V{D_1}\cap\V{D_2'}=\Set{u,v,w}$, $\V{D_1}\setminus\V{D_2'}\neq\emptyset$, $\V{D_2'}\setminus\V{D_1}\neq\emptyset$, and $\E{D_1}\cap\E{D_2'}=\Set{(u,v),(w,v),(w,u)}$.
	Let $D_2'$ have no edge with tail $v$, and no edge with head $w$ and note that this means that $(w,v)$ is butterfly contractible in $D_2'$, let $D_2$ be the digraph obtained from $D_2'$ by contracting $(w,v)$.
	Let $D$ be obtained from $D_0$ by deleting some (possible none) of the edges $(u,v),(w,v),(w,u)$.
	We say that $D$ is a \emph{$3$-sum} of $D_1$ and $D_2$.
	
	Let $D_0$ be a digraph, let $x,y,u,v\in\V{D_0}$ as well as $(x,y),(x,v),(u,y),(u,v)$, and assume that $D_0$ has a directed cycle containing precisely two of the edges $(x,y),(x,v),(u,y),(u,v)$.
	Let $D_1$ and $D_2'$ be such that $D_1\cup D_2'=D_0$, $\V{D_1}\cap\V{D_2'}=\Set{x,y,u,v}$, $\V{D_1}\setminus\V{D_2'}\neq\emptyset$, $\V{D_2'}\setminus\V{D_1}\neq\emptyset$, and $\E{D_1}\cap\E{D_2'}=\Set{(x,y),(x,v),(u,y),(u,v)}$.
	Let $D_2'$ have no edge with tail $y$ or $v$, and no edge with head $x$ or $u$ and note that this means that the edges $(x,y)$ and $(u,v)$ are butterfly contractible.
	Let $D_2$ be the digraph obtained from $D_2'$ by contracting the edges $(x,y)$ and $(u,v)$.
	Finally, let $D$ be obtained from $D_0$ by deleting some (possible none) of the edges $(x,y),(x,v),(u,y),(u,v)$.
	We say that $D$ is a \emph{$4$-sum} of $D_1$ and $D_2$.
	
	We say that a digraph $D$ is a \emph{small cycle sum} of two digraphs $D_1$ and $D_2$ if it is an $i$-sum of $D_1$ and $D_2$ for some $i\in[1,3]$.
\end{definition}

\begin{definition}[Strong Genus]\label{def:stronggenus}
Let $D$ be a digraph and $\Sigma$ be a surface.
An embedding $\mu\colon D\rightarrow\Sigma$ of $D$ into $\Sigma$ is \emph{strong} if for every vertex $v\in\V{D}$ there exist $r_v\in\R$ and an open disc $\zeta\subseteq\Sigma$ of radius $r_v$ centred at $v$ such that $\Brace{\CondSet{\Brace{u,v}}{\Brace{u,v}\in\E{D}},\CondSet{\Brace{v,u}}{\Brace{v,u}\in\E{D}}}$ is a butterfly in $\zeta$.
The smallest integer $h\in\N$ such that $D$ can be strongly embedded in $\Sigma_{\frac{h}{2}}$ or $\tilde{\Sigma}_h$ is called the \emph{strong genus} of $D$.
We denote the strong genus of $D$ by $\StrongGenus{D}$.
If $\StrongGenus{D}=0$, $D$ is said to be \emph{strongly planar}.
\end{definition}

Finally, let $H_{14}$ be the \hyperref[fig:heawood]{Heawood graph} and let $M$ be any perfect matching of $H_{14}$.
We denote the digraph $\DirM{H_{14}}{M}$ by $F_7$.

\begin{theorem}[\cite{robertson1999permanents}]\label{thm:nonevendigraphs}
	Let $D$ be a strongly $2$-connected digraph.
	Then $D$ is does not contain a member of $\Antichain{\Bidirected{K_3}}$ as a butterfly minor if and only if it can be obtained from a family of strongly $2$-connected strongly planar digraphs and $F_7$ by repeated applications of the small-cycle-sum operation.
\end{theorem}

For our matching theoretic definition of flatness we also made use of tight cuts to reduce the matching covered bipartite graphs involved to braces.
A similar thing needs to happen here as well for us to be able to make use of \cref{thm:nonevendigraphs}.
Recall \cref{thm:braces,thm:exttoconn} and observe that this means that the \hyperref[def:split]{split} of a strongly connected digraph $D$ has a non-trivial \hyperref[def:braces]{tight cut} if and only if $D$ itself has a cut vertex, that means $D$ has a non-trivial directed separation of order one.
It is a known fact that the non-trivial directed separations of a digraph $D$ and the non-trivial tight cuts of its split are in correspondence and thus we may define a directed version of the tight cut decomposition alongside the notion of `dibraces'.

\begin{definition}[Directed Tight Cut Contraction]\label{def:directedtightcut}
	Let $D$ be a digraph, and $\Brace{X,Y}$ be a non-trivial directed separation of order $1$ in $D$.
	Let $\Set{v}=X\cap Y$, we set
	\begin{align*}
		\ContractsTo{D}{X}{v_X}\coloneqq D-X+v_X+&\CondSet{\Brace{y,v_X}}{\Brace{y,v}\in\E{D}\text{ and }y\in Y}\\+&\CondSet{\Brace{v_X,y}}{\Brace{x,y}\in\E{D},~x\in X\text{, and }y\in Y}\text{, and}\\
		\ContractsTo{D}{Y}{v_Y}\coloneqq D-Y+v_Y+&\CondSet{\Brace{v_Y,x}}{\Brace{v,x}\in\E{D}\text{ and }x\in X}\\+&\CondSet{\Brace{x,v_Y}}{\Brace{x,y}\in\E{D},~x\in X\text{, and }y\in Y}.
	\end{align*}
	The two digraphs $D_X\coloneqq \ContractsTo{D}{X}{v_X}$ and $D_Y\coloneqq \ContractsTo{D}{Y}{v_Y}$ are called the \emph{$\Brace{X,Y}$-contractions} of $D$.
\end{definition}
Similar to the tight cut decomposition of matching covered graphs we may select any non-trivial directed separation of order one in a given digraph $D$ and form its two contractions.
Then, in case any of the two contractions still has a non-trivial directed separation of order one we can proceed with this process.
Eventually this yields a list of strongly $2$-connected digraphs uniquely determined\footnote{This follows from the famous result of Lovasz \cite{lovasz1987matching} on the uniqueness of the tight cut decomposition.} by $D$.
These strongly $2$-connected digraphs are called the \emph{dibraces} of $D$.

Since we have already seen the equivalence between the contraction of directed separations of order one and strongly $2$-connected digraphs with tight cut contractions and braces, by recalling the small-cycle-sum operation we can deduce analogues to the definitions necessary for $H$-flatness in bipartite graphs with perfect matchings.

Let $D$ be a strongly $2$-connected digraph without a butterfly minor from $\Antichain{\Bidirected{K_3}}$ and let $H$ be a strongly $2$-connected \hyperref[def:stronggenus]{strongly planar} digraph.
We say that $H$ is a \emph{summand} of $D$ if there exist strongly planar strongly $2$-connected digraphs $H_1,\dots,H_{\ell}$ such that $D$ can be constructed from the $H_i$ by means of small-cycle-sums, where $H=H_1$.

Let $D$ and $H$ be digraphs such that $H$ has exactly one dibrace\footnote{Recall that a dibrace of a digraph is a strongly $2$-connected butterfly minor that can be obtained purely by means of directed tight cut contractions.} $J$ which is not isomorphic to $\Bidirected{K_2}$.
We say that $H$ is a \emph{$J$-expansion}.
A dibrace $G$ of $D$ is said to be a \emph{host} of $H$, if $G$ contains a subgraph $H'$ that is a $J$-expansion and can be obtained from $H$ by means of directed tight cut contractions.
The digraph $H'$ is called the \emph{remnant} of $H$.

Let $D$, $H$, and $J$ be digraphs such that $H$ and $J$ are subgraphs of $D$.
We say that $H$ is \emph{$J$-bound} if there exists a subgraph $K$ of $D-J$ that is the union of strong components of $D-J$ such that $K\cup J$ is strongly connected and $H$ is a subgraph of $K\cup J$.
The digraph $K\cup J$ is called a \emph{$J$-base} of $H$ in $D$.

\begin{definition}[$P$-Flatness for Digraphs]
	Let $D$ be a digraph and $H$ be a strongly planar and strongly connected digraph that is a $J$-expansion of some strongly planar and strongly $2$-connected digraph $J$.
	Moreover, let $P$ be a collection of pairwise vertex disjoint faces of $H$ which each are bound by a directed cycle.
	At last, let $A\subseteq\V{D}$ be a set of vertices.
	Then $H$ is \emph{$P$-flat} in $D$ \emph{with respect to $A$} if
	\begin{enumerate}
		\item $H$ is a subgraph of $D'\coloneqq D-A$,
		\item some $P$-base of $H$ in $D'$ has a dibrace $D''$ that excludes all members of $\Antichain{\Bidirected{K_3}}$ as butterfly minors and that is a host of $H$, and
		\item $D''$ has a summand $G$ that contains the remnant $H'$ of $H$, such that every remnant of a face from $P$ bounds a face of $G$.  
	\end{enumerate}
\end{definition}

This allows us to formulate the anti-chain version of the directed flat wall theorem.

\begin{theorem}\label{thm:antichaindirectedflattwall}
	Let $r,t\in\N$ be positive integers.
	There exist functions $\DirectedApexBound\colon\N\rightarrow\N$ and $\DirectedWallBound\colon\N\times\N\rightarrow\N$ such that for every digraph the following is true:
	If $W$ is a cylindrical $\Fkt{\DirectedWallBound}{t,r}$-wall in $D$, then either
	\begin{enumerate}
		\item There exists $H\in\Antichain{\Bidirected{K_t}}$ such that $D$ has an $H$-butterfly minor grasped by $W$, or
		\item there exists a set $A\subseteq\V{D}$ with $\Abs{A}\leq\Fkt{\DirectedApexBound}{t}$ and a cylindrical $r$-wall $W'\subseteq W-A$ such that $W'$ is $\Perimeter{W'}$-flat in $D$ with respect to $A$.
	\end{enumerate}
\end{theorem}

\begin{proof}
	The theorem follows almost immediately from \cref{thm:matchingflatwall}.
	Let $D$ be a digraph, $t,r\in\N$ be positive integers, and $W$ be a cylindrical $\Fkt{\WallBound}{t,r}$-wall in $D$.
	Moreover, consider $B\coloneqq\Split{D}$ together with the perfect matching $M$ for which $\DirM{B}{M}=D$ holds.
	Then $\Split{W}$ is an $M$-conformal \hyperref[def:matchingwall]{matching $\Fkt{\WallBound}{t,r}$-wall} such that $M$ contains the canonical matching of $\Split{W}$.
	So \cref{thm:matchingflatwall} leaves us with two cases.
	In the first case we find a $K_{t,t}$-\hyperref[def:matchingminor]{matching minor} in $B$ grasped by $\Split{W}$.
	This implies the existence of some $H\in\Antichain{\Bidirected{K_t}}$, such that $D$ contains an $H$-butterfly minor grasped by $W$.
	Otherwise there exist an $M$-conformal set $A'\subseteq\V{B}$ of size at most $\Fkt{\ApexBound}{t}$ and an $M$-conformal matching $r$-wall $U'\subseteq\Split{W}$, such that $U'$ is \hyperref[def:matchingflat]{$\Perimeter{U'}$-flat} in $B$ with respect to $A'$.
	Since there is a bijection between the non-trivial \hyperref[def:braces]{tight cuts} in $B$ and the non-trivial directed separations of order one in $D$, we can choose $A$ to be the set of all edges of $M$ with both endpoints in $A'$, and our claim follows.
\end{proof}

Whilst a directed version of \cref{thm:matchingflatwall} can be obtained with a relatively straight forward argument, it is not clear whether one can adapt \cref{thm:weakmatchingstructure} into the digraphic setting.
This is because the reverse direction of \cref{thm:weakmatchingstructure} requires an assumption on all perfect matchings, while in $D$ we, a priori at least, can only talk about the unique perfect matching $M$ which is used to obtain $D$ from $\Split{D}$ as its $M$-direction.
Hence for a weak structure theorem for digraphs similar to \cref{thm:weakmatchingstructure} one first needs to find a way to express the structure of all perfect matching of $\Split{D}$ in a concise way by just considering $D$.
One direction this might go could be expanding \cref{thm:antichaindirectedflattwall} to all members of $\Antichain{W}$, where $W$ is a cylindrical $\Fkt{\DirectedWallBound}{t,r}$-wall.

\subsection{Generalising a Theorem of Thomassen}\label{sec:thomassen}

Before the resolution of the bipartite \textsc{Pfaffian Recognition} by McCuaig et al.\@ only partial results were known.
Among these is a result on symmetric digraphs\footnote{Recall that a digraph is symmetric if it can be obtained from an undirected graph by replacing every undirected edge by a digon.} by Thomassen.
This result shows an interesting feature of symmetric digraphs excluding the members of $\Antichain{\Bidirected{K_3}}$ as butterfly minors.

\begin{definition}[$C_4$-Cockade]
	A graph $G$ is a \emph{$C_4$-cockade} if it can be constructed by means of $2$-clique sums from a number of disjoint copies of $C_4$.
	A graph is a \emph{partial $C_4$-cockade} if it is a subgraph of a $C_4$-cockade.
\end{definition}

\begin{theorem}\cite{thomassen1986sign}]\label{thm:cockades}
	A symmetric digraph $\Bidirected{G}$ is excludes the members of $\Antichain{\Bidirected{K_3}}$ as butterfly minors if and only if $G$ is a partial $C_4$-cockade.
\end{theorem}

The generalisation we aim for is best formulated in terms of treewidth and we use the dual notion of brambles in our proof.

\begin{definition}[Treewidth]
	Let $G$ be a graph.
	A \emph{tree decomposition} of $G$ is a tuple $\Brace{T,\beta}$, where $T$ is a tree and $\beta\colon\V{T}\rightarrow 2^{\V{G}}$ is a function, called the \emph{bags} of $\Brace{T,\beta}$, such that the following properties hold:
	\begin{enumerate}
		\item $\bigcup_{t\in\V{T}}\Fkt{\beta}{t}=\V{G}$,
		\item for every $e\in\E{G}$ there exists $t_e\in\V{T}$ such that $e\subseteq\Fkt{\beta}{t_e}$, and
		\item for every $v\in\V{G}$, the set $\CondSet{t\in\V{T}}{v\in\Fkt{\beta}{t}}$ induces a subtree of $T$.
	\end{enumerate}
	The \emph{width} of $\Brace{T,\beta}$ is defined as $\Width{\Brace{T,\beta}}\coloneqq \max_{t\in\V{T}}\Abs{\Fkt{\beta}{t}}-1$.
	The \emph{treewidth} of $G$, denoted by $\tw{G}$, is the minimum width over all tree decompositions of $G$.
\end{definition}

	Let $G$ be a graph.
	A \emph{bramble} $\mathcal{B}=\Set{B_1,B_2,\dots,B_{\ell}}$ of $G$ is a family of connected and pairwise touching subgraphs $B_i$ of $G$.
	The \emph{order} of $\mathcal{B}$ is the size of a minimum hitting set for $B$.
	The \emph{bramble number}, denoted by $\BrambleNumber{G}$, is the maximum order over all brambles in $G$.

\begin{theorem}[\cite{seymour1993graph}]\label{thm:brambleduality}
	Let $G$ be a graph, and $k\in\N$ a positive integer.
	Then $G$ contains a bramble of order $k$ if and only if $\tw{G}\geq k-1$.
	Hence $\BrambleNumber{G}-1=\tw{G}$.
\end{theorem}

	There also exists a directed analogue of brambles.
	
	\begin{definition}[Directed Bramble]
		Let $D$ be a digraph.
		A \emph{directed bramble} in $D$ is a family $\mathcal{B}$ of strongly connected subgraphs of $D$ such that for every pair $B_1,B_2\in\mathcal{B}$, either $\V{B_1}\cap\V{B_2}\neq\emptyset$, or there exist edges $\Brace{u_1,v_2},\Brace{u_2,v_1}\in\E{D}$ with $u_i,v_i\in\V{B_i}$ for both $i\in[1,2]$.
		
		A \emph{cover} or \emph{hitting set} for $\mathcal{B}$ is a set $S\subseteq\V{D}$ such that $S\cap\V{D}\neq\emptyset$ for all $B\in\mathcal{B}$.
		
		The \emph{order} of $\mathcal{B}$ is the minimum size of a hitting set for $\mathcal{B}$.
		The \emph{directed bramble number} of a digraph $D$, denoted by $\DirectedBrambleNumber{D}$, is the maximum order of a bramble in $D$.
	\end{definition}
	
	\begin{theorem}[\cite{kreutzer2014width}]\label{thm:directedbrambleduality}
		Let $D$ be a digraph.
		Then $\DirectedBrambleNumber{D}\leq\dtw{D}\leq6\DirectedBrambleNumber{D}+1$.
	\end{theorem}

An immediate and easy to check corollary of \cref{thm:cockades} is the following:

\begin{corollary}\label{cor:twofcockades}
	Let $G$ be a graph.
	If $G$ is a partial $C_4$-cockade, then $\tw{G}\leq 2$.
\end{corollary}

By combining \cref{thm:cockades} with \cref{cor:twofcockades}, \cref{thm:pfaffian}, Lemma 4.7 from \cite{rabinovich2019cyclewidth} and Proposition 3.3 from \cite{giannopoulou2021excluding} one obtains the following interesting relation between the existence of $K_{3,3}$ in $\Split{\Bidirected{G}}$ and $\tw{G}$.

\begin{corollary}\label{cor:K33trees}
	Let $G$ be a graph.
	If $\Split{\Bidirected{G}}$ contains $K_{3,3}$ as a \hyperref[def:matchingminor]{matching minor} then $\tw{G}\geq 2$, and if $\tw{G}\geq 3$, then $\Split{\Bidirected{G}}$ contains $K_{3,3}$ as a matching minor.
\end{corollary}

So the existence of $K_{3,3}$ as a matching minor in $\Split{\Bidirected{G}}$ is closely linked to the treewidth of $G$.
The main result of this subsection is the following generalisation of \cref{cor:K33trees}.

\begin{theorem}\label{thm:Kttduality}
	There exists a function $f\colon\N\rightarrow\N$ such that for every $t\in\N$ and every graph $G$ the following two statements hold:
	\begin{enumerate}
		\item If $\Split{\Bidirected{G}}$ contains $K_{t,t}$ as a \hyperref[def:matchingminor]{matching minor}, then $\tw{G}\geq\frac{1}{2}t-1$, and
		\item if $\tw{G}\geq\Fkt{f}{t}$, then $\Split{\Bidirected{G}}$ contains $K_{t,t}$ as a matching minor.
	\end{enumerate}
\end{theorem}

We start by showing that any large enough wall $W$ contains $K_{t,t}$ as a \hyperref[def:matchingminor]{matching minor} in $\Split{\Bidirected{W}}$.
Since, by \cref{cor:undirwall}, every graph of large enough treewidth contains a large wall, \cref{thm:Kttduality} follows almost immediately.
Similarly, as $W$ is a planar graph, $\Bidirected{W}$ is a planar digraph and with $K_{t,t}$ being a matching minor of $\Split{\Bidirected{K_{t,t}}}$, \cref{lemma:excludingantichains} implies that $\Bidirected{W}$ must contain a member of $\Antichain{\Bidirected{K_t}}$ as a butterfly minor.
Since all butterfly minors of $\Bidirected{W}$ must be planar, we obtain the following as a corollary.

\begin{corollary}\label{thm:allantichainsareplanar}
	Let $D$ be a strongly connected digraph.
	Then $\Antichain{D}$ contains a planar digraph.
\end{corollary}

We think that \cref{thm:allantichainsareplanar} is of particular importance as it underlines the fundamental difference between the genus of a digraph and its \hyperref[def:stronggenus]{strong genus}.
In particular when it comes to the theory of butterfly minors.

\begin{lemma}\label{lemma:wallKtt}
	There exists a function $f_0\colon\N\rightarrow\N$ such that for every $t\in\N$, and every $\Fkt{f_0}{t}$-wall $W$, \hyperref[def:split]{$\Split{\Bidirected{W}}$} contains $K_{t,t}$ as a \hyperref[def:matchingminor]{matching minor}.
\end{lemma}

\begin{proof}
	First observe that every $k$-wall $W_k$ has treewidth at least $k$.
	To see this, we may construct a bramble of order $k+1$ as follows:
	Let $P$ be the right most path from top to bottom that contains two hub-vertices of every row of $W_k$ except for the bottom one, from which it is disjoint.
	Next let $B$ be the bottom row.
	At last, for each row $i$ except the bottom one let $T_i$ be the $i$th row except for the last two vertices together with the $i$th column, from left to right.
	Then each $T_i$ has an edge to $P$ and intersects $B$ and each other $T_j$.
	However, no vertex is contained in more than two of these subgraphs, and $P$ is disjoint from all of them.
	So a minimum hitting set for all of these subgraphs must contain a vertex from every row of $W_k$ and an additional vertex to cover $P$ which makes in total $k+1$ vertices.
	
	Moreover, note that this is also a directed bramble in $\Bidirected{W_K}$, hence $\dtw{\Bidirected{W_k}}\geq k$ by \cref{thm:directedbrambleduality}.
	Consequently, for every $t\in\N$, we are guaranteed by \cref{cor:directedwall} that every $\Fkt{\DirectedGrid}{2t}$-wall $W$ satisfies that $\Bidirected{W}$ contains a cylindrical $t$-wall.
	
	Now let $W$ be a $\Fkt{\DirectedGrid}{2\Fkt{\WallBound}{t,2}}$-wall.
	Then, as discussed above, $\Bidirected{W}$ contains a cylindrical $\Fkt{\WallBound}{t,2}$-wall $U$ as a subgraph.
	
	Suppose $\Split{\Bidirected{W}}$ does not contain $K_{t,t}$ as a matching minor.
	Let $B\coloneqq\Split{\Bidirected{W}}$, and let $M$ be the perfect matching of $B$ such that $\DirM{B}{M}=\Bidirected{W}$.
	Then $\Split{U}$ is an $M$-conformal subgraph of $B$.
	By \cref{thm:matchingflatwall} there must exist an $M$-conformal set $A\subseteq M$ of size at most $\Fkt{\ApexBound}{t}$, such that $B-M$ contains an $M$-conformal \hyperref[def:matchingwall]{matching $2$-wall} $U'\subseteq U$ which is \hyperref[def:matchingflat]{$\Perimeter{U'}$-flat} in $B$ with respect to $A$.
	Let $B''$ be the Pfaffian brace that is the host of $U'$ in a $\Perimeter{U'}$-base of $B-A$.
	
	Note that being a brace means being strongly $2$-connected in the setting of digraphs, which translates to $2$-connectivity in the setting of undirected graphs.
	Moreover, a directed tight cut contraction in a symmetric digraphs is isomorphic to deleting all non-separator vertices of one of the shores.
	Hence $\Undirected{\DirM{B''}{M}}$ is an actual subgraph of $W$.
	Let us call this subgraph $H$.
	For $H$ we have the following informations:
	\begin{enumerate}
		\item $H$ contains $\Undirected{\DirM{U'}{M}}$ as a subgraph, and
		\item $\Bidirected{H}$ excludes all members of $\Antichain{\Bidirected{K_3}}$ as butterfly minors by \cref{thm:pfaffian}.
	\end{enumerate}
	\begin{figure}[!ht]
		\centering
		\begin{tikzpicture}[scale=0.5]
			\pgfdeclarelayer{background}
			\pgfdeclarelayer{foreground}
			\pgfsetlayers{background,main,foreground}
			
			\node (C1) [v:ghost] {};
			\node (C2) [v:ghost,position=0:130mm from C1] {};
			
			\node (D1) [v:ghost,position=0:0mm from C1] {
				
				\begin{tikzpicture}[scale=1]
					\pgfdeclarelayer{background}
					\pgfdeclarelayer{foreground}
					\pgfsetlayers{background,main,foreground}
					
					\node (C) [v:ghost] {};
					
					\node (v1) [v:main,position=22.5:18mm from C] {};
					\node (v2) [v:main,position=45:18mm from C] {};
					\node (v3) [v:main,position=67.5:18mm from C] {};
					\node (v4) [v:main,position=90:18mm from C] {};
					\node (v5) [v:main,position=112.5:18mm from C] {};
					\node (v6) [v:main,position=135:18mm from C] {};
					\node (v7) [v:main,position=157.5:18mm from C] {};
					\node (v8) [v:main,position=180:18mm from C] {};
					\node (v9) [v:main,position=202.5:18mm from C] {};
					\node (v10) [v:main,position=225:18mm from C] {};
					\node (v11) [v:main,position=247.5:18mm from C] {};
					\node (v12) [v:main,position=270:18mm from C] {};
					\node (v13) [v:main,position=292.5:18mm from C] {};
					\node (v14) [v:main,position=315:18mm from C] {};
					\node (v15) [v:main,position=337.5:18mm from C] {};
					\node (v16) [v:main,position=0:18mm from C] {};
					
					\node (u1) [v:main,position=22.5:26mm from C] {};
					\node (u2) [v:main,position=45:26mm from C] {};
					\node (u3) [v:main,position=67.5:26mm from C] {};
					\node (u4) [v:main,position=90:26mm from C] {};
					\node (u5) [v:main,position=112.5:26mm from C] {};
					\node (u6) [v:main,position=135:26mm from C] {};
					\node (u7) [v:main,position=157.5:26mm from C] {};
					\node (u8) [v:main,position=180:26mm from C] {};
					\node (u9) [v:main,position=202.5:26mm from C] {};
					\node (u10) [v:main,position=225:26mm from C] {};
					\node (u11) [v:main,position=247.5:26mm from C] {};
					\node (u12) [v:main,position=270:26mm from C] {};
					\node (u13) [v:main,position=292.5:26mm from C] {};
					\node (u14) [v:main,position=315:26mm from C] {};
					\node (u15) [v:main,position=337.5:26mm from C] {};
					\node (u16) [v:main,position=0:26mm from C] {};
					
					\begin{pgfonlayer}{background}
						
						\draw [e:mainplus,->,bend right=12] (v1) to (v2);
						\draw [e:mainplus,->,bend right=12] (v2) to (v3);
						\draw [e:mainplus,->,bend right=12] (v3) to (v4);
						\draw [e:mainplus,->,bend right=12] (v4) to (v5);
						\draw [e:mainplus,->,bend right=12] (v5) to (v6);
						\draw [e:mainplus,->,bend right=12] (v6) to (v7);
						\draw [e:mainplus,->,bend right=12] (v7) to (v8);
						\draw [e:mainplus,->,bend right=12] (v8) to (v9);
						\draw [e:mainplus,->,bend right=12] (v9) to (v10);
						\draw [e:mainplus,->,bend right=12] (v10) to (v11);
						\draw [e:mainplus,->,bend right=12] (v11) to (v12);
						\draw [e:mainplus,->,bend right=12] (v12) to (v13);
						\draw [e:mainplus,->,bend right=12] (v13) to (v14);
						\draw [e:mainplus,->,bend right=12] (v14) to (v15);
						\draw [e:mainplus,->,bend right=12] (v15) to (v16);
						\draw [e:mainplus,->,bend right=12] (v16) to (v1);
						
						\draw [e:main,<-,color=Gray,bend left=12] (v1) to (v2);
						\draw [e:main,<-,color=Gray,bend left=12] (v2) to (v3);
						\draw [e:main,<-,color=Gray,bend left=12] (v3) to (v4);
						\draw [e:main,<-,color=Gray,bend left=12] (v4) to (v5);
						\draw [e:main,<-,color=Gray,bend left=12] (v5) to (v6);
						\draw [e:main,<-,color=Gray,bend left=12] (v6) to (v7);
						\draw [e:main,<-,color=Gray,bend left=12] (v7) to (v8);
						\draw [e:main,<-,color=Gray,bend left=12] (v8) to (v9);
						\draw [e:main,<-,color=Gray,bend left=12] (v9) to (v10);
						\draw [e:main,<-,color=Gray,bend left=12] (v10) to (v11);
						\draw [e:main,<-,color=Gray,bend left=12] (v11) to (v12);
						\draw [e:main,<-,color=Gray,bend left=12] (v12) to (v13);
						\draw [e:main,<-,color=Gray,bend left=12] (v13) to (v14);
						\draw [e:main,<-,color=Gray,bend left=12] (v14) to (v15);
						\draw [e:main,<-,color=Gray,bend left=12] (v15) to (v16);
						\draw [e:main,<-,color=Gray,bend left=12] (v16) to (v1);
						
						\draw [e:mainplus,->,bend right=12] (u1) to (u2);
						\draw [e:mainplus,->,bend right=12] (u2) to (u3);
						\draw [e:mainplus,->,bend right=12] (u3) to (u4);
						\draw [e:mainplus,->,bend right=12] (u4) to (u5);
						\draw [e:mainplus,->,bend right=12] (u5) to (u6);
						\draw [e:mainplus,->,bend right=12] (u6) to (u7);
						\draw [e:mainplus,->,bend right=12] (u7) to (u8);
						\draw [e:mainplus,->,bend right=12] (u8) to (u9);
						\draw [e:mainplus,->,bend right=12] (u9) to (u10);
						\draw [e:mainplus,->,bend right=12] (u10) to (u11);
						\draw [e:mainplus,->,bend right=12] (u11) to (u12);
						\draw [e:mainplus,->,bend right=12] (u12) to (u13);
						\draw [e:mainplus,->,bend right=12] (u13) to (u14);
						\draw [e:mainplus,->,bend right=12] (u14) to (u15);
						\draw [e:mainplus,->,bend right=12] (u15) to (u16);
						\draw [e:mainplus,->,bend right=12] (u16) to (u1);
						
						\draw [e:main,<-,color=Gray,bend left=12] (u1) to (u2);
						\draw [e:main,<-,color=Gray,bend left=12] (u2) to (u3);
						\draw [e:main,<-,color=Gray,bend left=12] (u3) to (u4);
						\draw [e:main,<-,color=Gray,bend left=12] (u4) to (u5);
						\draw [e:main,<-,color=Gray,bend left=12] (u5) to (u6);
						\draw [e:main,<-,color=Gray,bend left=12] (u6) to (u7);
						\draw [e:main,<-,color=Gray,bend left=12] (u7) to (u8);
						\draw [e:main,<-,color=Gray,bend left=12] (u8) to (u9);
						\draw [e:main,<-,color=Gray,bend left=12] (u9) to (u10);
						\draw [e:main,<-,color=Gray,bend left=12] (u10) to (u11);
						\draw [e:main,<-,color=Gray,bend left=12] (u11) to (u12);
						\draw [e:main,<-,color=Gray,bend left=12] (u12) to (u13);
						\draw [e:main,<-,color=Gray,bend left=12] (u13) to (u14);
						\draw [e:main,<-,color=Gray,bend left=12] (u14) to (u15);
						\draw [e:main,<-,color=Gray,bend left=12] (u15) to (u16);
						\draw [e:main,<-,color=Gray,bend left=12] (u16) to (u1);
						
						\draw [e:mainplus,->,bend right=12] (v1) to (u16);
						\draw [e:mainplus,->,bend right=12] (u3) to (v2);
						
						\draw [e:mainplus,->,bend right=12] (v5) to (u4);
						\draw [e:mainplus,->,bend right=12] (u7) to (v6);
						
						\draw [e:mainplus,->,bend right=12] (v9) to (u8);
						\draw [e:mainplus,->,bend right=12] (u11) to (v10);
						
						\draw [e:mainplus,->,bend right=12] (v13) to (u12);
						\draw [e:mainplus,->,bend right=12] (u15) to (v14);
						
						\draw [e:main,<-,color=Gray,bend left=12] (v1) to (u16);
						\draw [e:main,<-,color=Gray,bend left=12] (u3) to (v2);
						
						\draw [e:main,<-,color=Gray,bend left=12] (v5) to (u4);
						\draw [e:main,<-,color=Gray,bend left=12] (u7) to (v6);
						
						\draw [e:main,<-,color=Gray,bend left=12] (v9) to (u8);
						\draw [e:main,<-,color=Gray,bend left=12] (u11) to (v10);
						
						\draw [e:main,<-,color=Gray,bend left=12] (v13) to (u12);
						\draw [e:main,<-,color=Gray,bend left=12] (u15) to (v14);
						
					\end{pgfonlayer}
				\end{tikzpicture}
				
			};
			\node (D2) [v:ghost,position=0:0mm from C2] {
				
				\begin{tikzpicture}[scale=1]
					\pgfdeclarelayer{background}
					\pgfdeclarelayer{foreground}
					\pgfsetlayers{background,main,foreground}
					
					\node (C) [v:ghost] {};
					
					\node (v1) [v:main,position=22.5:18mm from C] {};
					\node (v2) [v:main,position=45:18mm from C] {};
					\node (v3) [v:main,position=67.5:18mm from C] {};
					\node (v4) [v:main,position=90:18mm from C] {};
					\node (v5) [v:main,position=112.5:18mm from C] {};
					\node (v6) [v:main,position=135:18mm from C] {};
					\node (v7) [v:main,position=157.5:18mm from C] {};
					\node (v8) [v:main,position=180:18mm from C] {};
					\node (v9) [v:main,position=202.5:18mm from C] {};
					\node (v10) [v:main,position=225:18mm from C] {};
					\node (v11) [v:main,position=247.5:18mm from C] {};
					\node (v12) [v:main,position=270:18mm from C] {};
					\node (v13) [v:main,position=292.5:18mm from C] {};
					\node (v14) [v:main,position=315:18mm from C] {};
					\node (v15) [v:main,position=337.5:18mm from C] {};
					\node (v16) [v:main,position=0:18mm from C] {};
					
					\node (v111) [v:ghost,position=22.5:0.25mm from v1] {};
					\node (v22) [v:ghost,position=45:0.25mm from v2] {};
					\node (v33) [v:ghost,position=67.5:0.25mm from v3] {};
					\node (v44) [v:ghost,position=90:0.25mm from v4] {};
					\node (v55) [v:ghost,position=112.5:0.25mm from v5] {};
					\node (v66) [v:ghost,position=135:0.25mm from v6] {};
					\node (v77) [v:ghost,position=157.5:0.25mm from v7] {};
					\node (v88) [v:ghost,position=180:0.25mm from v8] {};
					\node (v99) [v:ghost,position=202.5:0.25mm from v9] {};
					\node (v1010) [v:ghost,position=225:0.25mm from v10] {};
					\node (v1111) [v:ghost,position=247.5:0.25mm from v11] {};
					\node (v1212) [v:ghost,position=270:0.25mm from v12] {};
					\node (v1313) [v:ghost,position=292.5:0.25mm from v13] {};
					\node (v1414) [v:ghost,position=315:0.25mm from v14] {};
					\node (v1515) [v:ghost,position=337.5:0.25mm from v15] {};
					\node (v1616) [v:ghost,position=0:0.25mm from v16] {};
					
					\node (u1) [v:main,position=22.5:26mm from C] {};
					\node (u2) [v:main,position=45:26mm from C] {};
					\node (u3) [v:main,position=67.5:26mm from C] {};
					\node (u4) [v:main,position=90:26mm from C] {};
					\node (u5) [v:main,position=112.5:26mm from C] {};
					\node (u6) [v:main,position=135:26mm from C] {};
					\node (u7) [v:main,position=157.5:26mm from C] {};
					\node (u8) [v:main,position=180:26mm from C] {};
					\node (u9) [v:main,position=202.5:26mm from C] {};
					\node (u10) [v:main,position=225:26mm from C] {};
					\node (u11) [v:main,position=247.5:26mm from C] {};
					\node (u12) [v:main,position=270:26mm from C] {};
					\node (u13) [v:main,position=292.5:26mm from C] {};
					\node (u14) [v:main,position=315:26mm from C] {};
					\node (u15) [v:main,position=337.5:26mm from C] {};
					\node (u16) [v:main,position=0:26mm from C] {};
					
					\node (u44) [v:ghost,position=0:0.5mm from u4] {};
					\node (u1515) [v:ghost,position=70:0.5mm from u15] {};
					\node (u33) [v:ghost,position=160:0.5mm from u3] {};
					\node (u111) [v:ghost,position=285:0.5mm from u1] {};
					
					\begin{pgfonlayer}{background}
						
						\draw[e:marker,color=CornflowerBlue] (v111) to (v22) to (v33) to (v44) to (v55) to (v66) to (v77) to (v88) to (v99) to (v1010) to (v1111) to (v1212) to (v1313) to (v1414) to (v1515) to (v1616) to (v111);
						
						\node (uu16) [v:marked,scale=9,color=CornflowerBlue,position=0:0mm from u16] {};
						
						\draw[e:marker,color=CornflowerBlue] (u44) to (u5) to (u6) to (u7) to (u8) to (u9) to (u10) to (u11) to (u12) to (u13) to (u14) to (u1515);
						
						\draw[e:marker,color=CornflowerBlue] (u33) to (u2) to (u111);
						
						\draw [e:main] (v1) to (v2);
						\draw [e:main] (v2) to (v3);
						\draw [e:main] (v3) to (v4);
						\draw [e:main] (v4) to (v5);
						\draw [e:main] (v5) to (v6);
						\draw [e:main] (v6) to (v7);
						\draw [e:main] (v7) to (v8);
						\draw [e:main] (v8) to (v9);
						\draw [e:main] (v9) to (v10);
						\draw [e:main] (v10) to (v11);
						\draw [e:main] (v11) to (v12);
						\draw [e:main] (v12) to (v13);
						\draw [e:main] (v13) to (v14);
						\draw [e:main] (v14) to (v15);
						\draw [e:main] (v15) to (v16);
						\draw [e:main] (v16) to (v1);

						\draw [e:main] (u1) to (u2);
						\draw [e:main] (u2) to (u3);
						\draw [e:main] (u3) to (u4);
						\draw [e:main] (u4) to (u5);
						\draw [e:main] (u5) to (u6);
						\draw [e:main] (u6) to (u7);
						\draw [e:main] (u7) to (u8);
						\draw [e:main] (u8) to (u9);
						\draw [e:main] (u9) to (u10);
						\draw [e:main] (u10) to (u11);
						\draw [e:main] (u11) to (u12);
						\draw [e:main] (u12) to (u13);
						\draw [e:main] (u13) to (u14);
						\draw [e:main] (u14) to (u15);
						\draw [e:main] (u15) to (u16);
						\draw [e:main] (u16) to (u1);
						
						\draw [e:main] (v1) to (u16);
						\draw [e:main] (u3) to (v2);
						
						\draw [e:main] (v5) to (u4);
						\draw [e:main] (u7) to (v6);
						
						\draw [e:main] (v9) to (u8);
						\draw [e:main] (u11) to (v10);
						
						\draw [e:main] (v13) to (u12);
						\draw [e:main] (u15) to (v14);

					\end{pgfonlayer}
				\end{tikzpicture}
				
			};
			
			\begin{pgfonlayer}{background}
				
			\end{pgfonlayer}

		\end{tikzpicture}
		\caption{An elementary cylindrical $2$-wall as a subgraph of a symmetric digraph $\Bidirected{G}$ (left) and a bramble of order $4$ in $G$ (right).}
		\label{fig:undirecktedcylindricalwall}
	\end{figure}
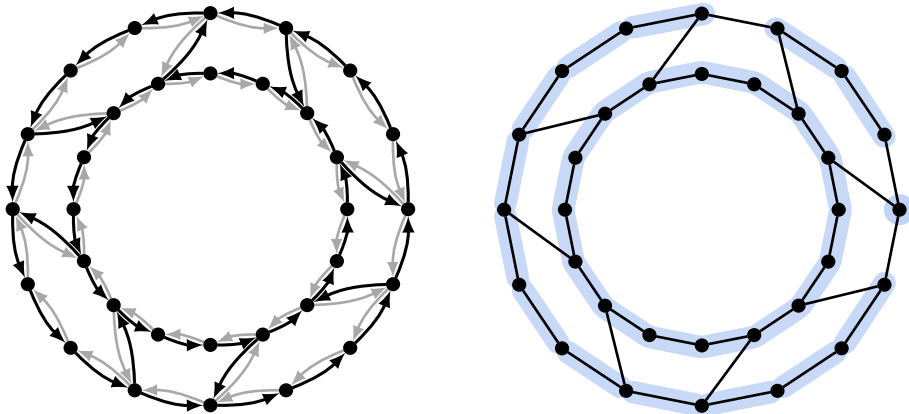
	Consider \cref{fig:undirecktedcylindricalwall}.
	Here we depict a cylindrical $2$-wall $Q$ as a subgraph of a symmetric digraph (on the left), and its underlying undirected graph (on the right).
	We consider the following family of subgraphs of $Q$: Let $v_1$ be a hub-vertex, i.\@e.\@ a vertex of degree $3$ in $Q$, on the outer cycle $C'$ of $Q$, then let $P$ be a shortest subpath of the outer cycle of $Q$ such that one endpoint of $P$ is a neighbour of $v_1$ and the other endpoint is a hub-vertex of $Q$.
	Next let $L\coloneqq C'-v_1-P$, and at last let $C\coloneqq Q-v_1-P-L$.
	The result are four pairwise disjoint connected subgraphs, one of them only consisting of $v_1$, of $Q$ that are pairwise joined by an edge.
	Hence the collection of these four subgraphs forms a bramble of order four, implying $\tw{H}\geq 3$.
	By \cref{cor:K33trees} this is a contradiction to $\Bidirected{H}$ excluding $\Antichain{\Bidirected{K_3}}$ as butterfly minors and therefore this contradicts \cref{thm:matchingflatwall} or $\Split{\Bidirected{W}}$ contains $K_{t,t}$ as a matching minor.
	Thus we may assume the latter and the proof is complete.
\end{proof}

As discussed above, \cref{thm:allantichainsareplanar} follows immediately.
Hence all that is left to do is to prove \cref{thm:Kttduality}.
This concludes the chapter.

\begin{proof}[Proof of \Cref{thm:Kttduality}]
	First suppose $\Split{\Bidirected{G}}$ contains $K_{t,t}$ as a \hyperref[def:matchingminor]{matching minor}.
	Then, by applying Lemma 4.7 from \cite{rabinovich2019cyclewidth} and Proposition 3.3 from \cite{giannopoulou2021excluding} we obtain that $\dtw{\Bidirected{G}}\geq\frac{1}{2}t-1$.
	Since $\dtw{\Bidirected{G}}=\tw{G}$ the first part of our claim follows.
	
	For the second part let $f_0$ be the function from \cref{lemma:wallKtt}, and let $\Fkt{f}{t}\coloneqq \Fkt{\UndirectedGrid}{2\Fkt{f_0}{t}}$.
	Then, by \cref{cor:undirwall}, $G$ contains an $\Fkt{f_0}{t}$-wall as a subgraph, and by \cref{lemma:wallKtt} this means that $\Split{\Bidirected{G}}$ must contain $K_{t,t}$ as a matching minor.
\end{proof}

\section{A Brief Summary of the Proof}

The next three sections are dedicated to the proof of \cref{thm:matchingflatwall}.
Here we follow the known proofs of the (Directed) Flat Wall Theorem with some slight alterations.
Specifically, given a bipartite graph $B$ with a perfect matching $M$ and a large $M$-conformal \hyperref[def:matchingwall]{matching wall} $W$ we show that
\begin{enumerate}
	\item If there are many pairwise disjoint \hyperref[def:alternatingpath]{internally $M$-conformal paths} that are internally disjoint from $W$ and have both of their endpoints in $W$ but they are far apart from each other in $W$, we find $K_{t,t}$ as an $M$-minor grasped by $W$.
	This is done in \cref{sec:longjumps} and we essentially adapt the tools introduced in \cite{giannopoulou2020directed} to achieve this.
	In particular, this part of the proof is done mostly in the setting of digraphs.
	
	\item The second part consists of two steps at once, both of which can be solved by the same technique, but since they are slightly different, we explain both:
	\begin{enumerate}
		\item In case there are many pairwise disjoint \hyperref[def:alternatingpath]{internally $M$-conformal paths} that are internally disjoint from $W$ whose endpoints both lie in $W$ and are close together but not in the same cell, we can find many pairwise disjoint matching minor models of $K_{3,3}$.
		These $K_{3,3}$-matching models yield many local crosses which can be used to construct a \hyperref[def:matchingminor]{matching minor} model of $K_{t,t}$ grasped by $W$.
		
		\item Finally, we know that every internally $M$-conformal path that is internally disjoint from $W$ must have both endpoints on the same cell of $W$.
		Hence we may associate with every cell of $W$ a bipartite matching covered graph that is otherwise disconnected from $W$.
		If many of these cells are, essentially non-Pfaffian, we can again find many pairwise disjoint models of $K_{3,3}$ which then can be used to construct a matching minor model of $K_{t,t}$ grasped by $W$.
	\end{enumerate}
\end{enumerate}
While the second and third steps are relatively similar to the proof of the undirected Flat Wall Theorem, they differ vastly from their directed analogues.
Both steps, (ii) and (iii), are discussed in \cref{sec:crosses}.
The actual proof of \cref{thm:matchingflatwall} and thus the combination of all three steps is then done in \cref{sec:proof}.

\section{Step 1: Remove Long Jumps}\label{sec:longjumps}

In general the main lemmas of this section will take as input a wall and a tiling, then they create an auxiliary graph and ask for a large amount of pairwise disjoint paths from one vertex set to another in the auxiliary graph.
These paths from the auxiliary graph will then correspond to `long jumps' over the wall which will allow for us to create a $\Bidirected{K_t}$-butterfly minor in the $M$-direction of our graph.
By \cref{lemma:mcguigmatminors} this means that we have found a $K_{t,t}$-\hyperref[def:matchingminor]{matching minor}.

We start, however, by extracting a utility lemma from \cite{giannopoulou2020directed} that helps us to create the $\Bidirected{K_t}$-butterfly minor.
For this we need further notation.

Recall the definition of a \hyperref[def:tile]{tile}.
The notion of \hyperref[def:tile]{tiles} allows us to add another layer of parametrisation on top of a cylindrical wall.
In some sense a tile can be seen as a generalisation of a cell that also contains a small acyclic wall inside to allow for additional routing.
The next few definitions are used to add further details to how our walls are divided into different regions and how tiles are used to achieve this additional layer of parametrisation.

\begin{definition}[Triadic Partitions]\label{def:triadicpartition}
	Let $k\in\N$ be a positive integer and $W=\Brace{Q_1,\dots,Q_{3k},\hat{P}_1,\dots,\hat{P}_{3k}}$ be a cylindrical $3k$-wall.
	The \emph{triadic partition} of $W$ is the tuple
	\begin{align*}
		\Triade=\Brace{W,k,W_1,W_2,W_3,W^1,W^2,W^3}
	\end{align*}
	such that for each $i\in[1,3]$, $W_i$ denotes the \hyperref[def:slice]{slice} of $W$ between $Q_{k\Brace{i-1}+1}$ and $Q_{ik}$, and $W^i$ denotes the strip of $W$ between the rows $k\Brace{i-1}+1$ and $ik$.
\end{definition}

\begin{definition}[Tiling]
	Let $k\in\N$ be a positive integer and $W=\Brace{Q_1,\dots,Q_{3k},\hat{P}_1,\dots,\hat{P}_{3k}}$ be a cylindrical $3k$-wall with its \hyperref[def:triadicpartition]{triadic partition} $\Triade=\Brace{W,k,W_1,W_2,W_3,W^1,W^2,W^3}$.
	
	A \emph{tiling} is a family of pairwise disjoint \hyperref[def:tile]{tiles} $\Tiling$.
	Let $W'$ be a \hyperref[def:slice]{slice} of $W$, we say that $\Tiling$ \emph{covers} $W'$ if every vertex $v\in\V{W'}$ with $\DegG{W}{v}=3$ is contained in some tile of $\Tiling$.
	The tiling $\Tiling$ is said to \emph{cover} $\Triade$ if it covers $W_2$.
	
	In most places we will use the following family of tilings:
	Let $f_w\colon\N\rightarrow\N$ be some function, $t\in\N$ a positive integer, and $\xi,\xi'\in[1,\Fkt{f_w}{t}+1]$.
	We define the \emph{column function} $\ColumnFunction\colon\N\rightarrow\N$ and the \emph{row function} $\RowFunction\colon\N\rightarrow\N$ as follows\footnote{Please note that $\ColumnFunction$ and $\RowFunction$ do also depend on $f_w$, $\xi$, $\xi'$, and $k$. However, it is more convenient to make these dependencies implicit.}:
	\begin{align*}
		\Fkt{\ColumnFunction}{p}&\coloneqq \Brace{k+2-\xi}+\Brace{p-1}\Brace{2\Fkt{f_w}{t}+1}\text{, and}\\
		\Fkt{\RowFunction}{q}&\coloneqq \xi'+\Brace{q-1}\Brace{2\Fkt{f_w}{t}+1}.
	\end{align*}
	We can now define our standard tiling for fixed $f_w$, $\xi$, and $\xi'$.
	\begin{align*}
		\Tiling_{W,k,\Fkt{f_w}{t},\xi,\xi'}\coloneqq\left\{T_{\Fkt{\ColumnFunction}{p},\Fkt{\RowFunction}{q},\Fkt{f_w}{t}} \mid p\in\left[1,\Ceil{\frac{k+\xi-1}{2\Fkt{f_w}{t}+1}}+1\right]\text{, }q\in\left[1,\Ceil{\frac{3k-\xi'-1}{2\Fkt{f_w}{t}+1}}+1\right]\right\}
	\end{align*}
\end{definition}

Note that every tiling $\Tiling_{W,k,\Fkt{f_w}{t}.\xi,\xi'}$ covers $\Triade$.
Moreover, every cell of $W_2$ that lies between the two paths of $\hat{P}_i$ for some $i\in[1,3k]$ is the centre of some tile $T'$ of some tiling $\Tiling'\in\Tiling_{W,k,\Fkt{f_w}{t}.\xi,\xi'}$.
Hence if we perform the mirror-image operation as described after the definition of tiles, we are able to find in total $2\Brace{\Fkt{f_w}{t}+1}^2$ many tilings that cover $W_2$, such that every cell of $W_2$ is the centre of some tile in one of these tilings.

We will use tilings in several different ways, and sometimes it is necessary to `zoom out' of our current wall, i.\@e.\@ to forget about some of the horizontal paths and vertical cycles in order to obtain a more streamlined version of our wall.

\begin{definition}[Walls from a Tiling]\label{def:tiling}
	Let $k,d\in\N$ be positive integers, $W=\Brace{Q_1,\dots,Q_{3k},\hat{P}_1,\dots,\hat{P}_{3k}}$ be a cylindrical $3k$-wall with its triadic partition $\Triade=\Brace{W,k,W_1,W_2,W_3,W^1,W^2,W^3}$, and $\Tiling$ be a tiling of width $d$ that covers $W_2$.
	Moreover, let $\widetilde{W}$ be the some slice of $W_2$ and let $I_Q$ be the largest set of integers such that $Q_i$ contains vertices of a tile from $\Tiling$ which intersects $\widetilde{W}$ for every $i\in I_Q$.
	Let $\Fkt{\widetilde{W}}{\Tiling}$ be the union of the cycles $Q_i$, $i\in I_Q$, and the paths $\InducedSubgraph{P_j^i}{\CondSet{Q_h}{h\in I_Q}}$ for every $\Brace{j,i}\in [1,3k]\times[1,2]$.
	We call $\Fkt{\widetilde{W}}{\Tiling}$ the \emph{extension of $\widetilde{W}$ that covers $\Tiling$}.
	
	Now, let $\Set{\Class_1,\dots,\Class_4}$ be a four colouring of $\Tiling$ and $i\in[1,4]$ be a fixed colour.
	Then let $J_Q\subseteq[1,3k]$ be the largest set of integers such that for all $j\in J_Q$ the vertical cycle $Q_j$ of $\Fkt{\widetilde{W}}{\Tiling}$ does not contain a vertex of some tile from $\Class_i$.
	Similarly, let $J_P\subseteq[1,3k]$ be the largest set of integers such that for every $j\in J_P$, none of the two paths from $\hat{P}_j$ contains a vertex of a tile from $\Class_i$.
	
	By $\InducedSubgraph{\widetilde{W}}{\Tiling,i}$ we denote the subgraph of $W$ induced by the union of the cycles $Q_i$, $i\in J_Q$, and the paths $\InducedSubgraph{P_j^i}{\CondSet{Q_h}{h\in J_Q}}$ for every $\Brace{j,i}\in J_P$.
	We say that $\InducedSubgraph{\widetilde{W}}{\Tiling,i}$ is the $i$th $\Tiling$-slice of $\widetilde{W}$.
\end{definition}

Note that in $\InducedSubgraph{\widetilde{W}}{\Tiling,i}$, we essentially cut out the tiles of $\Class_i$.
This operation gives us a slice $W'$ of some cylindrical wall for which the perimeter of every tile in $\Class_i$ has become the perimeter of some cell.
Next, we are going to find a tiling of $W'$ such that every tile of $\Class_i$ that belongs to $\widetilde{W}$ is captured by the centre of some tile in the new tiling.

\begin{definition}[Tier II Tiling]\label{def:tilingII}
	Let $t,k,k'\in\N$ be positive integers and $f\colon\N\rightarrow\N$ be some function where $k\geq k'$.
	Let $W$ be a cylindrical $3k$-wall with its \hyperref[def:triadicpartition]{triadic partition} $\Triade=\Brace{W,k,W_1,W_2,W_3,W^1,W^2,W^3}$, and \hyperref[def:tiling]{$\Tiling=\Tiling_{W,k,f,\xi,\xi'}$} for some $\xi,\xi'\in[1,\Fkt{f}{t}+1]$, as well as $\Set{\Class_1,\dots,\Class_4}$ be a four colouring of $\Tiling$ and $i\in[1,4]$ be a fixed colour.
	Moreover, let $\widetilde{W}$ be a \hyperref[def:slice]{slice} of $W_2$ of width $k'$ such that no \hyperref[def:tile]{tile} of $\Class_i$ contains a vertex of $\Perimeter{\widetilde{W}}$ and $\widetilde{\Tiling}$ be the collection of all tiles from $\Tiling$ that contain a vertex of $\widetilde{W}$.
	
	The \emph{tier II tiling for $\widetilde{W}$ and $i$ obtained from $\Tiling$} is defined as the unique tiling $\TierIITiling{\Tiling,i,f}{\widetilde{W}}$ of width $f$ of $\InducedSubgraph{\widetilde{W}}{\Tiling,i}$ such that every $T\in\Class_i\cap\widetilde{\Tiling}$ is in the interior of the centre of some tile in $\TierIITiling{\Tiling,i,f}{\widetilde{W}}$.
\end{definition}

Since every tile in $\Tiling$ consists of $2\Fkt{f}{t}+2$ path pairs, $\TierIITiling{\Tiling,i,f}{\widetilde{W}}$ is well defined and does in fact cover all of $\InducedSubgraph{\widetilde{W}}{\Tiling,i,f}$.

At last, we need to introduce the notion of a wall grasping a butterfly minor.
Let $D$ and $H$ be digraphs, and let $W$ be a cylindrical wall in $D$.
Let $B=\Split{D}$, and $M$ be the perfect matching of $B$ such that $\DirM{B}{M}=D$.
We say that $W$ \emph{grasps} an $H$-butterfly minor of $D$ if $\Split{W}$ $M$-grasps a $\Split{H}$-\hyperref[def:matchingminor]{matching minor} of $B$.

While \cref{thm:directedflatwall} could supply us with an intermediate wall which we could then refine further, we aim for a more self-contained proof wherever possible and feasible.
To accomplish this goal, we will use a single lemma from the original proof of \cref{thm:directedflatwall}, namely \cref{lemma:directedlongjumps}, together with a result on paths leaving and re-entering a fixed set of vertices.
Indeed, it is necessary to further refine \cref{lemma:directedlongjumps} for it to fit into the framework of our proof.
We start by introducing the preliminary results.

Let $D$ be a digraph and $X\subseteq\V{D}$.
A \emph{directed $X$-path} is a directed path $P$ of length at least one that has both endpoints in $X$ but is otherwise disjoint from $X$.

\begin{theorem}[\cite{giannopoulou2020directed}]\label{thm:xpaths}
	Let $D$ be a digraph and $X\subseteq\V{D}$.
	For all positive $k\in\N$, there are $k$ pairwise vertex disjoint directed $X$-paths in $D$, or there exists a set $S\subseteq\V{D}$ of size at most $2k$ such that every directed $X$-path in $D$ contains a vertex of $S$.
	
	Furthermore, there is a polynomial time algorithm which, given a digraph $D$ and a set $X\subseteq\V{D}$ as input, outputs $k$ pairwise disjoint directed $X$-paths, or a set $S\subseteq\V{D}$ of size at most $2k$ as above.
\end{theorem}

Let $k,w\in\N$ be positive integers, $W$ be a cylindrical $k$-wall and $W'$ be a slice of $W$.
A directed $\V{W'}$-path $P$ is called a \emph{jump over $W'$} if $\E{P}\cap \E{W'}=\emptyset$.
We say that a directed $\V{W'}$-path $P$ is a \emph{$w$-long jump over $W'$} if for all $\xi,\xi'\in[1,w+1]$ the endpoints of $P$ belong to distinct tiles $T_1$ and $T_2$ of the tiling $\Tiling_{W,k,w,\xi,\xi'}$.

\paragraph{A Lemma from \cite{giannopoulou2020directed}.}

The following is a combination of Lemmas 4.3 to 4.8 from \cite{giannopoulou2020directed} and a proof can be found in the proof of Lemma 4.9 in \cite{giannopoulou2020directed}.
The only difference of Lemma 4.9 from \cite{giannopoulou2020directed} and the statement below is, that we extract the last subcase of Case 1 in its proof as a potential outcome.

\begin{lemma}[\cite{giannopoulou2020directed}]\label{lemma:directedlongjumps}
	There exist functions $f_w\colon\N\rightarrow\N$, $f_P\colon\N\rightarrow\N$, and $f_W\colon\N\rightarrow\N$ such that for every $t\in\N$ the following holds:
	Let 
	\begin{itemize}
		\item $D$ be a digraph,
		\item $W$ be a cylindrical $3k$-wall with $k\geq\Fkt{f_W}{t}$ in $D$,
		\item $\Triade=\Brace{W,k,W_1,W_2,W_3,W^1,W^2,W^3}$ be the \hyperref[def:triadicpartition]{triadic partition} of $W$, and
		\item \hyperref[def:tiling]{$\Tiling=\Tiling_{W,k,\Fkt{f_w}{t}.\xi,\xi'}$} for some $\xi,\xi'\in[1,\Fkt{f_w}{t}+1]$.
	\end{itemize}
	If there exists a subfamily $\Tiling'$ of $\Tiling$ and a family $\Jumps$ of pairwise disjoint directed paths in $D$ with the following properties:
	\begin{enumerate}
		\item Every member of $\Jumps$ is internally disjoint from $W$ but has both endpoints on $W$,
		\item $\Abs{\Tiling'}=\Abs{\Jumps}=\Fkt{f_P}{t}$,
		\item for every $T_{\Fkt{\ColumnFunction}{p},\Fkt{\RowFunction}{q},\Fkt{f_w}{t}}\neq T_{\Fkt{\ColumnFunction}{p'},\Fkt{\RowFunction}{q'},\Fkt{f_w}{t}}\in\Tiling'$ we have $\max\Set{\Abs{p-p'},\Abs{q-q'}}\geq 2$,
		\item there exists a bijection $\Start\colon\Tiling'\rightarrow\Jumps$ ($\End\colon\Tiling'\rightarrow\Jumps$) such that the starting point (endpoint) of the path $\Fkt{\Start}{T}$ ($\Fkt{\End}{T}$) belongs to the \hyperref[def:tile]{centre} of $T$,
		\item $\V{\Fkt{\Start}{T}}\cap\V{\Tiling'}$ ($\V{\Fkt{\End}{T}}\cap\V{\Tiling'}$) contains exactly the endpoint of $\Fkt{\Start}{T}$ ($\Fkt{\End}{T}$) where $\V{\Tiling'}=\bigcup_{T'\in\Tiling'}\V{T'}$, and finally
		\item the endpoints (starting points) of the paths in $\Jumps$ are of mutual \hyperref[def:Wdistance]{$W$-distance} at least $4$. 
	\end{enumerate}
	Then one of the following is true.
	\begin{enumerate}
		\item[a)] $D$ has a $\Bidirected{K_t}$-butterfly minor grasped by $W$,
		\item[b)] there exists a family of tiles $\Tiling''\subseteq\Tiling'$ all contained in a single \hyperref[def:strip]{strip}  $S\subseteq W$ of height equal to the height of the tiles in $\Tiling$ such that
		\begin{itemize}
			\item we can number $\Tiling''=\Set{T_1,\dots,T_h}$ such that for each $i\in[1,h]$ $S-T_i$ has one component containing exactly the tiles $T_1,\dots,T_{i-1}$,
			\item $\Abs{\Tiling''}\geq \Fkt{f_P}{t}^{\frac{1}{4}}$,
			\item for every $i\in[1,h-1]$ the tiles $T_i$ and $T_{i+1}$ are separated in $W$ by a \hyperref[def:slice]{slice} of width equal to the width of the tiles in $\Tiling$, and
			\item there is a family $\Jumps'\subseteq\Jumps$ with $\Abs{\Jumps'}=\Abs{\Tiling''}$ such that for each $T\in\Tiling''$ we have $\Fkt{\Start}{T}\in\Jumps'$ ($\Fkt{\End}{T}\in\Jumps'$), and
			\item for each $i\in[1,h]$ the endpoint of $\Fkt{\Start}{T}$ (starting point of $\Fkt{\End}{T}$) lies in the component of $S-T_i$ that contains no tiles of $\Tiling''$ if $i=1$, or in the splice of $S$ separating $T_{i-1}$ and $T_i$.
		\end{itemize}
	\item[c)] there exists a family of tiles $\Tiling''\subseteq\Tiling'$ all contained in a single \hyperref[def:strip]{strip}  $S\subseteq W$ of height equal to the height of the tiles in $\Tiling$ such that $\Tiling''$ and $S$ meet the properties of outcome b) after taking the mirror image of $W$ along a vertical line.
	\end{enumerate}
\end{lemma}

By using the bounds obtained from the proofs in \cite{giannopoulou2020directed}, we get the following rough estimates for the functions $f_w$, $f_P$, and $f_W$:
\begin{enumerate}
	\item $\Fkt{f_w}{t}=2^9t^{10}$,
	\item $\Fkt{f_P}{t}=2^7t^8$, and
	\item $\Fkt{f_W}{t}=2^{32+t^{30}}$. 
\end{enumerate}

\paragraph{Swichting the Matching of a Wall}

While most of the routing necessary for our proofs will happen in the digraphic settings, mostly for convenience, matching theory will come in at several places.
The two most common techniques we are going to apply are changing the digraph we are working on by `flipping' the perfect matching along some directed cycles and tight cut contractions to remove large quantities of the (di)graphs we work on without changing their overall structure in an intrusive way.
Note that the later can also be done by applying directed tight cut contractions on the digraph directly.

An important part of this technique is the fact that we can use the `flipping' of the perfect matching along a directed cycle to essentially\footnote{Note that  the operation does more to the digraph than `just' reversing the direction of a cycle, it also affects all paths that enter or leave the cycle.} reverse its direction.

\begin{figure}[!ht]
	\centering

		};
		
		\node(Label1) [v:ghost,position=270:50mm from W1,align=center] {A slice $W$ of width $4$\\ of a cylindrical $8$-wall.};
		\node(Label2) [v:ghost,position=270:50mm from W2] {$B\coloneqq\Split{W}$ with the canonical matching $M$.};
		\node(Label3) [v:ghost,position=270:50mm from W3] {$B$ with the mixed matching $\Mix{W}$.};
		\node(Label4) [v:ghost,position=270:50mm from W4] {$\AB{W}$ and the twin walls of $W$.};
		
	\end{tikzpicture}
	\caption{The `flipping' of a cylindrical $4$-wall into two overlapping cylindrical $2$-walls that are oppositely oriented.}
	\label{fig:flippingawall}
\end{figure}

Let us briefly investigate the operation of switching the perfect matching $M$ along a horizontal cycle of a \hyperref[def:matchingwall]{matching wall} and its effect on the resulting $M'$-direction.

Let $k\in\N$ be a positive integer and $W$ be a cylindrical $4k$ wall, as well as $W'$ be a \hyperref[def:slice]{slice} of width $2k$ of $W$.
Let us number the vertical cycles $W'$ inherits from $W$ as $Q_1,\dots,Q_{2k}$ and let us write $\hat{P}_1,\dots,\hat{P}_{4k}$ for the subpaths of the vertical paths of $W$ that still are present in $W'$.

Then $B_{W'}\coloneqq\Split{W'}$ is a slice of width $k$ of a \hyperref[def:matchingwall]{matching $2k$-wall}, and we may assume $M$ to be its canonical matching.
Hence $\DirM{B_{W'}}{M}=W'$.
For each $i\in[1,2k]$ let $C_i$ be the $M$-conformal cycle of $B_{W'}$ such that $\DirM{C_i}{M}=Q_i$.
We define the \emph{mixed matching} of $B_{W'}$, denoted by $\Mix{W'}$, as
\begin{align*}
	\Mix{W'}\coloneqq M\Delta \bigcup_{\substack{i\in[1,2k]\\i\text{ odd}}}\E{C_i}.
\end{align*}
Thus $\Mix{W'}$ is obtained from $M$ by switching $M$ along every second vertical cycle starting with the first.
Let us denote by $\AB{W'}$ the digraph $\DirM{B_{W'}}{\Mix{W'}}$.
We say that $\A{W'}$ and $\B{W'}$ are the \emph{twin walls} of $W'$ as illustrated in \cref{fig:flippingawall}.

For the construction of $\A{W'}$ and $\B{W'}$, $W'$ must necessarily be of even width.
However, $\AB{W'}$ can be constructed for any slice, so in slight abuse of notation we will use $\AB{W'}$ to describe the digraph obtained from $\Split{W'}$ by switching the canonical matching along every second vertical cycle starting with the first.

Note that $\AB{W'}$ contains two cylindrical $k$-walls, one using the now flipped versions of the odd $Q_i$, denoted by $\A{W'}$, and the other using the still intact even $Q_i$, which is denoted by $\B{W'}$.
Moreover, if we start out with our embedding of $W$ and keep this embedding during all transformations, the vertical cycles of $\A{W'}$ go towards the top, while the vertical cycles of $\B{W'}$ still go to the bottom.
See \cref{fig:flippingawall} for an illustration.
Moreover, every path connecting one of the vertical cycles of $\CC{W'}$ to another one, for some $\CSymb\in\Set{\ASymb,\BSymb}$, must necessarily visit all vertical cycles of its twin that lie in between.

This construction allows us to move up- and downwards almost arbitrarily in a sufficiently large wall, a huge advantage ober the more rigid digraphic setting.
In a slice $W'$ of the cylindrical wall, every directed cycle must visit all of the vertical paths, but in $\AB{W'}$ we are able to find directed cycles locally, which means there are strongly connected subgraphs of $\AB{W'}$ which lie within some tiles of $W'$.

\paragraph{The Model of $\Bidirected{K_t}$}
Let us describe how the construction of the $\Bidirected{K_t}$-butterfly minor works.
We utilise essentially the same construction that was used in \cite{giannopoulou2020directed} to obtain \cref{lemma:directedlongjumps}.

Consider the slice $\widetilde{U}_1$ of width $d_1$ from the triadic partition of $U$.
Partition $\widetilde{U}_1$ into $t+1$ slices of equal width named $U'_1,\dots,U'_{t+1}$.
The first $t$ of these slices will hold the \emph{roots} of the $t$ vertices of $\Bidirected{K_t}$, where the root of the $i$th vertex will be contained in $U'_i$.
The root of the $i$th vertex will be a path $V_i$, subpath of some horizontal path of $W_i'$ from left to right, with $t$ incoming edges $e^i_j$ and $t$ outgoing edges $f^i_j$ in such a way that
\begin{itemize}
	\item the edges $e^i_j$ and $f^i_j$, $j\in[1,t-1]$, belong to vertical cycles of $U_i$, $i\in[1,t]$,
	\item the heads of $e^i_1,e^i_2,\dots,e^i_t$ appear in the order listed when traversing along $V_i$,
	\item the tails of $f^i_1,f^i_2,\dots,f^i_t$ appear in the order listed when traversing along $V_i$, and
	\item the head of $e^i_t$ appears on $V_i$ before the tail of $f^i_1$ when traversing along $V_i$.
\end{itemize}
We also require the slices $U'_i$ to appear in the order $U_1',U_2',\dots,U_t'$ from left to right, and that all $V_i$ belong to the same horizontal path.
If we are able to find a family of pairwise disjoint paths, internally disjoint from the roots, such that for every pair of distinct $i,j\in[1,t]$ the head of $f^i_j$ is linked to the tail of $e^j_i$, then the union of these paths together with the roots can be seen to form a butterfly minor model of $\Bidirected{K_t}$.
Indeed, each root can be contracted into a single vertex by the butterfly minor relation.
If we do this for each of the $t$ roots, then the resulting digraph is a subdivision of $\Bidirected{K_t}$.
For an illustration of a root see \cref{fig:aroot}.

\begin{figure}[!ht]
	\centering
	\begin{tikzpicture}[scale=0.8]
		\pgfdeclarelayer{background}
		\pgfdeclarelayer{foreground}
		\pgfsetlayers{background,main,foreground}
		
		\node(v1) [v:main] {};
		\node(v2) [v:main,position=0:12mm from v1] {};
		\node(v3) [v:main,position=0:12mm from v2] {};
		\node(v4) [v:main,position=0:12mm from v3] {};
		\node(v5) [v:main,position=0:12mm from v4] {};
		\node(v6) [v:main,position=0:12mm from v5] {};
		\node(v7) [v:main,position=0:12mm from v6] {};
		\node(v8) [v:main,position=0:12mm from v7] {};
		
		\node(e1) [v:main,position=90:12mm from v1] {};
		\node(e1L) [v:ghost,position=180:4.5mm from e1] {$e^i_1$};
		
		\node(e2) [v:main,position=90:12mm from v2] {};
		\node(e2L) [v:ghost,position=180:4.5mm from e2] {$e^i_2$};
		
		\node(e3) [v:main,position=90:12mm from v3] {};
		\node(e3L) [v:ghost,position=180:4.5mm from e3] {$e^i_3$};
		
		\node(e4) [v:main,position=90:12mm from v4] {};
		\node(e4L) [v:ghost,position=180:4.5mm from e4] {$e^i_4$};
		
		\node(f1) [v:main,position=270:12mm from v5] {};
		\node(f1L) [v:ghost,position=0:4.5mm from f1] {$f^i_1$};
		
		\node(f2) [v:main,position=270:12mm from v6] {};
		\node(f2L) [v:ghost,position=0:4.5mm from f2] {$f^i_2$};
		
		\node(f3) [v:main,position=270:12mm from v7] {};
		\node(f3L) [v:ghost,position=0:4.5mm from f3] {$f^i_3$};
		
		\node(f4) [v:main,position=270:12mm from v8] {};
		\node(f4L) [v:ghost,position=0:4.5mm from f4] {$f^i_4$};
		
		\begin{pgfonlayer}{background}
			
			\draw[e:main,->] (v1) to (v2);
			\draw[e:main,->] (v2) to (v3);
			\draw[e:main,->] (v3) to (v4);
			\draw[e:main,->] (v4) to (v5);
			\draw[e:main,->] (v5) to (v6);
			\draw[e:main,->] (v5) to (v7);
			\draw[e:main,->] (v7) to (v8);
			
			\draw[e:main,->] (e1) to (v1);
			\draw[e:main,->] (e2) to (v2);
			\draw[e:main,->] (e3) to (v3);
			\draw[e:main,->] (e4) to (v4);
			
			\draw[e:main,->] (v5) to (f1);
			\draw[e:main,->] (v6) to (f2);
			\draw[e:main,->] (v7) to (f3);
			\draw[e:main,->] (v8) to (f4); 
			
		\end{pgfonlayer}
	\end{tikzpicture}
	\caption{The the root of the model of the $i$th vertex of $\Bidirected{K_4}$.}
	\label{fig:aroot}
\end{figure}

\paragraph{A Refinement of \Cref{lemma:directedlongjumps}}

In this paragraph we aim to get rid of the second and third outcome of \cref{lemma:directedlongjumps}.
These outcomes produce a set of semi local long jumps, meaning that these long jumps originate in the same strip and have their other endpoint almost immediately to the left or the right of the tiles where they originate.
For the digraphic setting this poses a problem since this kind of crosses cannot be used to create a $\Bidirected{K_t}$-butterfly minor model.
Hence the bounded number of `cross rows` in the definition of barely flatness.
In \cite{giannopoulou2020directed} this is overcome by a different argument highly increasing the necessary functions.
However, In our case we may use the additional freedom of the matching setting and thus the idea of twin walls to enable us to use to gu up and down in the wall to make use of all these crosses.

\begin{lemma}\label{lemma:directedlongjumpsNew}
	There exist functions $f_w\colon\N\rightarrow\N$, $f_P\colon\N\rightarrow\N$, and $f_W\colon\N\rightarrow\N$ such that for every $t\in\N$ the following holds:
	Let 
	\begin{itemize}
		\item $D$ be a digraph,
		\item $W$ be a cylindrical $3k$-wall with $k\geq\Fkt{f_W}{t}$ in $D$,
		\item $\Triade=\Brace{W,k,W_1,W_2,W_3,W^1,W^2,W^3}$ be the \hyperref[def:triadicpartition]{triadic partition} of $W$, and
		\item \hyperref[def:tiling]{$\Tiling=\Tiling_{W,k,\Fkt{f_w}{t}.\xi,\xi'}$} for some $\xi,\xi'\in[1,\Fkt{f_w}{t}+1]$.
	\end{itemize}
	If there exists a subfamily $\Tiling'$ of $\Tiling$ and a family $\Jumps$ of pairwise disjoint directed paths in $D$ with the following properties:
	\begin{enumerate}
		\item Every member of $\Jumps$ is internally disjoint from $W$ but has both endpoints on $W$,
		\item $\Abs{\Tiling'}=\Abs{\Jumps}=\Fkt{f_P}{t}$,
		\item for every $T_{\Fkt{\ColumnFunction}{p},\Fkt{\RowFunction}{q},\Fkt{f_w}{t}}\neq T_{\Fkt{\ColumnFunction}{p'},\Fkt{\RowFunction}{q'},\Fkt{f_w}{t}}\in\Tiling'$ we have $\max\Set{\Abs{p-p'},\Abs{q-q'}}\geq 2$,
		\item there exists a bijection $\Start\colon\Tiling'\rightarrow\Jumps$ ($\End\colon\Tiling'\rightarrow\Jumps$) such that the starting point (endpoint) of the path $\Fkt{\Start}{T}$ ($\Fkt{\End}{T}$) belongs to the centre of $T$,
		\item $\V{\Fkt{\Start}{T}}\cap\V{\Tiling'}$ ($\V{\Fkt{\End}{T}}\cap\V{\Tiling'}$) contains exactly the endpoint of $\Fkt{\Start}{T}$ ($\Fkt{\End}{T}$) where $\V{\Tiling'}=\bigcup_{T'\in\Tiling'}\V{T'}$, and finally
		\item the endpoints (starting points) of the paths in $\Jumps$ are of mutual $W$-distance at least $4$. 
	\end{enumerate}
	Then $\Split{D}$ has a $K_{t,t}$-\hyperref[def:matchingminor]{matching minor} grasped by $\Split{W}$.
\end{lemma}

\begin{proof}
We may assume that case b) or c) of \cref{lemma:directedlongjumps} holds, as in case a) we are done by \cref{cor:Mmodels}. 
Moreover, we may assume to be in case b) since case c) can be handled analogously.
Let $S\subseteq W$ be the \hyperref[def:strip]{strip}, $\mathcal{T}''\subseteq \mathcal{T}'$, $\Abs{\mathcal{T}''}\geq\Fkt{f_P}{t}^{\frac{1}{4}}$, be the tiles, and $\mathcal{J}'\subseteq\mathcal{J}$, $\Abs{\mathcal{J}'}=\Abs{\mathcal{T}''}$, be the family of jumps from case b).
For each tile $T\in\mathcal{T}''$ let us define $s(T)$ to be the \hyperref[def:slice]{slice} of width $2\Fkt{f_P}{t}\leq \frac{1}{2}\Fkt{f_w}{t}$ such that each vertical path of $\Fkt{s}{T}$ intersects $T$ and $\Fkt{s}{T}$ contains the right most vertical path of $T$.
Please note that in $W-\Fkt{s}{T}$ there are still at least $2\Fkt{f_P}{t}$ vertical paths left that pass through $T$ to the right of the \hyperref[def:tile]{centre} of $T$.
Let us also define $\Fkt{c}{T}$ to be the \hyperref[def:slice]{slice} of $W$ of width $2\Fkt{f_P}{t}+2$ whose vertical paths all pass through $T$ and that has exactly $\Fkt{f_P}{t}+1$ vertical paths on either side of the centre of $T$.
So $\Fkt{c}{T}$ and $\Fkt{s}{T}$ are disjoint and whenever we have a family $\mathcal{P}$ of at most $\Fkt{f_P}{t}$ pairwise disjoint paths entering $T\cap\Fkt{c}{T}$ from the top we can select any path $P\in\mathcal{P}$ to be routed into the centre of $T$ while all other paths in $\mathcal{P}$ can be routed within $\Fkt{c}{T}$ to the bottom of $T$ while maintaining their disjointness.

Our goal is to make use of this property to achieve the following:
Suppose we have a family $\mathcal{P}$ of at most $\Fkt{f_P}{t}^{\frac{1}{4}}$ pairwise disjoint paths starting an the left perimeter of $W_2$.
We want to visit the tiles of $T\in \mathcal{T}''$ consecutively from left to right, each time sending some path of $\mathcal{P}$ into the \hyperref[def:tile]{centre} of $T$ and then into the corresponding jump from $\mathcal{J}'$.
If we manage to do this in such a way that in the end each path in $\mathcal{P}$ is assigned some long jump in $\mathcal{J}'$ and we can find a linkage from the endpoints of the jumps back to the left perimeter of $W_2$ without intersecting our previous routing, then we can use this to create any digraph on $t$ vertices as a butterfly minor by positioning the roots of the corresponding model within $W_1$ and attaching their in and out paths to the end- and start points of our routing on the left perimeter of $W_2$.
Hence this is enough to proof the lemma.

However, this cannot be done within the wall $W$ as it is because all vertical paths go downwards.
This means, to reach the next tile after we have visited some tile $T\in\mathcal{T}''$ we need to go all the way down to the bottom of $\Fkt{s}{T}$ and then use the flying edges of $W$ to reach the top part again.
By doing so it might become impossible for any path on the right of $\Fkt{s}{T}$ to reach the left perimeter of $W_2$ without crossing through the previously constructed routing.
To circumvent this problem we turn each $\Fkt{s}{T}$ into its two \hyperref[fig:flippingawall]{twin walls}.
This allows us to leave $\Fkt{c}{T}$ with our paths immediately after passing by the \hyperref[def:tile]{centre} of $T$ and continuing to the right until all remaining paths have reached their position in $\AB{\Fkt{s}{T}}$.
Here the remaining paths can now be routed \emph{upwards} within $\AB{\Fkt{s}{T}}$ to ensure they can enter the next tile from above, while its corresponding jump can be routed to the left without interference.
This strategy might change the perfect matching and thereby change the digraph $D$, however, we will still find a $\Bidirected{K_t}$-butterfly minor model grasped by the changed version of $W$.
Hence we will find a $K_{t,t}$-\hyperref[def:matchingminor]{matching minor} grasped by \hyperref[def:split]{$\Split{W}$} as desired.

We now describe how this routing is constructed.
Let $k\leq\Fkt{f_P}{t}$ be some positive integer and let $\pi$ be some permutation of $[1,k]$.
Furthermore, let us select two additional strips $S_1$ and $S_2$ of $W_2$, both of the same height as $S$, such that $S$ separates $S_1$ from $S_2$ in $W_2$ and $S_1$ is above $S$.
Let $s_1,\dots,s_k$ be the starting points of $k$ left-to-right paths, let us call them $P_1,\dots,P_k$, of $S_1$ such that for each $i\in[1,k]$ $s_i$ is the starting point of $P_i$ and for each $i\in[1,k-2]$ the path $P_{i+1}$ separates $P_i$ from $P_{i+2}$.
Additionally, let us assume the tiles of $\mathcal{T}''=\Set{T_1,T_2,\dots,T_{\Abs{\mathcal{T}''}}}$ to be numbered from left to right as in outcome b) of \cref{lemma:directedlongjumps}.
We start with $T_1$.
Visiting $T_1$ is done as follows:
In $\Fkt{c}{T_1}\cap T_1$ select $k-\Fkt{\pi}{1}$ distinct vertical paths $Q_k,\dots,Q_{k-\Fkt{\pi}{1}}$ occurring from left to right in the order listed such that $Q_{k-\Fkt{\pi}{1}}$ does not intersect the boundary of the centre of $T_1$ but still lies left of it.
Let $Q_{k-\Fkt{\pi}{1}-2},\dots,Q_1$ be another family of distinct vertical paths of $\Fkt{c}{T_1}\cap T_1$ occurring from left to right in the order listed, being disjoint from the centre of $T_1$ and occurring to the right of the centre $T_1$
Finally let $Q_{k-\Fkt{\pi}{1}-1}$ let be a shortest path in $D'$ starting in the start-vertex of the vertical path of $T_1$ that passes through the centre of $T_1$ and ending in the starting-point of the path $\Fkt{\Start}{T_1}$, where $D'$ is the digraph obtained from the union of $T_1$ and the subgraph of $D$ that is separated from $W$ by the boundary of the centre of $T_1$.
We may now find a linkage from $\Set{s_1,\dots,s_k}$ to the starting points of the $Q_i$ such that for each $i\in[1,k]$ the vertex $s_i$ is linked to $Q_i$ and the linkage stays entirely in the quadrant of $W_2$ defined by the upper path of $S$ and the left most path of $\Fkt{s}{T_1}$.
Moreover, we may choose this linkage to stay as much as possible within the strip $S_1$.
We may now continue with each of the paths routed so far until for each $i\in[1,k]\setminus\Set{\Fkt{\pi}{1}}$ the path with starting point $s_i$ meets $P_i$ again.
Let $W_2^2$ be the largest slice of $W_2$ containing the right perimeter of $W_2$ while being disjoint from $T_1$.
Note that
From here on, each of these paths may traverse along its respective $P_i$ to the right until it meets a vertex of the left perimeter of $W_2^2$ for the first time.

Note that we may now find a bijection $\varphi_1\colon[1,k-1]\rightarrow[1,k]\setminus\Set{1}$ such that $\Fkt{\varphi}{i}=i+1$ and a permutation $\pi_1$ of $[1,k-1]$ such that $\Fkt{\pi_1}{i}=\Fkt{\pi}{\Fkt{\varphi_1}{i}}$ for all $i\in[1,k-1]$.
This allows us to proceed with the construction of this part of the routing via induction.
In the end, for each $i\in[1,k]$ we have sent the $\Fkt{\pi}{i}$th path from the centre of $T_i$ to the left.
Let now $L_1,\dots,L_k$ be $k$ distinct right-to-left paths in $S_2$ occurring from top to bottom in the order listed.
For each $i\in[1,k]$ we may now find a path from the endpoint of $\Fkt{\Start}{T_i}$ to $L_i$ by moving along a single vertical path of $W_2$.
It follows from the construction in the previous step that each of these paths is disjoint from all parts of the routing so far.
We may now complete the second step by enhancing each of the paths build so far by following along the $L_i$ to their respective endpoints on the left perimeter of $W_2$.
Let $t_1,\dots,t_k$ be these endpoints such that $t_i$ is the endpoint of $L_i$ for all $i\in[1,k]$.
We have now found $k$ pairwise disjoint paths $R_1,\dots,R_k$ such that for each $i\in[1,k]$ $R_i$ starts in $s_i$ and ends in $t_{\Fkt{\pi}{i}}$.

To complete the proof let us find $t$ pairwise disjoint \hyperref[def:slice]{slices} of $W_1$, each of width $2\Fkt{f_w}{t}+2$ and let $U_1,\dots,U_t$ be these slices numbered according to their occurrence in $W_1$ from left to right.
For each $i\in[1,t]$ let us create a root $F_i$ of $\Bidirected{K_t}$ such that all roots together fit into a strip $S'$ of height $2\Fkt{f_w}{t}+2$ that can be combined with the strip $S$ to form a strip of $W_1\cup W_2$.
Moreover, let us create the root slightly different to how it was defined above:
Let the base path go from right to left, let the outgoing paths leave to the top on the left and let the incoming paths enter from the bottom on the right.
To make this possible we may create the $F_i$ in a way such that they are also all contained in a single slice of width $2\Fkt{f_w}{t}+2$ of $W_1$ and, moreover, such that the outgoing parts all belong to a subslice of with at most $\Fkt{f_w}{t}$.
For this subslice we switch the matching of the split of $W_1$ to form the twin wall.
Now the required upwards paths exist.
Let $S_1'$ be the strip of $W_1$ such that $S_1'\cup S_1$ is a strip of $W_1\cup W_2$.
We may now find a linkage from $\CondSet{f^j_i}{i\in[1,t]\text{, and }j\in[1,t]\setminus\Set{i}}$ to the set of starting points of left-to-right paths of $S_1$ on the left perimeter of $W_2$.
Similarly we may connect the vertices $t_i$ from above to the vertices in $\CondSet{e^j_i}{i\in[1,t]\text{, and }j\in[1,t]\setminus\Set{i}}$.
As proven above we may create any particular linkage of the $s_i$ to the $t_i$ within $W_2$ together with the jumps $\mathcal{J}'$ and thus we are able to construct a butterfly minor model of $\Bidirected{K_t}$ in some $M$-direction of $\Split{D}$.
\end{proof}

\paragraph{Further refining \Cref{lemma:directedlongjumpsNew} and removing long jumps}

\Cref{lemma:directedlongjumpsNew} is already a powerful tool.
However, it is not straightforward how to obtain the apex set $A$ just from its application.
Hence in the following, we aim for further refinement and the existence of the set $A$.

Let $f_w\colon\N\rightarrow\N$ be some function, $t,k\in\N$ two positive integers, and $\xi,\xi'\in[1,\Fkt{f_w}{t}+1]$.
Moreover, let $W=\Brace{Q_1,\dots,Q_{3k},\hat{P}_1,\dots,\hat{P}_{3k}}$ be a cylindrical $3k$-wall with its \hyperref[def:triadicpartition]{triadic partition} $\Triade=\Brace{W,k,W_1,W_2,W_3,W^1,W^2,W^3}$ and \hyperref[def:tiling]{$\Tiling=\Tiling_{W,k,\Fkt{f_w}{t},\xi,\xi'}$}.
A \emph{colouring} of $\Tiling$ is a partition of $\Tiling$ into four classes, namely $\Class_1,\Class_2,\Class_3$, and $\Class_4$ as follows:
For every $i\in\left[1,\Ceil{\frac{k+\xi-1}{2\Fkt{f_w}{t}+1}}+1\right]$ and every $j\in \left[1,\Ceil{\frac{3k-\xi'-1}{2\Fkt{f_w}{t}+1}}+1\right]$ we assign to $T_{\Fkt{\ColumnFunction}{i},\Fkt{\RowFunction}{j},\Fkt{f_w}{t}}$ the colour $\Brace{i\mod 2}+2\Brace{j\mod 2}+1$.
This means that tiles where $\Fkt{\ColumnFunction}{i}$ and $\Fkt{\RowFunction}{j}$ are even get colour $1$, the tiles where $\Fkt{\RowFunction}{j}$ is even but $\Fkt{\ColumnFunction}{i}$ is odd get $3$, and so on.
Hence every column is two-chromatic, every row is so as well, and between each pair of tiles from the same colour that share a row or a column, there is a tile of a different colour that separates those tiles in their respective row or column.
Additionally, if $T$ is some tile, then the eight tiles surrounding $T$ are all of different colour than $T$ itself.

\begin{definition}[Auxiliary Digraph Type I]\label{def:auxiliarydigraphI}
	Let $t,k,k',w\in\N$ be positive integers such that $k\geq k'\geq2\Fkt{f_W}{t}+4\Fkt{f_P}{t}\Brace{2w+1}$, $w\geq 2\Fkt{f_w}{t}$, and $\xi,\xi'\in\left[1,w+1\right]$.
	Let $D$ be a digraph containing a cylindrical $3k$-wall $W=\Brace{Q_1,\dots,Q_{3k},\hat{P}_1,\dots,\hat{P}_{3k}}$ with its \hyperref[def:triadicpartition]{triadic partition} $\Triade=\Brace{W,k,W_1,W_2,W_3,W^1,W^2,W^3}$, and a tiling \hyperref[def:tiling]{$\Tiling=\Tiling_{W,k,w,\xi,\xi'}$}.
	Let $i\in[1,4]$, $\Set{\Class_1,\dots,\Class_4}$ be a four colouring of $\Tiling$ and $W'\subseteq W$ be a \hyperref[def:slice]{slice} of width $k'$ of $W_2$.
	At last, let us denote by $\Tiling'$ the family of tiles from $\Tiling$ that share a vertex with $W'$.
	Similarly let $\Class_i'\coloneqq\Tiling'\cap\Class_i$.
	Then $\Fkt{D^1_i}{W'}$ is the digraph obtained from $D$ by performing the following construction steps for every $T\in\Class_i'$:
	\begin{enumerate}
		\item add new vertices $x_T^{\text{in}}$ and $x_T^{\text{out}}$,
		\item for every vertex $u$ in the centre of $T$ introduce the edges $\Brace{u,x_T^{\text{in}}}$ and $\Brace{x_T^{\text{out}},u}$, and then
		\item delete all internal vertices of $T$.
	\end{enumerate}
\end{definition}

We also need the following result:
Let $D$ be a digraph and $X,Y\subseteq\V{D}$. A \emph{half-integral $X$-$Y$-linkage} of order $k$ is a family $\mathcal{P}$ of directed $X$-$Y$-paths such that every vertex of $D$ is contained in at most two paths from $\mathcal{P}$.
By $\V{\mathcal{P}}$ we denote the set $\bigcup_{P\in\mathcal{P}}\V{P}$.

\begin{theorem}[\cite{giannopoulou2020directed}]\label{thm:halfintegral}
	Let $k\in\N$ be a positive integer, $D$ be a digraph, and $X,Y\subseteq\V{D}$.
	If $\mathcal{P}$ is a half-integral $X$-$Y$-linkage of order $2k$ in $D$, then there exists a family $\mathcal{J}$ of pairwise disjoint $X$-$Y$-paths such that $\V{\mathcal{J}}\subseteq\V{\mathcal{P}}$.
\end{theorem}

\begin{lemma}\label{lemma:longjumps1}
	Let $t,k,k',w\in\N$ be positive integers such that $k\geq k'\geq2\Fkt{f_W}{t}+2^{16}\Fkt{f_P}{t}+2$, $w\geq 2\Fkt{f_w}{t}+2^7\Fkt{f_P}{t}$, and $\xi,\xi'\in\left[1,w+1\right]$.
	Let $D$ be a digraph containing a cylindrical $3k$-wall $W=\Brace{Q_1,\dots,Q_{3k},\hat{P}_1,\dots,\hat{P}_{3k}}$ with its \hyperref[def:triadicpartition]{triadic partition} $\Triade=\Brace{W,k,W_1,W_2,W_3,W^1,W^2,W^3}$, and a tiling \hyperref[def:tiling]{$\Tiling=\Tiling_{W,k,w,\xi,\xi'}$}.
	Let $i\in[1,4]$, $\Set{\Class_1,\dots,\Class_4}$ be a four colouring of $\Tiling$ and $W'\subseteq W$ be a \hyperref[def:slice]{slice} of width $k'$ of $W_2$.
	
	Now let $\Tiling'$ be the family of all tiles of $\Tiling$ that are completely contained in $W'$ and let $\widetilde{W}$ be the smallest slice of $W$ that contains all tiles from $\Tiling'$.
	
	Consider the \hyperref[def:auxiliarydigraphI]{auxiliary digraph of type I} $\Fkt{D_i^1}{\widetilde{W}}$ and let $\Class_i'$ be as in the definition of $\Fkt{D_i^1}{\widetilde{W}}$.
	Define the sets 
	\begin{align*}
		X_{\text{I}}^{\text{out}}&\coloneqq\CondSet{x_T^{\text{out}}}{T\in\Class_i'}\text{, and}\\
		X_{\text{I}}^{\text{in}}&\coloneqq\CondSet{x_T^{\text{in}}}{T\in\Class_i'}.
	\end{align*}
	Additionally, we construct the set $Y_{\text{I}}$ as follows:
	Let $Q$ and $Q'$ be the two cycles of $\Perimeter{W_2}$.
	For every $j\in\left[ 1,\frac{3k}{4} \right]$, $Y_{\text{I}}$ contains exactly one vertex of $Q\cap P^1_{4j}$, $Q\cap P^2_{4j+2}$, $Q'\cap P^1_{4j}$, and $Q'\cap P^2_{4j+2}$ each.
	
	If there exists a family $\mathcal{L}$ of pairwise disjoint directed paths with $\Abs{\mathcal{L}}= 2^7\Fkt{f_P}{t}$ such that either
	\begin{itemize}
		\item $\mathcal{L}$ is a family of directed $X_{\text{I}}^{\text{out}}$-$Y_{\text{I}}$-paths, or
		\item $\mathcal{L}$ is a family of directed $Y_{\text{I}}$-$X_{\text{I}}^{\text{in}}$-paths,
	\end{itemize}
	then $\Split{D}$ has a $K_{t,t}$-\hyperref[def:matchingminor]{matching minor} grasped by \hyperref[def:split]{$\Split{W}$}.
\end{lemma}

\begin{proof}
	The proof is divided into several steps and we start with a brief outline.
	Our goal is to construct a cylindrical wall $U\subseteq W$ of sufficient size, together with a family of $\Fkt{f_P}{t}$ directed $U$-paths that meet the requirements of \cref{lemma:directedlongjumpsNew}.
	If we are able to do this, then \cref{lemma:directedlongjumpsNew} yields the desired outcome.
	
	Without loss of generality, let us assume $\mathcal{L}$ is a family of directed $X_{\text{I}}^{\text{out}}$-$Y_{\text{I}}$-paths.
	The other case can be seen using similar arguments.
	
	Let $L$ and $P$ be directed paths.
	We say that $P$ is a \emph{long jump of $L$} if $P$ is a $w$-long jump over $W$ and $P\subseteq L$.
	We also say that $P$ is a \emph{jump of $L$} if $P$ is a directed $W$-path.
	
	Towards our goal, we first show that we can use $\mathcal{L}$ to construct a half-integral $X_{\text{I}}^{\text{out}}$-$Y_{\text{I}}$-linking $\mathcal{L}_1$ such that
	\begin{enumerate}
		\item $\Abs{\mathcal{L}_1}=2^7\Fkt{f_P}{t}$,
		\item there exists a family $\mathcal{F}\subseteq\Tiling'$ with $\Abs{\mathcal{F}}\leq 2^7\Fkt{f_P}{t}$, and
		\item for every $L\in\mathcal{L}_1$, every endpoint $u$ of a jump of $L$ with $u\in\V{\widetilde{W}}$ belongs to a tile from $\Class_i'\cup\mathcal{F}$.
	\end{enumerate}
	Once this is achieved, we use \cref{thm:halfintegral} to obtain a family $\mathcal{L}_2$ of pairwise disjoint directed $X_{\text{I}}^{\text{out}}$-$Y_{\text{I}}$-paths of size $2^6\Fkt{f_P}{t}$ from $\mathcal{L}_1$.
	Afterwards, we remove the cycles and paths of $\widetilde{W}$ that meet tiles from $\mathcal{F}$ and obtain a new slice $\widetilde{W}'$ of some cylindrical wall.
	For this slice, we construct a \hyperref[def:tiling]{tiling} and a \hyperref[def:tilingII]{tier II tiling} as well as a half-integral linking $\mathcal{L}_4$ of size $2^6\Fkt{f_P}{t}$ from $\mathcal{L}_3$ that connects the \hyperref[def:tile]{centres} of some tiles in the tier II tiling to vertices of $\widetilde{W}'$ such that their endpoints are mutually far enough apart and every path in $\mathcal{L}_4$ is internally disjoint from the new wall.
	Another application of \cref{thm:halfintegral} then yields the family of long jumps necessary for an application of \cref{lemma:directedlongjumpsNew}.
	
	We start out with the construction of $\mathcal{L}_1$ and $\mathcal{F}$.
	For this, let $\mathcal{L}'\coloneqq \mathcal{L}$, $\mathcal{L}_1\coloneqq \emptyset$, and $\mathcal{F}\coloneqq\emptyset$.
	As long as $\mathcal{L}'$ is non-empty, perform the following actions:
	
	Select some path $L\in\mathcal{L}$.
	In case $L$ is internally disjoint from $\widetilde{W}$, add $L$ to $\mathcal{L}_1$ and remove it from $\mathcal{L}'$.
	Otherwise let $s_L$ be its starting point and let $v_L$ be the first vertex of $L$, when traversing along $L$ starting from $s_L$, that belongs to $\widetilde{W}$, but not to a tile from $\Class_i'$.
	\begin{enumerate}
		\item If $v_L$ does not belong to a tile from $\mathcal{F}$, let $T\in\Tiling\setminus\Class_i$ be the tile that contains $v_L$ and add $T$ to $\mathcal{F}$.
		Let $R$ be a shortest directed path from $v_L$ to $Y_{\text{I}}$ in $W$ such that $R$ avoids all vertices of $W$ that are contained in two different paths of $\mathcal{L}_1$ and that is internally disjoint from $Lv_L$.
		Now add $Lv_LR$ to $\mathcal{L}_1$ and remove $L$ from $\mathcal{L}'$.
		Note that such a path $R$ must exist since the paths in $\mathcal{L}$ are pairwise disjoint, we never used $T$ for such a re-routing before, and $w$ and $k'$ are chosen sufficiently large in proportion to $2^7\Fkt{f_P}{t}$.
		Also note that the path $R$ is exactly the part, where we might go from integral to half-integral, but since our paths were pairwise disjoint to begin with, we can be sure that $R$ does never meet a vertex contained in two distinct paths.
		
		\item So now suppose $v_L$ belongs to a tile $T$ from $\mathcal{F}$.
		Let us follow along $v_LL$ until the first time we encounter a vertex $u_L$ for which one of the following is true:
		\begin{enumerate}
			\item $u_L$ belongs to a tile $T$ from $\Tiling\setminus\Brace{\Class_i\cup\mathcal{F}}$, or
			\item every internal vertex of $u_LL$ belongs to $W-\widetilde{W}$ or to some tile from $\Class_i\cup\mathcal{F}$. 
		\end{enumerate}
		If a) is the case, repeat the instruction from i) but replace $v_L$ by $u_L$.
		In this case $T$ is added to $\mathcal{F}$.
		Otherwise b) must hold and here we may simply remove $L$ from $\mathcal{L}'$ and add it to $\mathcal{L}_1$.
	\end{enumerate}
	Now for every $L\in\mathcal{L}$ we added at most one tile to $\mathcal{F}$ and thus $\Abs{\mathcal{F}}\leq\Abs{\mathcal{L}}$.
	Moreover, from the construction it is clear that $\mathcal{L}_1$ is indeed a half-integral linkage from $X_{\text{I}}^{\text{out}}$ to $Y_{\text{I}}$.
	Also, please note that we may assume that every $L$ meets each tile in $\mathcal{F}$ in at most $2^7\Fkt{f_P}{t}+1$ horizontal path pairs and vertical cycles, since otherwise one could find a short cut through $W$ itself.
	
	Next, we may apply \cref{thm:halfintegral} to obtain a family $\mathcal{L}_2$ of pairwise disjoint directed $X_{\text{I}}^{\text{out}}$-$Y_{\text{I}}$-paths with $\V{\mathcal{L}_2}\subseteq\V{\mathcal{L}_1}$ and $\Abs{\mathcal{L}_2}=2^6\Fkt{f_P}{t}$.
	This completes the second step.
	
	For the third step, let us consider $W''\coloneqq \InducedSubgraph{\widetilde{W}}{\Tiling,i}$ together with the tiling $\Tiling''\coloneqq\TierIITiling{\Tiling,i,f}{\widetilde{W}}$ and a four-colouring $\Set{\widetilde{\Class}_1,\dots,\widetilde{\Class}_4}$.
	Note that by choice of $k'$ this means that $W''$ is a slice of width $k''\geq1\Fkt{f_W}{t}+2^7\Fkt{f_P}{t}\Brace{2w+1}+1$ of some cylindrical $3k''$-wall that is completely contained in $W$.
	For each $L\in\mathcal{L}_2$ let $T^1_L\in \Class_i$ such that the starting point $s_L$ of $L$ belongs to $T^1_L$.
	Let $K^1_L\in \Tiling''$ be the tile whose centre is the perimeter of $T^1_L$.
	Choose any vertex $s'_L$ of degree three in $W''$ that is not contained in any path of $\mathcal{L}_2$, and let $R_L$ be a directed path from $s_L$ to $s'_L$ within $T^1_L$.
	Let $\mathcal{L}_3'$ be the resulting, and potentially now again half-integral, family of directed paths.
	Now there must exist $j\in[1,4]$ such that at least $2^4\Fkt{f_P}{t}$ of the paths from $\mathcal{L}_3'$ start at the centre of a tile from $\widetilde{\Class}_j$.
	Let $\mathcal{L}_3''\subseteq\mathcal{L}_3'$ be a family of exactly $2^4\Fkt{f_P}{t}$ such paths.
	Next let us consider the family $\mathcal{F}$.
	Let $W'''$ be the subgraph of $W''$ induced by all vertical cycles and horizontal path pairs in $W''$ that do not contain a vertex of some tile in $\mathcal{F}$ that belongs to a path in $\mathcal{L}_3''$.
	Since $\Abs{\mathcal{F}}\leq 2^7\Fkt{f_P}{t}$ and each tile in $\mathcal{F}$ meets a path in $\mathcal{L}_3''$ in at most $2^7\Fkt{f_P}{t}+1$ such cycles and pairs of horizontal paths, it follows that $W'''$ is a slice of width $k'''\geq\Fkt{f_W}{t}+2$ of some cylindrical $3k'''$wall $W^*\subseteq W$.
	Moreover, $W^*$ can be partitioned into three \hyperref[def:slice]{slices} of width $k'''$ as in its \hyperref[def:triadicpartition]{triadic partition}, such that $W'''$ is the slice in the middle.
	Let us rename the paths and cycles of $W^*$ such that $W^*=\Brace{Q^*_1,\dots,Q^*_{3k'''},\hat{P^*}_1,\dots,\hat{P^*}_{3k'''}}$, and we construct the set $Y^*$ as follows:
	Let $Q^*$ and $'Q^*$ be the two cycles of $\Perimeter{W'''}$.
	For every $j\in\left[ 1,\frac{3k'''}{4} \right]$, $Y^*$ contains exactly one vertex of $Q^*\cap P^{*1}_{4j}$, $Q\cap P^{*2}_{4j+2}$, $'Q^*\cap P^{*1}_{4j}$, and $'Q^*\cap P^{*2}_{4j+2}$ each.
	Let $L\in\mathcal{L}_3''$ be any path and $t_L$ be the first vertex after its starting point $L$ shares with either $W'''$ or $W^*-W'''$.
	In case $t_L\in\V{W'''}$, simply add $Lt_L$ to $\mathcal{L}_3'''$.
	Otherwise, let $b_L$ be the endpoint of $L$ in $W^*-W'''$.
	Then we can find a path $R_L$ in $W$ from $b_L$ to a vertex $t^*_L$ of $Y^*$ such that $t^*_L$ is of $W^*$-distance at least $4$ to every endpoint of every path already in $\mathcal{L}_3'''$, $R_L$ is internally disjoint from $L$, and $R_L$ does not contain a vertex that is contained in two distinct paths from $\mathcal{L}_3'''$.
	Add $LR_L$ to $\mathcal{L}_3'''$.
	Finally, $\mathcal{L}_3'''$ is a half-integral linkage from the set starting points $S^*$ of the paths in $\mathcal{L}_3''$ to $Y^*$ of size $2^4\Fkt{f_P}{t}$, and thus by \cref{thm:halfintegral} we can find a family $\mathcal{L}_4$ of pairwise disjoint directed paths from $S^*$ to $Y^*$ with $\V{\mathcal{L}_4}\subseteq\V{\mathcal{L}_3'''}$ that is of size $2^3\Fkt{f_P}{t}$.
	It follows that all paths in $\mathcal{L}_4$ are internally disjoint from $W'''$.
	This concludes the fourth step.
	
	Let us consider the tiles of $\widetilde{\Class}_i$ whose centres contain a vertex of $S^*$.
	Since $W'''$ might be a proper subgraph of $W''$, $\Tiling''$ is not necessarily a tiling of $W'''$.
	Each such tile $T$, however, contains a tile $T'$ of width $\Fkt{f_w}{t}$ with the same \hyperref[def:tile]{centre}.
	Since $T$ can be surrounded by at most $8$ tiles from $\mathcal{F}$ in $W'$, we may find, among the $2^3\Fkt{f_P}{t}$ many such tiles, a family $\Jumps$ of $\Fkt{f_P}{t}$ tiles that are pairwise disjoint and thus, since they all are constrcuted from the family $\widetilde{\Class}_i$, they meet the distance requirements of the tiles in \cref{lemma:directedlongjumpsNew}.
	Hence we may apply \cref{lemma:directedlongjumpsNew} and obtain a
	$K_{t,t}$-\hyperref[def:matchingminor]{matching minor} grasped by \hyperref[def:split]{$\Split{W^*}$} $\Bidirected{K_t}$-butterfly minor grasped by $W^*$.
	Moreover, since $W^*\subseteq W$, this completes the proof of our lemma.
\end{proof}

With this we are able to handle all long jumps that attach to tiles of different colour.
Using a second auxiliary digraph, we find a similar way to handle those long jumps over $W$ that attach to tiles of the same colour by using \cref{thm:xpaths}.

\begin{definition}[Auxiliary Digraph Type II]\label{def:auxiliarydigraphII}
	Let $t,k,k',w\in\N$ be positive integers such that $k\geq k'\geq2\Fkt{f_W}{t}$, $w\geq 2\Fkt{f_w}{t}$, and $\xi,\xi'\in\left[1,w+1\right]$.
	Let $D$ be a digraph containing a cylindrical $3k$-wall $W=\Brace{Q_1,\dots,Q_{3k},\hat{P}_1,\dots,\hat{P}_{3k}}$ with its \hyperref[def:triadicpartition]{triadic partition} $\Triade=\Brace{W,k,W_1,W_2,W_3,W^1,W^2,W^3}$, and a \hyperref[def:tiling]{tiling} $\Tiling=\Tiling_{W,k,w,\xi,\xi'}$.
	Let $i\in[1,4]$, $\Set{\Class_1,\dots,\Class_4}$ be a four colouring of $\Tiling$ and $W'\subseteq W$ be a \hyperref[def:slice]{slice} of width $k'$ of $W_2$ such that no \hyperref[def:tile]{tile} of $\Class_i$ contains a vertex of the perimeter of $W'$.
	Then $\Fkt{D^2_i}{W'}$ is the digraph obtained from $D$ by performing the following construction steps:
	
	for every $T\in\Class_i$, such that $T$ contains a vertex of $W'$, we do the following:
	\begin{enumerate}
		\item add a new vertex $x_T$, and
		\item for every vertex $v$ that belongs to the interior or the centre of $T$, introduce the edges $\Brace{x_T,v}$ and $\Brace{v,x_T}$.
	\end{enumerate}
	Once this is done, delete all vertices of $W'$ that do not belong to tiles of $\Class_i'$
	Let $X_{\text{II}}^i$ be the collection of all newly introduced vertices $x_T$.
\end{definition}

\begin{lemma}\label{lemma:longjumpspahse2}
	Let $t,k,k',w\in\N$ be positive integers, and $\xi,\xi'\in\left[1,w+1\right]$ where $w\geq2\Fkt{f_w}{t}$.
	Let $D$ be a digraph containing a cylindrical $3k$-wall $W_0=\Brace{Q_1,\dots,Q_{3k},\hat{P}_1,\dots,\hat{P}_{3k}}$, where $k\geq k'\geq4\Fkt{f_W}{t}^2$, with its \hyperref[def:triadicpartition]{triadic partition} $\Triade=\Brace{W_0,k,W_1,W_2,W_3,W^1,W^2,W^3}$, a \hyperref[def:slice]{slice} $W\subseteq W_2$ of width $k'$, a \hyperref[def:tiling]{tiling} $\Tiling=\Tiling_{W,k,w,\xi,\xi'}$, a four colouring $\Set{\Class_1,\dots,\Class_4}$, and a fixed colour $i\in[1,4]$.
	
	Then $D$ has a $K_{t,t}$-\hyperref[def:matchingminor]{matching minor} grasped by \hyperref[def:split]{$\Split{W_0}$}, or there exists a set $\MarkedTiles{2}_{i,\xi,\xi'}\subseteq\Tiling$ with $\Abs{\MarkedTiles{2}_{i,\xi,\xi'}}\leq 8\Fkt{f_P}{t}$ and a set $Z^2_{i,\xi,\xi'}\subseteq\V{D-W}$ with $\Abs{Z^2_{i,\xi,\xi'}}\leq 8\Fkt{f_P}{t}$ such that every directed $\V{W_0}$-path $P$ in $D-Z^2_{i,\xi,\xi'}$ whose endpoints belong to different tiles of $\Class_i$ contains a vertex of some tile in $\MarkedTiles{2}_{i,\xi,\xi'}$.
\end{lemma}

\begin{proof}
	Let $W'$ be the largest \hyperref[def:slice]{slice} of $W$ such that no tile of $\Class_i$ contains a vertex of $\Perimeter{W'}$.
	Let us consider the \hyperref[def:auxiliarydigraphI]{auxiliary digraph} $\Fkt{D^2_i}{W'}$ with the set $X_{\text{II}}^i$ of newly added vertices.
	By applying \cref{thm:xpaths} to the set $X_{\text{II}}^i$ in $\Fkt{D^2_i}{W'}$, we either find a set $Z$ of size at most $8\Fkt{f_P}{t}$ that hits all directed $X_{\text{II}}^i$-paths, or there exists a family $\mathcal{J}'$ of $4\Fkt{f_P}{t}$ pairwise disjoint directed $X_{\text{II}}^i$-paths in $\Fkt{D^2_i}{W'}$.
	
	Let us assume the latter.
	Then, by construction of $\Fkt{D^2_i}{W'}$, no path in $\mathcal{J}'$ contains a vertex of $W_2$.
	Back in the digraph $D$, let us consider the \hyperref[def:tilingII]{tier II tiling} $\Tiling''\coloneqq \TierIITiling{\Tiling,i,w}{W'}$ of width $w$ of $W''\coloneqq \InducedSubgraph{W'}{\Tiling,i}$.
	Notice that, by choice of $k$, $W''$ still contains a cylindrical $\Fkt{f_W}{t}$-wall $W'''$ such that the perimeter of every tile $T\in\Class_i$, for which $x_T$ is an endpoint of some path in $\mathcal{J}'$, bounds a cell of $W'''$.
	Let $\Tiling''$ be a tiling of $W'''$ such that the perimeter of every $T\in\Class_i$, for which $x_T$ is an endpoint of a path in $\mathcal{J}'$, is the centre of some tile in $\Tiling''$.
	We now consider a four colouring $\Set{\Class_1',\dots,\Class_4'}$ of $\Tiling''$.
	Then there must exist $j\in[1,4]$ and a family $\mathcal{J}''$ of size $\Fkt{f_P}{t}$ such that the starting points of every path in $\mathcal{J}''$ belongs to a tile of $\Class_i$ whose perimeter is the centre of a tile in $\Class_j'$.
	For every $J''\in\mathcal{J}''$ let $T_1,T_2\in\Class_i$ be the two tiles such that $J'$ is a directed $x_{T_1}$-$x_{T_2}$-path.
	We can now find a directed path $J$ that starts on the perimeter of of $T_1$, ends on the perimeter of $T_2$, and is internally disjoint from $W'''$.
	Hence we find a family $\mathcal{J}$ of pairwise disjoint directed $W'''$-paths whose endpoints all lie on the centres of distinct tiles of $\Tiling''$ and that all start at the centres of tiles from $\Class_j'$.
	So we may apply \cref{lemma:directedlongjumpsNew} to find a $K_{t,t}$-\hyperref[def:matchingminor]{matching minor} grasped by \hyperref[def:split]{$\Split{W'''}$}.
	Moreover, with $W'''\subseteq W\subseteq W_0$, we have found a $K_{t,t}$-matching minor grasped by $\Split{W_0}$.
	
	Therefore we may assume that we find a set $Z$ of size at most $8\Fkt{f_P}{t}$ that hits all directed $X_{\text{II}}^i$-paths.
	Let $Z^2_{i,\xi,\xi'}\coloneqq Z\cap \V{D}$, and $\MarkedTiles{2}_{i,\xi,\xi'}\coloneqq\CondSet{T\in\Tiling}{x_T\in Z}$.
	Since $\Abs{Z}\leq 2\Fkt{f_P}{t}$, the bounds on the two sets follow immediately.
	Moreover, since $Z$ meets every directed $X_{\text{II}}^i$-path in $D'$, every directed path with endpoints in distinct tiles of $\Class_i$ which is otherwise disjoint from $W_0$ must contain a vertex from $Z^2_{i,\xi,\xi'}$ or meet a tile from $\MarkedTiles{2}_{i,\xi,\xi'}$.
\end{proof}

\section{Step 2: Crosses and Bipartite $K_{3,3}$-Free Graphs}\label{sec:crosses}

In the previous section, we have taken care of long jumps over our wall.
The proofs presented there are closely related to those necessary for the Directed Flat Wall Theorem.
This section is dedicated to the second step: removing local crosses from our wall.

\paragraph{Non-Pfaffian Cells and Tiles}

The two techniques of changing the perfect matching and applying \hyperref[def:braces]{tight cut contractions}, are of immediate importance for the first step towards controlling the crossings over our wall.

Let $w\in\N$ be some positive integer.
In the following we will say that a slice $W'$ of some cylindrical wall $W$ in some digraph $D$ is \emph{clean} if there is no $w$-long jump over $W'$ in $D$.
We say that the slice $W'$ is \emph{proper} if it does not contain the perimeter of $W$.

A first step is to localise crossings to be able to use them for routing.
Ideally, we want our crossings to occur `over', or, more precisely, within the `attachment' of a single cell in our wall.
However, this is not necessarily possible, as there might still exist short jumps even over a clean slice.
So instead, we consider crossings over tiles and then use a tier II tiling to force these crossings into a single cell.
For this, we need the following:
\begin{itemize}
	\item We need a proper definition of an `attachment',
	\item there needs to be a way to use crossings, or in other words conformal bisubdivisions of $K_{3,3}$ as seen in \cref{thm:4cycleK33}, even if these are not $M$-conformal, and
	\item we need to show that the existence of a short jump with both endpoints in a tile immediately forces the existence of a conformal $K_{3,3}$ subdivision within the attachment of the tile.
\end{itemize}

It is important to note that, since we are interested in conformal bisubdivisions of $K_{3,3}$, it suffices to work on a brace.
We require some additional observations on matching walls in bipartite graphs with perfect matchings and their braces.

\begin{lemma}\label{lemma:wallremnants}
	Let $k\in\N$ be positive integer with $k\geq 2$, $W$ be a \hyperref[def:matchingwall]{matching $k+2$-wall} with canonical matching $M$, $G_k$ be a subdivision the cylindrical grid of order $k$, and $W'$ be the unique proper \hyperref[def:slice]{slice} of width $k$ of $W$.
	Then the following statements are true:
	\begin{enumerate}
		\item $W$ is an expansion of a brace $J$,
		\item $J$ contains a remnant $H$ of $W'$, and
		\item $H$ is a cubic graph that is an expansion of $\Split{G_k}$.
	\end{enumerate} 
\end{lemma}

\begin{proof}
	We start by proving that $W$ has a unique brace that is not isomorphic to $C_4$.
	To do this, let $X\subseteq\V{W}$ be any set that induces a non-trivial tight cut in $W$.
	Moreover, let us assume that $\Complement{X}$ contains more degree-$3$-vertices of $W$ than $X$, while, in case both shores contain the same number of such vertices, let $X$ be chosen arbitrarily.
	Note that $X$ cannot contain vertices from two distinct vertical cycles of $W$, since we may switch their matchings independently and thus could always force at least two matching edges to lie in $\CutG{W}{X}$.
	Similarly, for every pair $P_1$, $P_2$ of horizontal paths, we can always find an $M$-conformal cycle $C$ in $W$ that contains both $P_1$ and $P_2$.
	Hence $X$ cannot contain vertices from both $P_1$ and $P_2$, as otherwise it would either contain vertices from two distinct vertical cycles, or we could switch the matching along $C$ to force at least two edges of the new matching into $\CutG{W}{X}$.
	Hence $\InducedSubgraph{W}{X}$ must be one of three things:
	\begin{enumerate}
		\item an induced subpath of some vertical cycle of $W$,
		\item an induced subpath of some horizontal path of $W$, or
		\item a subdivided star that contains exactly one degree-$3$-vertex $s$ of $W$ and $s$ lies at the centre of the star.
	\end{enumerate}
	In the first two cases, contracting $\Complement{X}$ clearly results in a cycle of even length and thus all of its braces are isomorphic to $C_4$.
	In the last case observe that, in order for $X$ to induce a \hyperref[def:braces]{tight cut}, the leaf vertices of the star must all have the same colour by \cref{obs:tightcutminoritymajority} and thus, since we are forced to have an imbalance of exactly one, the centre of the star must be part of the \hyperref[obs:tightcutminoritymajority]{majority} as well.
	
	\begin{figure}[!ht]
		\centering
		\begin{tikzpicture}[scale=0.9]
			\pgfdeclarelayer{background}
			\pgfdeclarelayer{foreground}
			\pgfsetlayers{background,main,foreground}
			
			\node (a1) [v:main] {};
			
			\node (b1) [v:mainempty,position=45:13mm from a1] {};
			\node (b2) [v:mainempty,position=0:13mm from a1] {};
			\node (b3) [v:mainempty,position=315:13mm from a1] {};
			
			\node (b4) [v:mainempty,position=0:33mm from a1] {};
			
			\node (a2) [v:main,position=135:13mm from b4] {};
			\node (a3) [v:main,position=180:13mm from b4] {};
			\node (a4) [v:main,position=225:13mm from b4] {};
			
			\begin{pgfonlayer}{background}
				
				\draw[e:coloredborder] (a1) to (b1);
				\draw[e:coloredborder] (a2) to (b4);
				\draw[e:coloredborder] (a3) to (b2);
				\draw[e:coloredborder] (a4) to (b3);
				
				\draw[e:colored,color=BostonUniversityRed] (a1) to (b1);
				\draw[e:colored,color=BostonUniversityRed] (a2) to (b4);
				\draw[e:colored,color=BostonUniversityRed] (a3) to (b2);
				\draw[e:colored,color=BostonUniversityRed] (a4) to (b3);
				
				\draw[e:main] (a1) to (b2);
				\draw[e:main] (a1) to (b3);
				\draw[e:main] (a3) to (b4);
				\draw[e:main] (a4) to (b4);
				\draw[e:main] (a2) to (b1);
				
			\end{pgfonlayer}
		\end{tikzpicture}
		\caption{The tight cut contraction of a star.}
		\label{fig:contractedstar}
	\end{figure}
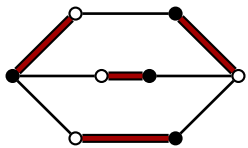
	So when contracting $\Complement{X}$, this case yields a bisubdivision of the graph in \cref{fig:contractedstar}.
	It is straightforward to see that every brace of this graph must be isomorphic to $C_4$.
	As we have seen, one of the two tight cut contractions of any chosen non-trivial tight cut in $W$ yields exclusively braces which are isomorphic to $C_4$.
	We have also seen that no such shore can contain more then one degree-$3$-vertex of $W$, and thus the other shore must still contain a remnant of some brace $J$ which contains all other degree-$3$-vertices.
	Hence $W$ contains a brace $J$ which has at least as many degree-$3$-vertices as $W$ has, and this brace $J$ must be unique.
	Our claim now follows from the uniqueness of the tight cut decomposition \cite{lovasz1987matching}.
	
	Next, observe that there are exactly two vertical cycles $C_1$ and $C_2$ in $W$ which contain degree-$3$-vertices that are linked by paths of even length along $C_i$, whose internal vertices are all of degree two in $W$.
	These two cycles are in fact exactly those which form the perimeter of $W$.
	Hence the second and third point from the assertion follow immediately.
\end{proof}

\begin{definition}[Attachment]\label{def:attachment}
	Let $B$ be a bipartite graph with a perfect matching $M$, $D\coloneqq \DirM{B}{M}$, and $W$ be a proper slice of some cylindrical wall $W'$ in $D$.
	Let $B'$ be a host of $\Split{W'}$ in $B$, and let $\widetilde{W}$ be the remnant of $\Split{W}$ in $B'$.
	
	Now let $H\subseteq B$ and $M'\in\Perf{B}$ be chosen such that
	\begin{enumerate}
		\item $\Split{W}$ is $M'$-conformal in $B$,
		\item $\E{B-\Split{W'}}\cap M'\subseteq M$, and
		\item $H$ is an induced $M'$-conformal subgraph of $\Split{W'}$ such that its outer face, i.\@e.\@ the face of $H$ that contains $\Split{W'}-H$ in the canonical embedding of $\Split{W}$, is an $M'$-conformal cycle.
	\end{enumerate}
	Let $H'$ be the remnant of $H$ in $B'$.
	The \emph{attachment of $H$ over $W$ in $B$}, denoted by $\Attachment{B,\Split{W}}{H}$, is the elementary component of $B'-\Brace{\widetilde{W}-H'}$ that contains the outer face of $H'$.
\end{definition}

\begin{lemma}[\cite{giannopoulou2021two}]\label{lemma:planarnonbraceparts}
	Let $B$ be a planar brace and $C$ a conformal and separating cycle in $B$, then $\Inner{B}{C}$ and $\Outer{B}{C}$ are matching covered.
\end{lemma}

Let $W$ be a \hyperref[def:matchingwall]{matching wall} and $H$ be an induced $M'$-conformal subgraph of $W$ whose outer face is an $M$-conformal cycle $C$ where $M'\in\Perf{W}$.
Then $C$ is conformal and separating in $W$ and thus, by \cref{lemma:planarnonbraceparts}, $H$ is matching covered.
That means, for every $H\subseteq W$ for which an attachment exists, we know that $H'$, as in the definition above, belongs to this attachment.

Since we are interested in tilings and their centres especially, we need to show that the attachment of a cell and the attachment of a tile are well defined objects.
Moreover, we would like to know that the attachment of a tile contains the attachment of its centre.
To this end, we introduce the following two lemmas.

\begin{lemma}\label{lemma:conformaltiles}
	Let $h,k,w\in\N$ be positive integers and $W$ be a proper \hyperref[def:slice]{slice} of width $k$ of some \hyperref[def:matchingwall]{matching $h$-wall} where $h\geq k+2$.
	Let $M$ be the canonical matching of $W$ and \hyperref[def:tiling]{$\Tiling=\Tiling_{W,k,w,\xi,\xi'}$} for $\xi,\xi'\in[1,w+1]$.
	At last let $T\in\Tiling$ be any \hyperref[def:tile]{tile} that is completely contained in $W$ and let $Q$ be the vertical cycle of $W$ that contains the leftmost part of $\Perimeter{T}$.
	Then $M_T\coloneqq M\Delta \E{Q}$ is a perfect matching of $W$ such that $\Perimeter{T}$ contains an $M_T$-conformal cycle $C_T$ and $T-\Brace{\Perimeter{T}-C_T}$ is an $M_T$-conformal matching covered subgraph of $W$ whose outer face is bounded by $C_T$.
\end{lemma}

\begin{figure}[!ht]
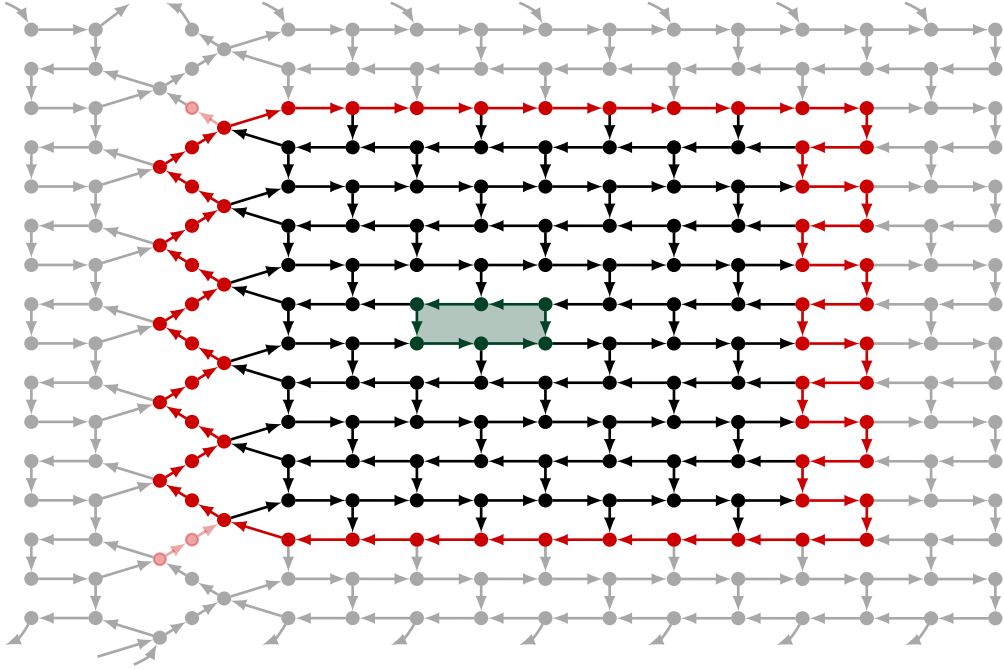

	\centering

	\caption{The tile $T$ in the $M_T$-direction of $W$.}
	\label{fig:aconformaltile}
\end{figure}

\begin{proof}
	To see that the statement is true, observe that the boundary of $T$ consists of a cycle and two vertices, one of them being the starting vertex $s$ of the subpath of the upper path in $\Perimeter{T}$, and the other one being the endpoint $t$ of the subpath of the lower path in $\Perimeter{T}$.
	By definition, $Q$ is the left path of $\Perimeter{T}$ and $L\coloneqq Q\cap\Perimeter{T}$.
	Then both $s$ and $t$ belong to $L$.
	See \cref{fig:atile} for an illustration.
	Now consider $L'\coloneqq L-s-t$ and note that $\Split{L'}$ contains exactly two vertices which do not belong to the unique cycle of $\Split{\Perimeter{T}}$.
	Let $P$ be the path obtained from $\Split{L'}$ by removing these two vertices.
	Since we removed both endpoints of an $M$-conformal path, where $M$ is the canonical matching of $W$, the result is an internally $M$-conformal path.
	Let $C$ be the unique cycle of $\Split{\Perimeter{T}}$, then, as $\Perimeter{T}-Q$ is a directed path, $C-P$ is another internally $M$-conformal path.
	Now consider the perfect matching $M_T$.
	Then every internally $M$-conformal subpath of $Q$ has now become an $M_T$-conformal path.
	So now $C-P$ is internally $M_T$-conformal, while $P$ is $M_T$-conformal.
	Hence $C$ is $M_T$-conformal as required.
	See \cref{fig:aconformaltile} for an illustration.
\end{proof}

\begin{figure}[!ht]
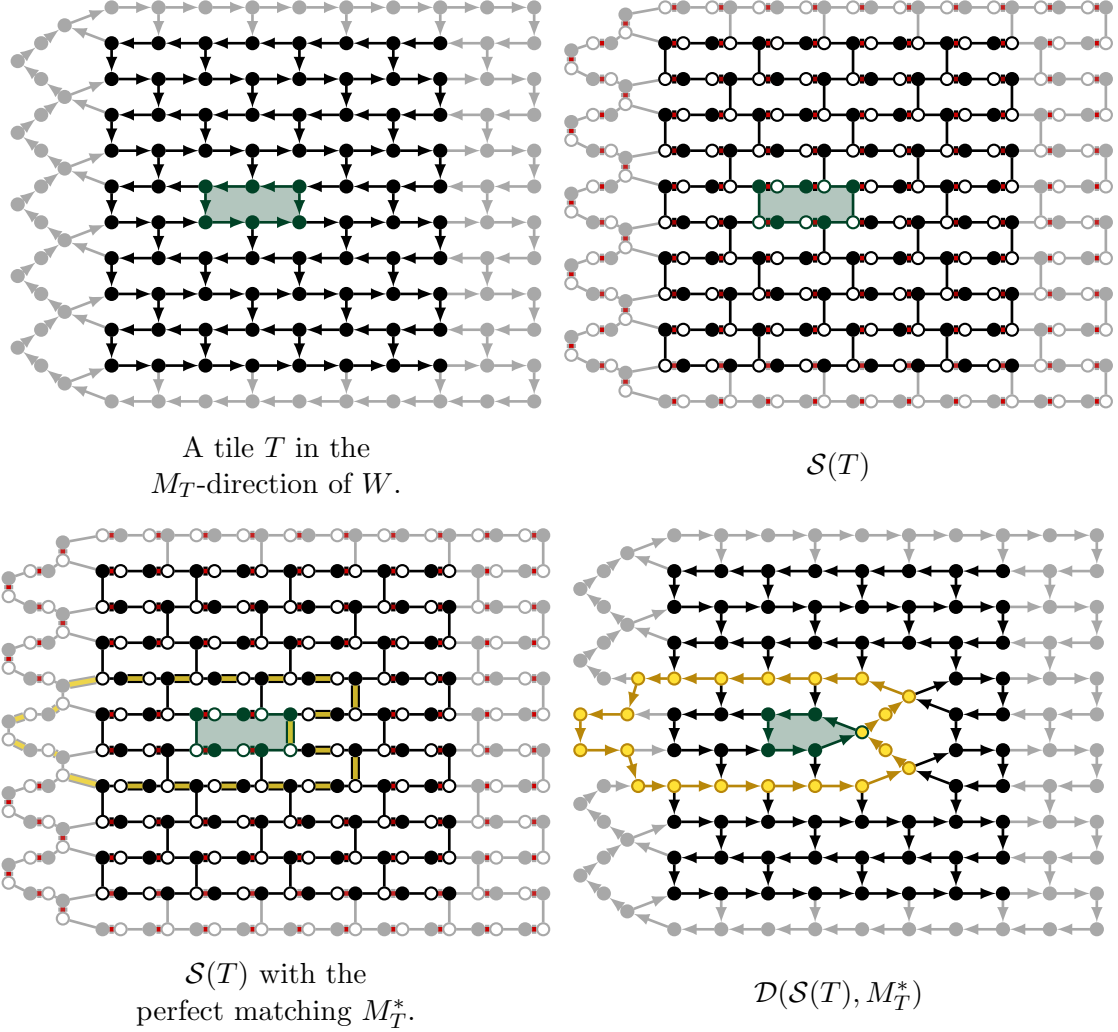

	\centering

		};
		
		\node (LabelLU) [v:ghost,position=270:34mm from LU,align=center] {A tile $T$ in the\\ $M_T$-direction of $W$.};
		\node (LabelRU) [v:ghost,position=270:34mm from RU,align=center] {$\Split{T}$};
		\node (LabelLD) [v:ghost,position=270:34mm from LD,align=center] {$\Split{T}$ with the\\ perfect matching $M_T^*$.};
		\node (LabelRD) [v:ghost,position=270:34mm from RD,align=center] {$\DirM{\Split{T}}{M_T^*}$};
		
	\end{tikzpicture}
	\caption{The tile $T$ in the $M_T^*$-direction of $W$ and its split.}
	\label{fig:aconformaltilewithconformalcentre}
\end{figure}

\begin{lemma}[\cite{mccuaig2004polya}]\label{lemma:conformalfaces}
	Let $B$ be a bipartite and planar matching covered graph, then every facial cycle of $B$ is conformal.
\end{lemma}

\begin{observation}\label{obs:conformalcentres}
	Let $h,k,w\in\N$ be positive integers and $W$ be a proper slice of width $k$ of some \hyperref[def:matchingwall]{matching $h$-wall} where $h\geq k+2$.
	Let $M$ be the canonical matching of $W$, and let $C$ be a cell of $W$.
	Then $C$ is a conformal cycle of $W$.
\end{observation}

\begin{proof}
	As $W$ is a slice of a matching wall, it is planar and matching covered.
	Since $C$ bounds a face of $W$, our claim follows immediately from \cref{lemma:conformalfaces}. 
\end{proof}

Hence the \hyperref[def:attachment]{attachments} of tiles and their \hyperref[def:tile]{centres} are indeed well defined.
Next, we want to see that the attachment of the centre of a tile is contained in the attachment of the tile itself.

\begin{lemma}\label{lemma:centreattachments}
	Let $h,k,w\in\N$ be positive integers and $W$ be a proper and clean slice of width $k$ of some \hyperref[def:matchingwall]matching $h$-wall in a bipartite matching covered graph $B$, where $h\geq k+2$.
	Let $M$ be a perfect matching of $W$ that contains the canonical matching of $W$ and \hyperref[def:tiling]{$\Tiling=\Tiling_{W,k,w,\xi,\xi'}$} for $\xi,\xi'\in[1,w+1]$.
	At last let $T\in\Tiling$ be any tile that is completely contained in $W$, and let $C_T$ be the cycle that bounds the centre of $T$.
	Then $\Attachment{B,W}{C_T}\subseteq\Attachment{B,W}{T}$.
\end{lemma}

\begin{proof}
	By \cref{lemma:conformaltiles} we already know that $T$ is an $M_T$-conformal and matching covered subgraph of $W$.
	As $C_T$ bounds a face of $T$, which is matching covered, \cref{lemma:conformalfaces} still guarantees us the existence of a perfect matching $M_T^*$ of $B$ such that $T$ and $C_T$ are $M_T^*$-conformal and the $M_T^*$-direction of $T$ is strongly connected.
	Let $A_T^*$ be the $M_T^*$-direction of the \hyperref[def:attachment]{attachment} of $T$, while $A_{C_T}^*$ is the $M_T^*$-direction of the attachment of $C_T$.
	Then for every vertex $v$ from $A_{C_T}^*$ there exists a directed path from $v$ to $\DirM{C_T}{M_T^*}$ and a directed path from $\DirM{C_T}{M_T^*}$ to $v$, since $A_{C_T}^*$ must be strongly connected by the definition of attachments.
	See \cref{fig:aconformaltilewithconformalcentre} for an example of the construction of $M_T^*$.
	This, however, means that $\Split{\DirM{T}{M_T^*}\cup A_{C_T}^*}$ must be matching covered.
	Let $B'$ be the host of $W$ in $B$, and let $\widetilde{W}$, $\widetilde{T}$ be the remnants of $W$ and $T$ in $B'$.
	Then $\Attachment{B,W}{C_T}$ must be contained in the elementary component of $B'-\Brace{\widetilde{W}-\widetilde{T}}$ that contains $\widetilde{T}$ and thus we are done.
\end{proof}

Suppose $C_T$ is the centre of some tile $T$ in a slice $W$ as above.
If $\Attachment{B,W}{C_T}$ is non-Pfaffian it must contain a conformal bisubdivision of $K_{3,3}$.
However, to be able to use \cref{thm:4cycleK33} and \cref{lemma:goodcrossesmeanK33}, we need to show that these conformal bisubdivisions of $K_{3,3}$ cannot be separated from the remnant $\widetilde{T}$ of $T$ by a non-trivial \hyperref[def:braces]{tight cut} in $\Attachment{B,W}{T}$. 

\begin{lemma}\label{lemma:centralK33}
	Let $h,k,w\in\N$ be positive integers and $W$ be a proper and clean slice of width $k$ of some matching $h$-wall in a bipartite matching covered graph $B$, where $h\geq k+2$.
	Let $M$ be a perfect matching of $W$ that contains the canonical matching of $W$ and $\Tiling=\Tiling_{W,k,w,\xi,\xi'}$ for $\xi,\xi'\in[1,w+1]$.
	At last, let $T\in\Tiling$ be any tile that is completely contained in $W$, and let $C_T$ be the cycle that bounds the centre of $T$.
	Then every non-trivial tight cut in $\Attachment{B,W}{T}$ has a shore $X$ such that $\V{\Attachment{B,W}{C_T}}\subseteq X$, and $\ContractsTo{\Attachment{B,W}{T}}{\Complement{X}}{v_{\Complement{X}}}$ has a brace that is the host of $T$ in $\Attachment{B,W}{T}$.
\end{lemma}

\begin{proof}
	Let $B'$ be the host of $W$ in $B$ and let $\widetilde{W}$, $\widetilde{T}$, and $\widetilde{C_T}$ be the remnants of $W$, $T$, and $C_T$ in $B'$ respectively.
	Now let $A\coloneqq\Attachment{B,W}{T}$ and consider a non-trivial tight cut $\CutG{A}{X}$ in $A$.
	Note that neither $X$ nor $\Complement{X}$ can contain vertices purely from $A-\widetilde{T}$, since otherwise $\CutG{B'}{X}$ would be a non-trivial tight cut in $B'$.
	Indeed, $X$ or $\Complement{X}$ must contain vertices of the remnant of $\Perimeter{T}$ in $B'$.
	Suppose $\Complement{X}$ is that shore.
	Let us call a cell of $T$ an \emph{inner cell} if it does not contain vertices of the perimeter of $T$.
	If $\Complement{X}$ contains vertices of the remnant of an inner cell of $T$, then there exist two disjoint conformal cycles of $T$, both of which have remnants in $B'$ with edges in $\CutG{B'}{X}$, and thus $\CutG{B}{X}$ cannot be a tight cut.
	At last, suppose $X$ contains vertices of $\Attachment{B,W}{C_T}$, then we may find an  $M_T$-conformal cycle $C'$ in $\Attachment{B,W}{C_T}$ that contains vertices of $\widetilde{C_T}$ and vertices of $\Complement{X}$, but which is disjoint from the remnant of $\Perimeter{T}$.
	This again yields a contradiction to $\CutG{B'}{X}$ being tight and thus our claim follows.
\end{proof}

The last remaining piece before we can attempt to create the desired crossings is:
What if the attachment of every cell is Pfaffian, but there still is a short jump with both endpoints in the interior of $T$?

\begin{lemma}\label{lemma:shortjumpsmeanK33}
	Let $h,k,w\in\N$ be positive integers and $W$ be a proper and clean slice of width $k$ of some \hyperref[def:matchingwall]{matching $h$-wall} in a bipartite matching covered graph $B$, where $h\geq k+2$ and $w\geq 2$.
	Let $M$ be a perfect matching of $W$ that contains the canonical matching of $W$ and \hyperref[def:tiling]{$\Tiling=\Tiling_{W,k,w,\xi,\xi'}$} for $\xi,\xi'\in[1,w+1]$.
	At last, let $T\in\Tiling$ be any tile that is completely contained in $W$ and $J$ be a short jump in $\DirM{B}{M}$ over $\DirM{W}{M}$ with both endpoints in the interior of $\DirM{T}{M}$.
	Then, if $B'$ is the host of $W$ in $T$ and $\widetilde{T}$ is the remnant of $T$ in $B'$, $T'$ contains a conformal bisubdivision $H$ of $K_{3,3}$ such that all six degree-$3$-vertices of $H$ belong to $\widetilde{T}$.
\end{lemma}

\begin{figure}[!ht]
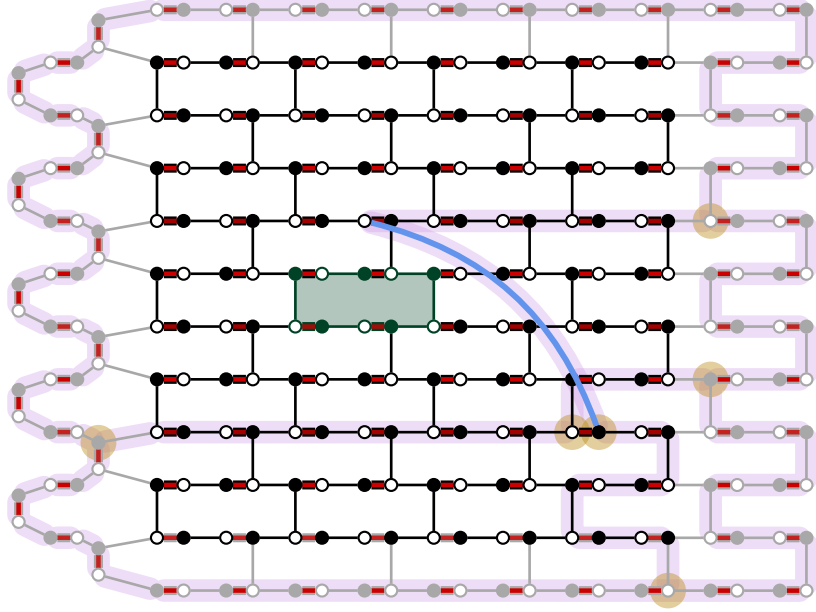

	\centering

	\caption{An $M_T$-conformal bisubdivision of $K_{3,3}$ constructed using a short jump.}
	\label{fig:conformalK33viashortjump}
\end{figure}

\begin{proof}
	Consider the perfect matching $\widetilde{M}$ of $\widetilde{T}$ that is maintained from $M_T$ through he \hyperref[def:braces]{tight cut contractions} applied.
	First of all note that, since $J$ is a short jump in $\DirM{B}{M}$ over $\DirM{W}{M}$, and by \cref{lemma:wallremnants} all faces of $W$ are preserved in $\widetilde{W}$, which is the remnant of $W$ in $B'$, it corresponds to an internally $\widetilde{M}$-conformal path $\widetilde{J}$ with endpoints $a_1\in V_1\cap\V{\widetilde{T}}$ and $b_1'\in V_2\cap\V{\widetilde{T}}$ such that $a_1$ and $b_1'$ do not belong to the same face of $\widetilde{T}$, and $\widetilde{J}$ is internally disjoint from $\widetilde{T}$.
	Let $a_1'$ and $b_1$ be chosen such that $a_1b_1,a_1'b_1'\in\widetilde{M}$.
	Since $\widetilde{T}$ is matching covered there exist $a_2\in V_2\cap \V{\Perimeter{\widetilde{T}}}$ and $b_3\in V_1\cap\V{\Perimeter{\widetilde{T}}}$ such that
	\begin{enumerate}
		\item $a_2b_2,a_3b_3\in\widetilde{M}$ are distinct,
		\item the internally $\widetilde{M}$-conformal subpath $P_1$ of $\Perimeter{\widetilde{T}}$ with endpoints $a_2$ and $b_3$ is non-trivial, and
		\item $P_2\coloneqq \Perimeter{\widetilde{T}}-P_1$ is as short as possible such that there exist disjoint internally $\widetilde{M}$-conformal paths $L_a$, $L_b$ where $L_a$ connects $b_1$ to $a_2$ and $L_b$ connects $a_1'$ to $b_3$.
	\end{enumerate}
	Now, since $w\geq 2$, there exist $b_4$ and $a_5$ on $P_1$ such that the path $R$ connecting $b_4$ and $a_5$ in $P_1$ is $\widetilde{M}$-conformal, and there are internally $\widetilde{M}$-conformal and disjoint paths $Q_a$ and $Q_b$ in $\widetilde{T}$ that satisfy:
	\begin{enumerate}
		\item $Q_a$ connects $a_1$ and $b_4$,
		\item $Q_b$ connects $b_1$ and $a_5$, and
		\item $Q_a$ and $Q_b$ avoid $L_a$ and $L_b$.
	\end{enumerate}
	By the choices of these paths, we have found in total nine pairwise internally disjoint paths such that 
	\begin{itemize}
		\item $a_1$ and $b_1$ are joined by a $\widetilde{M}$-conformal path; that is the edge $a_1b_1$,
		\item $a_2$ and $b_3$ are joined by a $\widetilde{M}$-conformal path; namely the path $P_2$ together with the edges $a_2b_2$ and $b_3a_3$,
		\item $a_5$ and $b_4$ are joined by the $\widetilde{M}$-conformal path $R$,
		\item $a_1$ is joined by the internally $\widetilde{M}$-conformal path $\widetilde{J}b_1'a_1'L_b$ to $b_3$ and by the internally $\widetilde{M}$-conformal path $Q_a$ to $b_4$,
		\item $b_1$ is joined by the internally $\widetilde{M}$-conformal path $L_a$ to $a_2$, and by $Q_b$ to $a_5$, and
		\item $a_5$ is joined to $b_3$ by an internally $\widetilde{M}$-conformal subpath of $P_1$, while $b_4$ and $a_2$ are joined by the remaining internally $\widetilde{M}$-conformal subpath of $P_1$ that does not contain $a_5$.
	\end{itemize}
	Hence overall, we obtain an $\widetilde{M}$-conformal bisubdivision of $K_{3,3}$.
	See \cref{fig:conformalK33viashortjump} for an illustration.
	
	Another, less constructive, way to see that there must exist a conformal $K_{3,3}$ bisubdivision is to observe that $\widetilde{T}$ together with $\widetilde{J}$ is matching covered and non-planar.
	Moreover, similar to the proofs of \cref{lemma:wallremnants} and \cref{lemma:centralK33}, any non-trivial tight cut must either sit within a bisubdivided edge of $T$, form a star around a vertex of degree at least three, or consist entirely of internal vertices of $\widetilde{J}$.
	Hence $\widetilde{T}\cup\widetilde{J}$ is an expansion of a non-planar brace $B''$ that closely resembles $\widetilde{T}$ with a single short jump over it.
	Thus $B''$ is not the Heawood graph and therefore it must contain a conformal bisubdivision of $K_{3,3}$ by \cref{thm:trisums}.
\end{proof}

\paragraph{A Tier II Cross}

By now we know that non-Pfaffian \hyperref[def:attachment]{attachments} of cells and short jumps force conformal $K_{3,3}$-bisubdivisions to exist.
Next we show how to use this knowledge to create a perfect matching $M'$ such that the $M'$-direction of our digraph contains two crossing paths within the centre of a tile.
More precisely, we show how to create a cross over a single tile and then place a tier II tile around it to obtain a tile with a cross in its centre.
Hence it suffices to show that we can locally manipulate the perfect matching of the attachment of a cell in order to create a cross over the cell, if we know that the attachment of the centre contains a structure similar to a tile itself.

Given the \hyperref[def:tile]{centre} $C$ of a tile $T$ in a \hyperref[def:matchingwall]{matching wall} $W$ with canonical matching $M$, we are particularly interested in the following four vertices:
There are exactly two vertical cycles of $W$ that meet $C$, let us call them $Q_j$ and $Q_{j+1}$, where $Q_j$ lies left of $Q_{j+1}$ in the canonical embedding.
Let $U_i\coloneqq Q_i\cap C$ for each $i\in[j,j+1]$.
Then $U_j$ is $M$-conformal, while $U_{j+1}$ is internally $M$-conformal.
Each of the $U_i$ has exactly one endpoint in the upper path of $C$ and one endpoint in the lower path of $C$.
Moreover, let
\begin{itemize}
	\item $b^1_{T}\in V_2$ be the endpoint of $U_j$ on the upper path of $C$,
	\item $a^2_{T}\in V_1$ be the endpoint of $U_j$ on the lower path of $C$,
	\item $a^1_{T}\in V_1$ be the endpoint of $U_{j+1}$ on the upper path of $C$, and
	\item $b^2_{T}\in V_2$ be the endpoint of $U_{j+1}$ on the lower path of $C$.
\end{itemize}
Additionally there exist $c^1_{T}$ and $c^2_{T}$ on $U_{j+1}$ such that $a^1_{T}c^2_{T}\in M$, and $b^2_Tc^1_{T}\in M$.
From the definition it follows that $c^1_{T},c^2_{T}\notin\V{C}$.

We say that a tile $T$ is \emph{non-Pfaffian} if either the \hyperref[def:attachment]{attachment} of the \hyperref[def:tile]{centre} of $T$ is a non-Pfaffian bipartite graph, or there exists a short jump over $W$ with both endpoints on the interior of $T$.

\begin{lemma}\label{lemma:crossdirection}
	Let $h,k,w\in\N$ be positive integers and $W$ be a proper and clean slice of width $k$ of some matching $h$-wall in a bipartite matching covered graph $B$, where $h\geq k+2$.
	Let $M$ be a perfect matching of $W$ that contains the canonical matching of $W$ and $\Tiling=\Tiling_{W,k,w,\xi,\xi'}$ for $\xi,\xi'\in[1,w+1]$ with a four colouring $\Set{\Class_1,\dots,\Class_4}$, and for some $i\in[1,4]$ let $W'\coloneqq\InducedSubgraph{W}{\Tiling,i}$.
	Consider \hyperref[def:tilingII]{$\TierIITiling{\Tiling,i,w}{W}$} and let us select $T\in\Class_i$ as well as $T'\in \TierIITiling{\Tiling,i,w}{W'}$ such that the centre of $T'$ contains $T$ in its interior in the canonical embedding of $W$.
	Let $C$ be the cycle of $B$ that is the centre of $T'$ and let $I_{T'}$ be the union of $C$ and the component of $W-C$ that contains $T$.
	For $i\in[1,2]$ let $a^i_{T'}$ and $b^i_{T'}$ be defined as above.
	
	If $T$ is non-Pfaffian, there exist a perfect matching $N_T$ of $B$ such that $\CondSet{e\in M}{e\cap \V{W-I_{T'}}\neq\emptyset}\subseteq N_T$, and there exist vertex disjoint paths $R_1$ and $R_2$ in $B$ such that
	\begin{enumerate}
		\item $R_1$ and $R_2$ are internally vertex disjoint from $W'$,
		\item $R_1$ and $R_2$ are fully contained in $\Attachment{B,W'}{T'}$
		\item both paths are $N_T$-alternating,
		\item $R_1$ has endpoints $a^1_{T'}$ and $a^2_{T'}$ and the edge of $R_1$ that is incident with $a^2_{T'}$ lies in $N_T$, and
		\item $R_2$ has endpoints $b^1_{T'}$ and $b^2_{T'}$ and the edge of $R_2$ that is incident with $b^2_{T'}$ lies in $N_T$.
	\end{enumerate}
\end{lemma}

\begin{proof}
	By \cref{lemma:conformaltiles} we know that $I_{T'}$ is matching covered.
	Moreover, since $T$ is a conformal subgraph of $I_{T'}$, we must also have $\Attachment{B,W}{T}\subseteq\Attachment{B,W'}{T'}$.
	Therefore $\Attachment{B,W'}{T'}$ is non-Pfaffian as well.
	
	Let us add the edges $a^i_{T'}b^j_{T'}$ to $\Attachment{B,W'}{T'}$ for all $i,j\in[1,2]$ and let $G$ be the resulting bipartite graph.
	By \cref{thm:bipartiteextendibility}, $G$ is still matching covered.
	
	We claim that the host $B'$ of $I_{T'}$ in $G$ is also non-Pfaffian, and that it contains four distinct vertices, each of them representing a vertex from $\Set{a^1_{T'},a^2_{T'},b^1_{T'},b^2_{T'}}$ in the following sense:
	We say that a vertex $u$ of $B'$ \emph{represents} a vertex $v$ of $I_{T'}$ if $u,v\in V_i$ for some $i\in[1,2]$, and there exists a tight cut $\CutG{B}{X}$ in $G$ such that $v\in X$ and $u$ is the vertex of $B'$ obtained from contracting $X$.
	
	Note that in this case, every edge $e$ incident with $u$ in $B'$ can be replaced by a path $P_v$, which has exactly one endpoint outside of $X$, in $G$ that shares the same properties regarding any perfect matching $M'$ of $G$ as $e$ does regarding the remainder of $M'$ in $B'$.
	That is, if $e\in\Remainder{M'}{B'}$, then $P_v$ is $M'$-conformal, and otherwise $P_v$ is internally $M'$-conformal.
	
	Towards the validity of our claim first let $C'$ be the centre of $T$ and suppose $\Attachment{B,W}{C'}$ is non-Pfaffian.
	In this case observe that $I_{T'}$ indeed meets all requirements of a tile in $W$ and thus we may call upon \cref{lemma:centralK33} to see that $B'$ must be non-Pfaffian.
	If the attachment of the centre of $T$ is Pfaffian, there must exist a short jump $J$ over $W$ with both endpoints in $T$.
	Then, by \cref{lemma:shortjumpsmeanK33} and its proof, one can see that $I_{T'}\cup J$ contains a conformal bisubdivision of $K_{3,3}$ that contains the cycle $C$.
	Hence again $B'$ is non-Pfaffian.
	
	One can observe that every non-trivial tight cut in $I_{T'}$ must either be a cut around a bisubdivided claw, or along a bisubdivided edge of $W$.
	Moreover, in $\Attachment{B,W'}{T'}$ the only non-trivial tight cuts can occur on the outer face of $T'$ since any other non-trivial tight cut would correspond to a non-trivial tight cut in the brace that was used to construct $\Attachment{B,W'}{T'}$.
	Note that every pair among the four vertices $\Set{a^1_{T'},a^2_{T'},b^1_{T'},b^2_{T'}}$ can be separated on $C$ by two degree-$3$-vertices of $I_{T'}$ which are not in $\Set{a^1_{T'},a^2_{T'},b^1_{T'},b^2_{T'}}$.
	Additionally, for each $i\in[1,2]$, the two vertices in $\Set{a^1_{T'},a^2_{T'},b^1_{T'},b^2_{T'}}\cap V_i$ can be separated by $\Set{a^1_{T'},a^2_{T'},b^1_{T'},b^2_{T'}}\cap V_{3-i}$ on $C$.
	Hence we may assume for each $x\in \Set{a^1_{T'},a^2_{T'},b^1_{T'},b^2_{T'}}$ to have a vertex $u_x\in\V{B'}$ that represents $x$ and that all of these vertices are pairwise distinct.
	
	By our addition of the fresh edges to obtain $G$, we now have the edges $a^i_{T'}b^j_{T'}$ in $B^+$ for all $i,j\in[1,2]$.
	Let $\hat{C}$ be the four-cycle consisting exactly of these four edges, and let $B^+$ be the resulting graph.
	By \cref{thm:bipartiteextendibility}, $B^+$ is still a brace.
	Hence we may use \cref{thm:4cycleK33} to find a conformal bisubdivision $H$ of $K_{3,3}$ in $B^+$ that contains $\hat{C}$ as a subgraph.
	Let $N''$ be the perfect matching of $H$ that contains the edge $u_{a^1_{T'}}u_{b^2_{T'}}$ but not the edge $u_{a^2_{T'}}u_{b^1_{T'}}$.
	As $H$ is a bisubdivision of a brace and $\hat{C}$ is a subgraph of $H$, $N''$ must exist.
	Now let $N'$ be a perfect matching of $B^+$ that contains $N''$.
	As $H$ is a bisubdivision of $K_{3,3}$ there exist paths $R_1'$ and $R_2'$ with the following properties in $H$:
	\begin{enumerate}
		\item $R_1'$ and $R_2'$ are vertex disjoint and $N'$ alternating,
		\item $R_1'$ has endpoints $u_{a^1_{T'}}$ and $u_{a^2_{T'}}$, and the edge of $R_1'$ that is incident with $a^2_{T'}$ lies in $N'$, and
		\item $R_2'$ has endpoints $u_{b^1_{T'}}$ and $b_{a^2_{T'}}$, and the edge of $R_2'$ that is incident with $b^2_{T'}$ lies in $N'$.
	\end{enumerate}
	Let $N$ be a perfect matching of $\Attachment{B,W'}{T'}+\hat{C}$ such that $N'=\Remainder{N}{B''}$.
	The paths $R_1'$ and $R_2'$ can now be extended to the desired paths $R_1$ and $R_2$ in $\Attachment{B,W'}{T'}+\hat{C}$.
	
	Note that the only edge of $N$ that is not an edge of $B$ is $a^1_{T'}b^2_{T'}$.
	Also the only edges of $\E{W'}\cap M$ with exactly one endpoint in $I_{T'}$ are the two edges $a^1_{T'}c^2_{T'}$ and $b^2_{T'}c^1_{T'}$.
	Hence
	\begin{align*}
		N_T\coloneqq \Brace{M\setminus N}\cap\E{W'}\cup\Brace{N\setminus\Set{a^1_{T'}b^2_{T'}}}
	\end{align*}
	is a perfect matching of $B$ with properties as required by the assertion.
	In particular, since the subpath of $C$ that is parallel to the edge $a^1_{T'}b^2_{T'}$ is internally $M$-conformal, every edge of $M\cap {W-I_{T'}}$ belongs to $N_T$.
\end{proof}

\begin{figure}[!ht]
	\centering

	\caption{The $N_T$-direction of the tile $T'$, the cross obtained through \cref{lemma:crossdirection}, and two paths through the tile that use the cross. Note that we cannot guarantee anything specific about the matching within $I_{T'}$ besides to the matching edges that leave $I_{T'}$ and the crossing itself.}
	\label{fig:atilewithacrossinthecentre}
\end{figure}

While we cannot really control how the perfect matching $N_T$ changes the structure of $I_{T'}$ in the $N_T$-direction of $B$ compared to its $M$-direction, we still know that most of our cylindrical wall is intact.
Indeed, by our choices of the four vertices which are connected via the paths $R_1$ and $R_2$ and the fact that $c^2_T\in V_2$ and $c^1_T\in V_1$, we may follow along two parallel vertical cycles from the top of $T'$, then move towards its centre, enter the two $N_T$-directions of $R_1$ and $R_2$, and switch between the two cycles.
See \cref{fig:atilewithacrossinthecentre} for an illustration.
This means that we are able to switch the relative position of two disjoint paths that are routed through our cylindrical wall given that these paths enter $T'$ from the top.
This means that, given enough such tiles $T'$ that are mutually far enough apart within $W$, we are able to re-order any given family of pairwise disjoint paths in any order we see fit.
Of course within a cylindrical wall any directed path is only ever allowed to go left, right or downwards, but we may take some slice free of crossings and switch the perfect matching $M$ along all vertical cycles of this slice.
By doing so we create a perfect matching $M'$, where we can have a large family of pairwise disjoint paths move `upwards within $W$' without changing their relative positions to one another.

With this we have established all tools necessary for our proof of \cref{thm:matchingflatwall}.

\section{Proof of \Cref{thm:matchingflatwall}}\label{sec:proof}

This section is entirely dedicated to the proof of \cref{thm:matchingflatwall}.
The proof goes through several phases and establishes a few subclaims along the way as follows:
\begin{enumerate}
	\item[$\bullet$] In the beginning we have a large \hyperref[def:matchingwall]{matching wall} $W$ with its canonical matching $M$ of $B$.
	From here on we mostly work in the digraphic setting and consider the cylindrical wall $U_0\coloneqq\DirM{W}{M}$ in $D=\DirM{B}{M}$.
	\item[] \textbf{Phase I}
	\item[] Phase I is an iterative process consisting of two steps which are applied to a slice $U_i$ of $U_0$, where $U_i$ is the end result of the $i$th round of Phase I:
	\begin{description}
		\item[\textbf{Step I}] Here we simultaneously apply \cref{lemma:longjumps1} for both parametrisations of $U_0$, every tiling defined by choices of $\xi,\xi'\in[1,w+1]$, and every colour class.
		This either results in a $K_{t,t}$-\hyperref[def:matchingminor]{matching minor} grasped by \hyperref[def:Split]{$\Split{U_0}$}, or in bounded size sets $F'_i$ and $\mathcal{F}'_i$ which are a set of vertices and a set of \hyperref[def:tile]{tiles} respectively.
		The tiles in $\mathcal{F}'_i$ are considered \emph{marked}.
		We then find a large slice $U'_i$ of $U_{i-1}$ without vertices of marked tiles.
		\item[\textbf{Step II}] Next we apply \cref{lemma:longjumpspahse2} for all four colours of every \hyperref[def:tiling]{tiling} defined by $\xi,\xi'\in[1,w+1]$ and again either find a $K_{t,t}$-\hyperref[def:matchingmonir]{matching minor} grasped by \hyperref[def:split]{$\Split{U_0}$}, or in bounded size sets $F''_i$ and $\mathcal{F}''_i$, which are a set of vertices and a set of tiles respectively.
		We then find a slice $U_i$ of $U_i'$ of appropriate size that does not contain vertices of a marked tile.
	\end{description}
	So at the end of Phase I we have found bounded size sets $F_I$ and $\mathcal{F}_I$ of vertices and tiles respectively which are now considered \emph{marked}, and we have found a large slice $U_I$ of $U_0$ which is free of marked vertices and tiles.
	\item[$\bullet$] We can now show that there either exists a $K_{t,t}$-matching minor grasped by $Split{U_0}$, or there is no long jump over $U_I$ in $D-F_I$ that avoids the tiles from $\mathcal{F}_I$.
	Moreover, the number of tiles in $\mathcal{F}_I$ grows with every round and thus, after a certain threshold of such tiles has been surpassed, we can find many pairwise disjoint long jumps over $U_0$ within a single long jump over $U_I$, which implies the existence of a large $K_{t,t}$-\hyperref[def:matchingminor]{matching minor}.
	\item[] \textbf{Phase II}
	\item[] Now we know that only short jumps can exist over $U_I$ in $D-F_I$.
	We divide $U_I$ into $\Fkt{\mathcal{O}}{t^4}$ slices, each of which can be partitioned into three appropriately sized slices.
	Then either one of these middle slices is already flat with respect to its perimeter, or we may apply \cref{lemma:crossdirection} to locally change the perfect matching of a single tile for each of these middle slices and therefore we find $\Fkt{\mathcal{O}}{t^2}$ disjoint crosses that are relatively far apart from each other.
	These can then be used to create $\Bidirected{K_t}$.
\end{enumerate}
So in Phase II we either find the flat wall as desired, or $\Bidirected{K_t}$ in the $M'$-direction of our graph, where $M'$ is the perfect matching obtained by creating the crossings through \cref{lemma:crossdirection} and switching $M$ along some of the vertical cycles of $U_I$.

Throughout the proof we will collect tiles from different tilings.
In general we will say that a tile is \emph{marked} if it either contains a vertex of some separator obtained from \cref{thm:directedlocalmenger} and \cref{lemma:longjumps1}, or through \cref{lemma:longjumpspahse2}, or one of these separators contains the fresh vertex created for it in the construction of some auxiliary graph (type I or II).
This means that the overall number of tiles that are marked exceeds the size of the sets $\mathcal{F}_i$, still the number of columns that can contain marked tiles is bounded, which is enough for our purposes.
A vertex of $U$ is said to be \emph{marked} if it belongs to a separator obtained from \cref{thm:directedlocalmenger} and \cref{lemma:longjumps1}, or through \cref{lemma:longjumpspahse2}.
In each of the steps we introduce families of marked tiles and vertices and a slice of $U$ is said to be \emph{clear} if it does not contain vertices that are marked.

\begin{proof}[Proof of \Cref{thm:matchingflatwall}]
	Let $r,t\in\N$ be positive integers, $B$ be a bipartite graph with a perfect matching $M$, $D\coloneqq\DirM{B}{M}$, and $W$ be an $M$-conformal matching $\Fkt{\WallBound}{t,r}$-wall, where $\Fkt{\WallBound}{t,r}$ will be determined throughout the proof.
	To do this we will introduce a constant $d_i$ for each step, for which we will make more and more assumptions in the form of lower bounds.
	Let $U\coloneqq\DirM{W}{M}$ be the $M$-direction of $W$, then $U$ is a cylindrical $\Fkt{\WallBound}{t,r}$-wall in $D$.
	\begin{align*}
		\text{Let us assume $\Fkt{\WallBound}{t,r}\geq 3d_1$.}
	\end{align*}
	Then let $U$ be the cylindrical $3d_1$-wall $U=\Brace{Q_1,\dots,Q_{3k},\hat{P}_1,\dots,\hat{P}_{3k}}$ with its triadic partition $\mathcal{U}=\Brace{U,d_1,\widetilde{U}_1,\widetilde{U}_2,\widetilde{U}_3,\widetilde{U}^1,\widetilde{U}^2,\widetilde{U}^3}$.
	Throughout the proof let us fix
	\begin{align*}
		w\coloneqq2\Fkt{f_w}{t}+2^7\Fkt{f_P}{t}.
	\end{align*}
	Recall that we may use \cref{lemma:mcguigmatminors} to transform a $\Bidirected{K_t}$-butterfly minor grasped by $U$ in $D$ into a $K_{t,t}$-matching minor grasped by $W$ in $B$.
	Hence, whenever we find such a butterfly minor in $D$ we may consider the corresponding case to be closed.
	
	\paragraph{Phase I}
	Phase I is divided into $2^{11+8\Fkt{f_P}{t}}\Fkt{f_P}{t}$ \emph{rounds}, each of which produces a \hyperref[def:slice]{slice} $U_i$ of $U_0\coloneqq \widetilde{U}_2$ which is clean with respect to all tiles that have not been marked up to this point.
	We also obtain sets $F_i\subseteq\V{D-\bigcup_{j=1}^{i-1}F_j}$ and $\mathcal{F}_i$ of marked vertices and tiles respectively for each round $i\in[1,2^{11+8\Fkt{f_P}{t}}\Fkt{f_P}{t}]$, and in each round $i$ we will work on the digraph $D_i\coloneqq D_{i-1}-F_i$.
	In this context, whenever we ask for a clean slice of the current slice $U_{i-1}$ or $U'_i$ we ask for a slice $U'$ such that there do not exist $\xi,\xi'\in[1,w+1]$ whose corresponding \hyperref[def:tiling]{tiling} of $U_0$ has a \hyperref[def:tile]{tile} $T$, that is marked or contains a vertex of any separator set found so far, which satisfies $\V{T}\cap\V{U'}\neq\emptyset$.
	
	Each round is divided into two steps, \emph{Step I} and \emph{Step II}.
	Let $i\in[1,2^{11+8\Fkt{f_P}{t}}\Fkt{f_P}{t}]$.
	From the previous round or the initial state, we may assume to have obtained the following sets and graphs, which serve as the \emph{input} for round $i$:
	\begin{itemize}
		\item $F_{I,i-1}\subseteq\V{D}$ such that $\Abs{F_{I,i-1}}\leq \Brace{i-1}\Brace{2^{11}\Brace{w+1}^2\Fkt{f_P}{t}+2^5\Brace{w+1}^2\Fkt{f_P}{t}}$,
		\item $D_{i-1}\coloneqq D-F_{I,i-1}$,
		\item $\mathcal{F}_{I,i-1}\subseteq \bigcup_{\xi,\xi'\in[1,w+1]\Tiling_{U_0,d_1,w,\xi,\xi'}}$ of size at most $\Brace{i-1}\Brace{2^{11}\Brace{w+1}^2\Fkt{f_P}{t}+2^5\Brace{w+1}^2\Fkt{f_P}{t}}$ and
		\item a slice $U_{i-1}$ of $U_0$ of width $\Brace{\Brace{2^{12}\Brace{w+1}^3\Fkt{f_P}{t}+1}\Brace{2^6\Brace{w+1}^3\Fkt{f_P}{t}+1}}^{2^{11+8\Fkt{f_P}{t}}\Fkt{f_P}{t}-i+1}d_2$ that is clean with respect to $F_{I,i-1}$ and $\mathcal{F}_{I,i-1}$ such that, in case $i-1\geq 1$, every long jump over $U_{i-1}$ in $D_{i-1}$ contains a vertex of some tile in $\mathcal{F}_{I,j}\setminus\mathcal{F}_{I,j-1}$ for every $j\in[1,i-1]$.
	\end{itemize}
	For $i=1$ the inputs are the graphs $D_0\coloneqq D$, $U_0$, and two empty sets.
	
	After round $i$ is complete we require the following sets and graphs as its \emph{output}:
	\begin{itemize}
		\item $F_{I,i}\subseteq\V{D}$ such that $\Abs{F_{I,i}}\leq i\Brace{2^{11}\Brace{w+1}^2\Fkt{f_P}{t}+2^5\Brace{w+1}^2\Fkt{f_P}{t}}$,
		\item $D_{i}\coloneqq D-F_{I,i}$,
		\item $\mathcal{F}_{I,i}\subseteq \bigcup_{\xi,\xi'\in[1,w+1]\Tiling_{U_0,d_1,w,\xi,\xi'}}$ of size at most $i\Brace{2^{11}\Brace{w+1}^2\Fkt{f_P}{t}+2^5\Brace{w+1}^2\Fkt{f_P}{t}}$ and
		\item a slice $U_i$ of $U_0$ of width $\Brace{\Brace{2^{12}\Brace{w+1}^3\Fkt{f_P}{t}+1}\Brace{2^6\Brace{w+1}^3\Fkt{f_P}{t}+1}}^{2^{11+8\Fkt{f_P}{t}}\Fkt{f_P}{t}-i}d_2$ that is clean with respect to $F_{I,i}$ and $\mathcal{F}_{I,i}$ such that every long jump over $U_{i}$ in $D_{i}$ contains a vertex of some tile in $\mathcal{F}_{I,j}\setminus\mathcal{F}_{I,j-1}$ for every $j\in[1,i]$.
	\end{itemize}
	
	To be able to find a slice of width $d_2$ in the end we therefore must fix
	\begin{align*}
		d_1\geq \Brace{\Brace{2^{12}\Brace{w+1}^3\Fkt{f_P}{t}+1}\Brace{2^6\Brace{w+1}^3\Fkt{f_P}{t}+1}}^{2^{11+8\Fkt{f_P}{t}}\Fkt{f_P}{t}}d_2,
	\end{align*}
	and we further assume $d_2\geq 2^{16}\Fkt{f_W}{t}^2$ to make sure we can apply \cref{lemma:longjumps1} and \cref{lemma:longjumpspahse2} in every round.
	Note that this is not our final lower bound on $d_2$, just an intermediate assumption.
	
	Next we describe the steps we perform in every round.
	Let $i\in[1,2^{11+8\Fkt{f_P}{t}}\Fkt{f_P}{t}]$ and suppose we are given sets $F_{I,i-1}$, $\mathcal{F}_{I,i-1}$ and graphs $D_{i-1}$, $U_{i-1}$ as required by the input conditions for round $i$.
	
	\emph{Step I}:
	Let $k_i\coloneqq\Brace{\Brace{2^{12}\Brace{w+1}^3\Fkt{f_P}{t}+1}\Brace{2^6\Brace{w+1}^3\Fkt{f_P}{t}+1}}^{2^{11+8\Fkt{f_P}{t}}\Fkt{f_P}{t}-i+1}d_2$.
	For each of the two possible parametrisations of $U$, we consider for every possible choice of $\xi,\xi'\in[1,w+1]$, the \hyperref[def:tiling]{tiling} $\Tiling\coloneqq\Tiling_{U_{i-1},k_i,w,\xi,\xi'}$ together with its four colouring $\Set{\Class_1,\dots,\Class_4}$.
	For each $j\in[1,4]$ we consider the smallest slice $U'$ of $U$ that contains all vertices which belong to some tile of $\Tiling$.
	Consider the auxiliary digraph $\Fkt{D_j^1}{U'}$ of type I obtained from $D_{i-1}$ and define the sets
	\begin{align*}
		X_{\text{I}}^{\text{out}}&\coloneqq\CondSet{x_T^{\text{out}}}{T\in\Class_j}\text{, and}\\
		X_{\text{I}}^{\text{in}}&\coloneqq\CondSet{x_T^{\text{in}}}{T\in\Class_j}.
	\end{align*}
	Additionally we construct the set $Y_{\text{I}}$ as follows:
	Let $Q$ and $Q'$ be the two cycles of $\Perimeter{U_0}$.
	For every $j\in\left[ 1,\frac{3d_1}{4} \right]$, $Y_{\text{I}}$ contains exactly one vertex of $Q\cap P^1_{4j}$, $Q\cap P^2_{4j+2}$, $Q'\cap P^1_{4j}$, and $Q'\cap P^2_{4j+2}$ each.
	Then remove all vertices of $Y_{\text{I}}$ that do not belong to $D_{i-1}$.
	Note that, by choice of $d_1$ and the bound on $F_{I,i-1}$, this does not decrease the size of $Y_{\text{I}}$ dramatically.
	
	Then, if there is a family of $2^7\Fkt{f_P}{t}$ pairwise disjoint directed $X_{\text{I}}^{\text{out}}$-$Y_{\text{I}}$-paths in $\Fkt{D_j^1}{U'}$, \cref{lemma:longjumps1} provides us with the existence of a $K_{t,t}$-\hyperref[def:matchingminor]{matching minor} grasped by \hyperref[def:split]{$\Split{U}$} and we are done.
	So we may assume that there does not exist such a family and thus we may find a set $Z_1\subseteq\V{\Fkt{D_j^1}{U'}}$ of size at most $2^7\Fkt{f_P}{t}$ that meets all such paths by \cref{thm:directedlocalmenger}.
	With a similar argument, we either find a $K_{t,t}$-matching minor grasped by $\Split{U}$, or a set $Z_2\subseteq\V{\Fkt{D_j^1}{U'}}$ of size at most $2^7\Fkt{f_P}{t}$ that meets all directed $Y_{\text{I}}$-$X_{\text{I}}^{\text{in}}$-paths in $\Fkt{D_j^1}{U'}$.
	Let $\pi\in[1,2]$ indicate which of the two parametrisations of $U$ we are currently considering.
	Then we can define the following two sets:
	\begin{align*}
		Z_{\pi,\xi,\xi',j}\coloneqq& \Brace{Z_1\cup Z_2}\cap\V{D}\text{, and}\\
		\mathcal{Z}_{\pi,\xi,\xi',j}\coloneqq& \CondSet{T\in\Tiling}{\Brace{\V{T}\cup\Set{x_T^{\text{out}},x_T^{\text{in}}}}\cap\Brace{Z_1\cup Z_2}\neq\emptyset}.
	\end{align*}
	Note that $\max\Set{\Abs{Z_{\pi,\xi,\xi',j}},\Abs{\mathcal{Z}_{\pi,\xi,\xi',j}}}\leq 2^8\Fkt{f_P}{t}$.
	
	If we do not find a $K_{t,t}$-matching minor grasped by $\Split{U}$ at any point, the sets $Z_{\pi,\xi,\xi',j}$ and $\mathcal{Z}_{\pi,\xi,\xi',j}$ are well defined for every possible choice of $\pi\in[1,2]$, $\xi,\xi'\in[1,w+1]$, and $j\in[1,4]$.
	We collect these sets into the following to sets, which make up the first step towards the output of round $i$:
	\begin{align*}
		F'_{I,i}\coloneqq& \bigcup_{\pi\in[1,2]}\bigcup_{\xi,\xi'\in[1,w+1]}\bigcup_{j\in[1,4]}Z_{\pi,\xi,\xi',j}\text{, and}\\
		\mathcal{F}'_{I,i}\coloneqq& \bigcup_{\pi\in[1,2]}\bigcup_{\xi,\xi'\in[1,w+1]}\bigcup_{j\in[1,4]}\mathcal{Z}_{\pi,\xi,\xi',j}.
	\end{align*}
	Consequently we have that $\max\Set{\Abs{F'_{I,i}},\Abs{\mathcal{F}'_{I,i}}}\leq 2^{11}\Brace{w+1}^2\Fkt{f_P}{t}$.
	
	We may now find a clear slice $U'_i\subseteq U_{i-1}$ of width $\Brace{2^{12}\Brace{w+1}^3\Fkt{f_P}{t}+1}^{2^{11+8\Fkt{f_P}{t}}\Fkt{f_P}{t}-i}\Brace{2^6\Brace{w+1}^3\Fkt{f_P}{t}+1}^{2^{11+8\Fkt{f_P}{t}}\Fkt{f_P}{t}-i+1}d_2$ which does not contain a marked vertex.
	Note that we lose the additional factor of $2\Brace{w+1}$ since we remove whole tiles of width $w$ from $U_{i-1}$.
	Let $D_i'\coloneqq D_{i-1}-F'_{I,i}$.
	This concludes Step I of round $i$.
	
	\begin{claim}\label{claim:stepIclaim}
		Every long jump $J$ over $U'_i$ in $D'_i$ whose endpoints belong to tiles of different colour contains a vertex of a tile from $\mathcal{F}_{I,j}\setminus \mathcal{F}_{I,j-1}$ for every $j\in[1,i-1]$, and it contains a vertex of a tile from $\mathcal{F}'_{I,i}$.
	\end{claim}
	
	\emph{Proof of \Cref{claim:stepIclaim}:} Suppose $J$ is also a long jump over $U_{i-1}$ in $D_{i-1}$, then, as $J$ still exists in $D'_{I,i}$, the tile in whose centre $J$ starts, or the tile in whose centre $J$ ends for some choices of $\pi\in[1,2]$, $\xi,\xi'\in[1,w+1]$, and $j\in[1,4]$, must be marked and therefore cannot belong to $U'_i$.
	
	Hence $J$ must contain some vertex of $U_{i-1}$ as an internal vertex.
	Let $T_{\mathsf{s}}$ be the tile of $U'_i$ in whose centre $J$ starts, and let $T$ be the first tile from the same tiling of $U_{i-1}$, that $J$ meets after $T_{\mathsf{s}}$.
	Let $J'$ be the shortest subpath of $J$ with endpoints in $T_{\mathsf{s}}$ and $T$.
	Then $J'$ is a long jump over $U_{i-1}$ in $D_{i-1}$.
	Therefore, by our assumptions on the input of the $i$th round of Phase I, the first part of our claim is satisfied.
	Moreover, if $T$ has a different colour than $T_{\mathsf{s}}$, $T$ must be marked.
	So suppose $T$ has the same colour as $T_{\mathsf{s}}$.
	Nonetheless, since $T_{\mathsf{s}}$ and the tile $T_{\mathsf{t}}$ which contains the endpoint of $J$ in the current tiling have different colours, $J$ contains a directed subpath $J''$ which is a long jump over $U_{i-1}$ and attaches to tiles of different colour.
	Hence our claim follows.\hfill$\blacksquare$
	
	With this we are ready for Step II of round $i$.
	
	\emph{Step II}:
	For this step let $k_i\coloneqq \Brace{2^{12}\Brace{w+1}^3\Fkt{f_P}{t}+1}^{2^{11+8\Fkt{f_P}{t}}\Fkt{f_P}{t}-i}\Brace{2^6\Brace{w+1}^3\Fkt{f_P}{t}+1}^{2^{11+8\Fkt{f_P}{t}}\Fkt{f_P}{t}-i+1}d_2$.
	We are mainly concerned with the digraph $D'_i$.
	In Step II it suffices to fix one parametrisation of $U$ since the construction of the \hyperref[def:auxiliarydigraphII]{type II auxiliary digraph} leaves the complete interior of same-colour-tiles intact instead of only their \hyperref[def:tile]{centres}.
	For every pair of $\xi,\xi'\in[1,w+1]$ let us consider the \hyperref[def:tiling]{tiling} $\Tiling\coloneqq\Tiling_{U'_i,k_i,w,\xi,\xi'}$ together with its four colouring $\Set{\Class_1,\dots,\Class_4}$.
	Then, for every $j\in[1,4]$ we call upon \cref{lemma:longjumpspahse2} which either provides us with a $K_{t,t}$-matching minor grasped by $\Split{U}$, and therefore closes the proof, or it produces two sets
	\begin{align*}
		Z^{2}_{\xi,\xi',j}\subseteq&\V{D'_i}\text{ of size at most }2^3\Fkt{f_P}{t}\text{, and}\\
		\mathcal{Z}^{2}_{\xi,\xi',j}\subseteq&\Tiling\text{ of size at most }2^3\Fkt{f_P}{t},
	\end{align*}
	such that every directed $\V{U'_i}$-path whose endpoints belong to different tiles of $\Class_j$, contains a vertex of some tile in $\mathcal{Z}^{2}_{\xi,\xi',j}$.
	This allows us to form the sets for the second step of round $i$:
	\begin{align*}
		F''_{I,i}\coloneqq& \bigcup_{\xi,\xi'\in[1,w+1]}\bigcup_{j\in[1,4]}Z^2_{\xi,\xi',j}\text{, and}\\
		\mathcal{F}''_{I,i}\coloneqq& \bigcup_{\xi,\xi'\in[1,w+1]}\bigcup_{j\in[1,4]}\mathcal{Z}^2_{\xi,\xi',j}.
	\end{align*}
	As a result we obtain $\max\Set{\Abs{F''_{I,i}},\Abs{\mathcal{F}''_{I,i}}}\leq 2^5\Brace{w+1}^2\Fkt{f_P}{t}$, and we are able to produce the two sets which will be passed on to the next round.
	\begin{align*}
		F_{I,i}\coloneqq& F'_{I,i}\cup F''_{I,i}\cup F_{I,i-1}\text{, and}\\
		\mathcal{F}_{I,i}\coloneqq& \mathcal{F}'_{I,i}\cup \mathcal{F}''_{I,i}\cup \mathcal{F}_{I,i-1}.
	\end{align*}
	The bounds on $F_{I,i}$ and $\mathcal{F}_{I,i}$ follow immediately from the bounds on $F'_{I,i}$ and $F''_{I,i}$, $\mathcal{F}'_{I,i}$ and $\mathcal{F}''_{I,i}$, and the assumptions on the input of round $i$ respectively.
	
	The pigeon hole principle allows us to find a clear slice $U_i\subseteq U'_i$ of width $\Brace{\Brace{2^{12}\Brace{w+1}^3\Fkt{f_P}{t}+1}\Brace{2^6\Brace{w+1}^3\Fkt{f_P}{t}+1}}^{2^{11+8\Fkt{f_P}{t}}\Fkt{f_P}{t}-i}d_2$ which does not contain a marked vertex.
	Similar to Step I, we lose the additional factor of $2\Brace{w+1}$ since we remove whole tiles of width $w$ from $U'_i$.
	Let $D_i\coloneqq D'_i-F''_{I,i}$.
	This concludes Step II of round $i$.
	
	\begin{claim}\label{claim:stepIIclaim}
		Every long jump over $U_i$ in $D_i$ contains a vertex of some tile in $\mathcal{F}_{I,j}\setminus\mathcal{F}_{I,j-1}$ for every $j\in[1,i]$.
	\end{claim}
	
	\emph{Proof of \Cref{claim:stepIIclaim}:}
	Let $J$ be a long jump over $U_i$ in $D_i$, and let $\Tiling$ be a tiling of $U_0$ defined by $w$ and some $\xi,\xi'\in[1,w+1]$ such that $J$ starts at the centre of some tile $T_{\mathsf{s}}\in\Tiling$.
	Suppose all tiles of $\Tiling$ which contain vertices of $J$ belong to the same colour.
	Then $J$ must have existed during the corresponding part of Step II of round $i$ and thus either $T_{\mathsf{s}}$ or $T_{\mathsf{t}}\in\Tiling$, which is the tile that contains the endpoint of $J$, must have been marked.
	Therefore $J$ must contain at least one tile of a colour different from the one of $T_{\mathsf{s}}$.
	Moreover, we may assume $T_{\mathsf{s}}$ and $T_{\mathsf{t}}$ to be of the same colour as otherwise we would be done by \cref{claim:stepIclaim}.
	Next suppose $J$ is also a long jump over $U_{i-1}$, then again $J$ would have been considered during Step II as a long jump connecting two tiles of the same colour and thus $T_{\mathsf{s}}$ or $T_{\mathsf{t}}$ would have been marked.
	Therefore $J$ must contain a vertex of some tile from $U_{i-1}$.
	Let $J'$ be a shortest subpath from $T_{\mathsf{s}}$ to some tile $T$ of $U_{i-1}$, then $J'$ is a long jump over $U_{i-1}$ and thus $J$ contains a vertex of some tile of $\mathcal{F}_{I,j}\setminus\mathcal{F}_{I,j-1}$ for every $j\in[1,i-1]$ by our assumptions on the input of round $i$.
	If $T$ has a different colour than $T_{\mathsf{s}}$, then $T$ would have been marked in Step I of round $i$, and if $T$ shares the colour of $T_{\mathsf{s}}$ it must have been marked in Step II of round $i$.
	Either way our claim follows. \hfill$\blacksquare$
	
	From \cref{claim:stepIIclaim} it follows that we satisfy all requirements for the output of round $i$ and thus round $i$ is complete.
	We continue until we finish round $i=2^{11+8\Fkt{f_P}{t}}\Fkt{f_P}{t}$ and obtain the following four objects as its output:
	\begin{itemize}
		\item a slice $U_I\coloneqq U_{2^{11+8\Fkt{f_P}{t}}\Fkt{f_P}{t}}$ of width $d_2$,
		\item a set $A\coloneqq F_{I,2^{11+8\Fkt{f_P}{t}}\Fkt{f_P}{t}}$ of size at most $t^{28}2^{60+2^{10}t^8}$,
		\item a digraph $D_I\coloneqq D_{2^{11+8\Fkt{f_P}{t}}\Fkt{f_P}{t}}=D-A$, and
		\item a sequence $\mathcal{F}_{I,1}\subseteq \mathcal{F}_{I,2}\subseteq\dots\subseteq \mathcal{F}_{I,2^{11+8\Fkt{f_P}{t}}\Fkt{f_P}{t}}$ such that every long jump over $U_I$ in $D_I$ contains a vertex of some tile in $\mathcal{F}_{I,i}\setminus\mathcal{F}_{I,i-1}$ for every $i\in[1,2^{11+8\Fkt{f_P}{t}}\Fkt{f_P}{t}]$.
	\end{itemize}
	This brings us to the final claim of Phase I.
	
	\begin{claim}\label{claim:nolongjumps}
		If there is a long jump over $U_I$ in $D_I$, then there is a $K_{t,t}$-\hyperref[def:matchingminor]{matching minor} grasped by $\Split{U}$ in $D$.
	\end{claim}
	
	\emph{Proof of \Cref{claim:nolongjumps}:}
	Let $J$ be a long jump over $U_I$ in $D_I$.
	We fix a parametrisation of $U$, $\xi,\xi'\in[1,w+1]$, and $c\in[1,4]$ such that there exists a tile $T_{\mathsf{s}}\in\Tiling\coloneqq\Tiling_{U_0,d_1,w,\xi,\xi'}$ of colour $c$ whose centre contains the starting point of $J$.
	Let $T_{\mathsf{t}}\in\Tiling$ be the tile that contains the endpoint of $J$.
	As $J$ is a long jump, note that $T_{\mathsf{s}}\neq T_{\mathsf{t}}$.
	
	We now create a family $\mathcal{L}_0$ of $2^{9+8\Fkt{f_P}{t}}\Fkt{f_P}{t}$ pairwise disjoint subpaths of $J$ with the following properties:
	\begin{enumerate}
		\item for every $L\in\mathcal{L}_0$, let $T_{L,1}$ and $T_{L,2}$ be the tiles of $\Tiling$ that contain the starting point $s_L$ and the endpoint $t_L$ of $L$ respectively, then there exist distinct $i_{L,1},i_{L,2}\in[1,2^{11+8\Fkt{f_P}{t}}\Fkt{f_P}{t}]$ such that $s_L$ is a vertex of a tile from $\mathcal{F}_{I,i_{L,1}}\setminus\mathcal{F}_{I,i_{L,1}-1}$, and $t_L$ is a vertex of some tile in $\mathcal{F}_{I,i_{L,2}}\setminus\mathcal{F}_{I,i_{L,2}-1}$, and
		\item if $L,L'\in\mathcal{L}_0$ are distinct, then $\Set{i_{L,1},i_{L,2}}\cap\Set{i_{L',1},i_{L',2}}=\emptyset$.
	\end{enumerate}
	Let us initialise $\mathcal{I}_0\coloneqq[1,2^{11+8\Fkt{f_P}{t}}\Fkt{f_P}{t}]$ and for every subset $\mathcal{I}'\subseteq\mathcal{I}_0$, we define the family $\mathcal{F}_{\mathcal{I}'}\coloneqq \bigcup_{i\in\mathcal{I}'}\mathcal{F}_{I,i}\setminus\mathcal{F}_{I,i-1}$.
	Please note that every internal vertex of $J$ that belongs to $U$ must belong to some tile from $\mathcal{F}_{I,2^{11+8\Fkt{f_P}{t}}\Fkt{f_P}{t}}$, since otherwise we could find a directed path from the centre of $T_{\mathsf{s}}$ to the perimeter of $U_0$ contradicting the construction in Step I of Phase I, or we would have a directed path between two tiles of the same colour, where both of them are unmarked.
	This second outcome contradicts the construction in Step II of Phase I.
	
	Consider the shortest subpath of $J$ that shares its starting point with $J$ and is a long jump over $U$.
	Let $T_{L_1,1}$ be the tile of $\Tiling$ where this path ends and let $s_{L_1}$ be the first vertex of $J$ for which its successor along $J$ does not belong to $T_{L_1,1}$.
	Note that there exists $i_{L_1,1}\in\mathcal{I}_0$ such that $T_{L_1,1}\in\mathcal{F}_{I,i_{L_1,1}}\setminus \mathcal{F}_{I,i_{L_1,1}-1}$ by the discussion above.
	Let $L_1$ be the shortest subpath of $J$ that starts in $s_{L_1}$ and ends in a vertex $t_{L_1}$ for which $i_{L_1,2}\in\mathcal{I}_0\setminus\Set{i_{L_1,1}}$ exists such that $t_{L_1}$ belongs to a tile from $\mathcal{F}_{I,i_{L_1,2}}\setminus \mathcal{F}_{I,i_{L_1,2}-1}$.
	Let $T_{L_1,2}$ be the tile of $\Tiling$ that contains $t_{L_1}$.
	We add $L_1$ to $\mathcal{L}_0$ and set $\mathcal{I}_1\coloneqq \mathcal{I}_0\setminus\Set{i_{L_1,1},i_{L_1,2}}$.
	By our choice of $i_{L_1,1}$ and $i_{L_1,2}$, the path $t_{L_1}J$ still contains a vertex of a tile from $\mathcal{F}_{I,j}\setminus\mathcal{F}_{I,j-1}$ for every $j\in\mathcal{I}_1$.
	
	Now let $q\in[2,2^{9+8\Fkt{f_P}{t}}\Fkt{f_P}{t}]$ and assume that the paths $L_1,\dots,L_{q-1}$ together with the tiles, indices and the set $\mathcal{I}_{q-1}$ have already been constructed.
	Follow along $J$, starting from $t_{L_{q-1}}$, until the next time we encounter the last vertex $s_{L_q}$ of some tile from $\mathcal{F}_{\mathcal{I}_{q-1}}$ before $J$ leaves said tile again.
	Let $T_{L_q,1}\in\Tiling$ be the tile that contains $s_{L_q}$, and let $i_{L_q,1}\in\mathcal{I}_{q-1}$ be the integer such that $s_{L_q}$ belongs to a tile of $\mathcal{F}_{I,i_{L_q,q}}\setminus \mathcal{F}_{I,i_{L_q,q}-1}$.
	Then let $L_q$ be the shortest subpath of $J$ that starts in $s_{L_q}$ and ends in a vertex $t_{L_q}$ which belongs to a tile from $\mathcal{F}_{I,i_{L_q,2}}\setminus \mathcal{F}_{I,i_{L_q,2}-1}$, where $i_{L_q,2}\in\mathcal{I}_{q-1}\setminus\Set{i_{L-q,1}}$.
	We choose $T_{L_q,2}\in\Tiling$ to be the tile that contains $t_{L_q}$ and set $\mathcal{I}_q\coloneqq \mathcal{I}_{q-1}\setminus\Set{i_{L_q,1},i_{L_q,2}}$.
	As before notice that $t_{L_q}J$ still contains a vertex from some tile in $\mathcal{F}_{I,j}\setminus\mathcal{F}_{I,j-1}$ for every $j\in\mathcal{I}_q$.
	Add $L_q$ to $\mathcal{L}_0$.
	
	With every iteration we remove exactly two members from $\mathcal{I}$ and, as $\Abs{\mathcal{I}}=2^{11+8\Fkt{f_P}{t}}\Fkt{f_P}{t}$, this means that by the time we reach some $q$ for which $\mathcal{I}_q=\emptyset$, we have indeed constructed $2^{10+8\Fkt{f_P}{t}}\Fkt{f_P}{t}$ paths as required.
	
	Now there must exist $c'\in[1,4]$ and a family $\mathcal{L}_1\subseteq \mathcal{L}_0$ of size $2^{8+8\Fkt{f_P}{t}}\Fkt{f_P}{t}$ such that each path $L\in\mathcal{L}_1$ has at least one endpoint in $\Class_{c'}$.
	And, in an immediate second step, we can find a family $\mathcal{L}_2\subseteq\mathcal{L}_1$ of size $2^{7+8\Fkt{f_P}{t}}\Fkt{f_P}{t}$ such that every path in $\mathcal{L}_2$ starts in a tile of $\Class_{c'}$, or every path in $\mathcal{L}_2$ ends in a tile of $\Class_{c'}$.
	Without loss of generality we may assume that every path in $\mathcal{L}_2$ starts in a tile of $\Class_{c'}$ since the other case follows with similar arguments.
	
	Let $\widetilde{U}'$ be the smallest slice of $U$ that contains all tiles from $\Class_{c'}$, but no tile from $\Class_{c'}$ meets the perimeter of $\widetilde{U}$.
	Then let $\widetilde{\Tiling}\coloneqq\TierIITiling{\Tiling,c',w}{\widetilde{U}'}$ be the \hyperref[def:tilingII]{tier II tiling} of $\widetilde{U}\coloneqq\InducedSubgraph{\widetilde{U}'}{\Tiling,c',w}$.
	Since the paths in $\mathcal{L}_2$ are pairwise disjoint, we can extend each $L\in\mathcal{L}_2$ such that it starts on the centre of the tile of $\widetilde{\Tiling}$ which encloses its endpoint in $U$, while making sure that the resulting family of paths is still at least half-integral.
	Similarly, wherever necessary, we may extend the paths through $U$ such that each of them also ends in a tile of $\widetilde{\Tiling}$.
	Indeed, we can even guarantee that the endpoints of the resulting paths are mutually at $\widetilde{U}$-distance at least $4$.
	Let $\mathcal{L}_3$ be the resulting half integral linkage.
	
	Next consider the four colouring $\Set{\widetilde{\Class_1},\dots,\widetilde{\Class_4}}$ of $\widetilde{\Tiling}$.
	Then there exists $\widetilde{c}\in[1,4]$ and a family $\mathcal{L}_4\subseteq\mathcal{L}_3$ of size $2^{5+8\Fkt{f_P}{t}}\Fkt{f_P}{t}$ such that every path in $\mathcal{L}_4$ starts at the centre of some tile from $\Class_{\widetilde{c}}$.
	It follows from the construction of $\mathcal{L}_0$ that no two paths in $\mathcal{L}_4$ start in the same tile.
	
	By a similar argument, there exists a family $\mathcal{L}_5\subseteq\mathcal{L}_4$ of size $2^{4+8\Fkt{f_P}{t}}\Fkt{f_P}{t}$ such that either none, or all paths in $\mathcal{L}_5$ end in tiles of $\widetilde{\Class_{\widetilde{c}}}$.
	
	In the first case we can extend every path in $\mathcal{L}_5$ towards the perimeter of $\widetilde{U}$ such that the resulting family $\mathcal{L}_6$ of paths remains at worst half-integral, and the endpoints of the resulting paths are mutually at $\widetilde{U}$-distance at least $4$.
	Now \cref{thm:halfintegral} provides us with a family $\mathcal{L}_7$ of size $2^{3+8\Fkt{f_P}{t}}\Fkt{f_P}{t}$ such that $\V{\mathcal{L}_7}\subseteq\V{\mathcal{L}_6}$, and the paths in $\mathcal{L}_7$ are pairwise vertex disjoint.
	Hence \cref{lemma:longjumps1} yields the existence of a $K_{t,t}$-matching minor grasped by $\Split{\widetilde{U}}$ and our claim follows.
	
	In the second case we consider two subcases.
	Let $\mathcal{X}$ be the family of all tiles of $\widetilde{\Tiling}\setminus\widetilde{\Class_{\widetilde{c}}}$ which contain an internal vertex of some path in $\mathcal{L}_5$ but no endpoint of any path in $\mathcal{L}_5$.
	
	Suppose $\Abs{\mathcal{X}}\geq 2^8\Fkt{f_P}{t}$ and recall that, if $L$ and $P$ are directed paths, we say that $P$ is a \emph{long jump of $L$} if $P$ is a $w$-long jump over $U$ and $P\subseteq L$.
	We also say that $P$ is a \emph{jump of $L$}, if $P$ is a directed $\V{U}$-path.
	Then we can use the technique from the first part of the proof of \cref{lemma:longjumps1} to construct a half-integral family $\mathcal{L}_6$ such that
	\begin{enumerate}
		\item $\Abs{\mathcal{L}_6}=2^{4+8\Fkt{f_P}{t}}\Fkt{f_P}{t}$, and
		\item for every $L\in\mathcal{L}_6$, every endpoint $u$ of a jump of $L$ with $u\in\V{\widetilde{U}}$ belongs to a tile from $\widetilde{\Class_{\widetilde{c}}}\cup\mathcal{X}$.
	\end{enumerate}
	We can then apply \cref{thm:halfintegral} to obtain a family $\mathcal{L}_7$ of size $2^{3+8\Fkt{f_P}{t}}\Fkt{f_P}{t}$ with $\V{\mathcal{L}_7}\subseteq\V{\mathcal{L}_6}$ such that the paths in $\mathcal{L}_7$ are pairwise disjoint and link the same two sets of vertices as the paths in $\mathcal{L}_6$ do.
	Finally \cref{lemma:longjumps1} yields the existence of a $K_{t,t}$-matching minor grasped by $\Split{U}$.
	
	So we may assume $\Abs{\mathcal{X}}<2^8\Fkt{f_P}{t}$.
	In this case, we may find a subwall $\widetilde{U}'$ of $\widetilde{U}$ of order $d_1-2^8\Fkt{f_P}{t}\Brace{2w+1}$ that does not contain a vertex of any tile in $\mathcal{X}$ by removing, for every tile $T\in\mathcal{X}$, all edges and vertices of the horizontal cycles and vertical paths of $T$, which are not used by other cycles or paths.
	For each tile we remove during this procedure, we remove a row and a column of tiles and thereby might reduce the number of distinct tiles which contain starting vertices of paths in $\mathcal{L}_5$ by a factor of $\frac{1}{2}$.
	However, since $\Abs{\mathcal{X}}<2^8\Fkt{f_P}{t}$, we can still find, after potentially expanding the start and end sections of some paths to again reach the slightly shifted perimeters of their tiles, a half-integral family $\mathcal{L}_6$ of size $2^4\Fkt{f_P}{t}$ of paths that start and end in tiles of $\widetilde{\Class_{\widetilde{c}}}$ and that are otherwise disjoint from $\widetilde{U}'$.
	By using \cref{thm:halfintegral} we can transform this family into an integral family $\mathcal{L}_7$ of size $2^3\Fkt{f_P}{t}$ and thus an application of \cref{lemma:longjumpspahse2} yields a $K_{t,t}$-matching minor grasped by $\Split{U}$\hfill$\blacksquare$

	Concluding Phase I, \cref{claim:nolongjumps} either yields a $K_{t,t}$-matching minor grasped by $\Split{U}$ and therefore closes the case, or $U_I$ is in fact clean, meaning that $U_I$ has no long jump in $D_I$.
	The set $A$ is a set of vertices of $D$ which means that it can be seen as a set of edges from $M$ in $B$.
	Hence we may treat $A$ as an $M$-conformal set of vertices in $B$ of size at most $t^{28}2^{61+2^{10}t^8}$.
	Consequently we may bound the function $\ApexBound$ from the statement of \cref{thm:matchingflatwall} as follows:
	\begin{align*}
		\Fkt{\ApexBound}{t}\leq t^{28}2^{61+2^{10}t^8}.
	\end{align*}
	
	\paragraph{Phase II}
	
	Until now all of our efforts were focused on the removal of long jumps, so the only jumps that remain must be short, i.\@e.\@ their endpoints in $U_I$ must have $U_I$-distance at most $w-1$.
	To meet the requirements from Phase I and have enough space left in $U_I$, let us make the following assumption:
	\begin{align*}
		d_2\geq t^4\Brace{r+2^{32+t^{30}}}.
	\end{align*}
	
	Let us partition $U_I$ into $2t^4$ many slices as follows:
	First we can partition $U_I$ into $t^4$ slices $S_i$ of width $r+2^{32+t^{30}}$.
	Each $S_i$ is then partitioned into a slice $H_i$ of width $r+2^{15+t^30}$ that contains the left perimeter cycle of the of $S_i$, and a slice $G_i$ of width $2^{31+t^{30}}$ containing the right cycle of the perimeter of $S_i$.
	For every $i\in[1,t^4]$ we may now further partition $H_i$.
	Let $N_{i,L}\subseteq H_i$ be the slice of width $2^{30+t^{30}}$ containing the left cycle of $\Perimeter{H_i}$, let $N_{i,R}$ be the slice of width $2^{30+t^{30}}$ containing the right cycle of $\Perimeter{H_i}$, and let $N_i\coloneqq H_i-N_{L,i}-N_{R,i}$ be the remaining slice of width $r$.
	
	\begin{claim}\label{claim:crossesmeanKt}
		If for all $i\in[1,t^4]$ the slice $\Split{N_i}$ is \textbf{not} \hyperref[def:matchingflat]{$\Perimeter{\Split{N_i}}$-flat} in $B$ with respect to $A$, then there exists a $K_{t,t}$-matching minor grasped by $W$ in $B$.
	\end{claim}
	
	\emph{Proof of \Cref{claim:crossesmeanKt}:}
	To start, notice that, since $U_I$ is clean in $D_I$, every internally $M$-conformal path in $B-A$ whose endpoints belong to $\Split{U_I}$ and which is internally disjoint from $\Split{U_I}$ must correspond to a short jump over $U_I$ in $D_I$.
	
	The goal of this proof is to adjust the perfect matching $M$ of $B$ in such a way that the $M'$-direction of the new perfect matching $M'$ yields the desired $\Bidirected{K_t}$-butterfly minor.
	
	Recall our construction of the butterfly minor model of $\Bidirected{K_t}$.

	For every $i\in[1,t^4]$, let $N_i'$ be the slice of width $4w+2+r$ of $H_i$, such that $N_i'-N_i$ consists of exactly two components which both are slices of width $2w+1$.
	Since $\Split{N_i}$ is not $\Perimeter{\Split{N_i}}$-flat in $B$ with respect to $A$, there must exist $\xi,\xi'$ and tiling $\Tiling=\Tiling_{N_i',4w+2+r,w,\xi_,\xi'}$ containing a tile $T_i$ such that $\Split{T_i}$ is non-Pfaffian in $B-A$.
	Note that the centre of $T_i$ must lie within $N_i$, since otherwise we either could choose a different $T_i$, or $\Split{N_i}$ would be $\Perimeter{\Split{T_i}}$-flat in $B$ with respect to $A$.
	
	We may adjust the parametrisation of $U$ and its embedding, such that the top most strip of height $16\Fkt{f_P}{t}+2w$ does not contain a vertex of $T_i$ for any $i\in[1,t^4]$.
	This is to make sure that the top most strip of height $8\Fkt{f_P}{t}$ is not used by any of the preliminary constructions below and can therefore later be used to construct our $\Bidirected{K_t}$-butterfly minor.
	
	Let $U_I''$ be the subgraph of $U_I$ obtained as the union of all vertical cycles and horizontal paths of $U_I$ which do not contain vertices of the $T_i$.
	Note that the resulting graph is again a slice of some wall.
	Then we may find a tile $T_i'$ of width $w$ in $U_I''$ for every $i\in[1,t^4]$, such that the interior of the centre of $T_i'$ in $U_I$ contains the entire tile $T_i$.
	Observe that, since the centres of the $T_i$ lie within the $N_i$, the $T_i'$ are completely contained in the $N_i'$ respectively.
	Let $s_{i,L}$ be the top left, and $s_{i,R}$ be the top right vertex of the perimeter of $T_i'$ as well as $t_{i,L}$ be the bottom left and $t_{i,R}$ be the bottom right vertex of the perimeter of $T_i'$.
	By \cref{lemma:crossdirection} we can find a perfect matching $M_i$ of $\Split{N_i'}\cup\Attachment{B-A,\Split{U_I}}{T_i}$ such that
	\begin{itemize}
		\item after slightly re-routing the horizontal paths of $U_I''$ in $\DirM{\Split{D_I}}{M\Delta M_i}$, $U_I''$ still exists, and
		\item there are paths vertex disjoint directed paths $L_i$ and $R_i$ such that
		\begin{itemize}
			\item $L_i$ starts in $s_{i,L}$ and ends in $t_{i,R}$,
			\item $R_i$ starts in $s_{i,R}$ and ends in $t_{i,L}$, and
			\item $L_i$ and $R_i$ are vertex disjoint from $U_I''-\DirM{\Split{T_i'}}{M_i}$.
		\end{itemize}
	\end{itemize}
	Note that for $i\neq j\in[1,t^4]$, we know that $M_i$ and $M_j$ are disjoint.
	Moreover, $\Split{N_i'}\cup\Attachment{B-A,\Split{U_I}}{T_i}$ and $\Split{N_j'}\cup\Attachment{B-A,\Split{U_I}}{T_j}$ must be disjoint since otherwise we could find an internally $M$-conformal path $P$ starting in $\Split{T_i}$ and ending in $\Split{T_j}$ which is internally disjoint from $\Split{U_I}$.
	In $D_I$ such a path $P$ would correspond to a long jump over $U_I$ by the construction of $H_i$ and $H_j$, but by our assumption, there cannot exist such a long jump.
	Let $M'\coloneqq M\Delta\bigcup_{i=1}^{t^4}M_i$ and $D_{II}'\coloneqq\DirM{\Split{D_I}}{M'}$.
	Hence we can find a slice $U_{II}'$ corresponding to the construction of $U_I''$ for all $i\in[1,t^4]$ simultaneously in $D_{II}'$, and the paths $L_i$, $R_i$ exist for all $i\in[1,t^4]$ such that for all $i\neq j\in[1,t^4]$, $L_i$ and $R_i$ are disjoint from $L_j$ and $R_j$.
	
	For every $i\in[1,t^4]$, let $M'_i$ be the perfect matching of $\Split{G_i}$ obtained by switching $M$ along all horizontal cycles of $\Split{G_i}$.
	At last let $M''$ be the perfect matching of $\Split{U'_{t+1}}$ obtained by switching $M$ along all of its horizontal cycles.
	Then we may set $M_{II}\coloneqq M'\Delta\Brace{\bigcup_{i=1}^{t^4} M'_i \cup M''}$ and consider $D_{II}\coloneqq\DirM{B}{M_{II}}$.
	Let $U_{II}$ be the union of $U-U_I$, the slices of $U_II$ that correspond to the $M_i$-directions of the $H_i$ together with the paths $R_i$ and $L_i$, and the $M_{II}$-directions of the $G_i$.
	Note that $U_{II}$ is indeed a subgraph of $D_{II}$ that closely resembles a cylindrical wall, with a few holes, $t^4$ pairs of crosses over small areas over the wall, and some sections, corresponding the the $G_i$, in which the vertical cycles are oriented upwards instead of downwards.
	See \cref{fig:switchedwall} for a rough sketch of the digraph $U_{II}$.
	
	\begin{figure}[!ht]
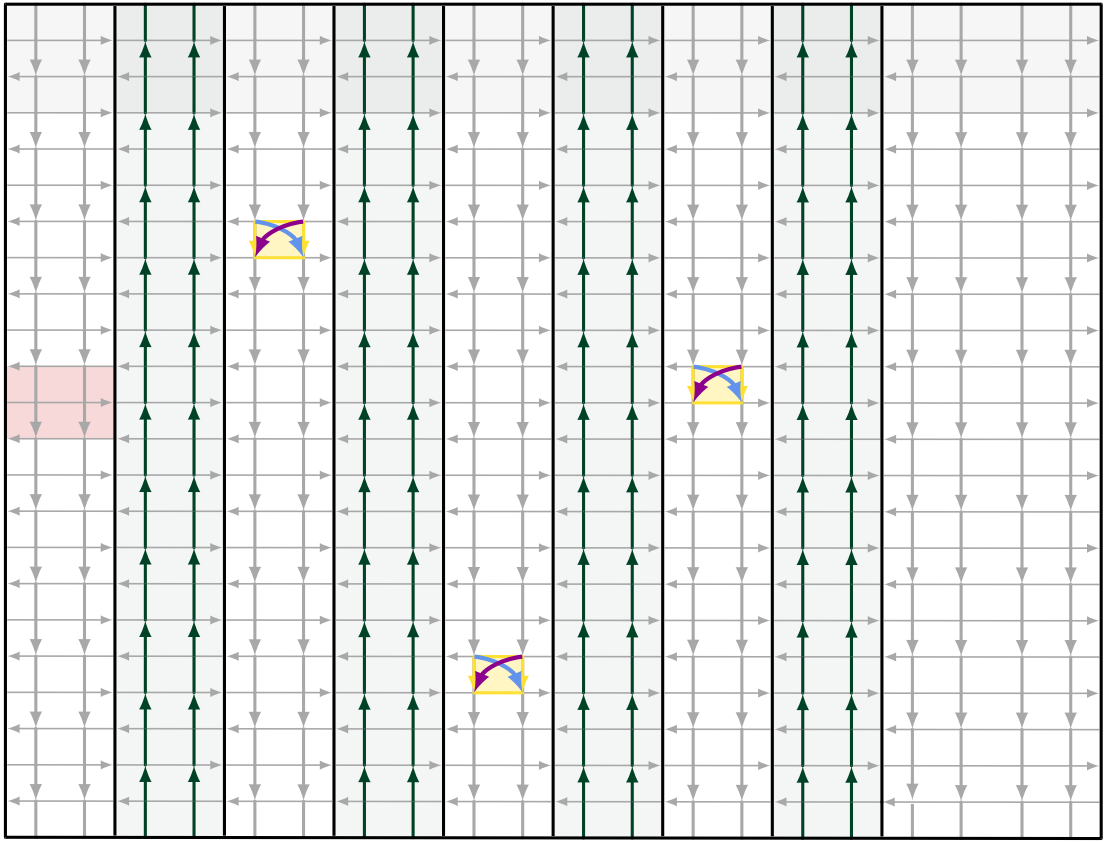

		\centering

		\caption{A sketch of the digraph $U_{II}$ with its crosses, the top strip of bounded height we aim to keep unused, and the section of the left slice which is used for the roots of $\Bidirected{K_t}$ (highlighted in red).}
		\label{fig:switchedwall}
	\end{figure}
	
	As the final step of the proof we construct a $\Bidirected{K_t}$-butterfly minor grasped by $U_{II}$ in $D_{II}$.
	Since $\Split{D_{II}}=B$, \cref{lemma:mcguigmatminors} can then be applied to find a $K_{t,t}$-matching minor grasped by $W$ in $B$, which will conclude the proof of our claim.
	
	Let us position the roots of the vertices for our $\Bidirected{K_t}$ on the path $P^1_{\frac{3}{2}d_1}$.
	For each pair $\Brace{i,j}\in[1,t]^2$, we will construct a path $J_{i,j}$ starting at the head of $f_i^j$ and eventually ending at the tail of $e_j^i$.
	Since the roots are ordered from left to right, we start out with a family of disjoint paths in $U_{II}$ which are ordered from left to right as
	\begin{align*}
		J_{1,1},J_{1,2},\dots,J_{2,1},\dots,J_{t-1,t},J_{t,1},\dots,J_{t,t}
	\end{align*}
	while moving downwards in the current embedding of $U_{II}$.
	We grow the paths $J_{i,j}$ until each of them ends on some $P^1_{\Fkt{h_1}{i,j}}$ where $\Fkt{h_1}{i,j}>\Fkt{h_1}{i',j'}$, if $\Brace{i,j}$ is lexicographically smaller than $\Brace{i',j'}$.
	Then we grow the paths along their respective $P^1_{\Fkt{h_1}{i,j}}$, until each of them ends on a switched vertical cycle $Q_{\Fkt{h_2}{i,j}}$ of $U'_{t+1}$ where $\Fkt{h_2}{i,j}>\Fkt{h_2}{i',j'}$, whenever $\Brace{i,j}$ is lexicographically smaller then $\Fkt{h_2}{i,j}$.
	Afterwards we move each path along its respective horizontal cycle upwards, until each $J_{i,j}$ reaches some $P^1_{\Fkt{h_3}{i,j}}$, where $\Fkt{h_3}{i',j'}<\Fkt{h_3}{i,j}$, if $J_{i,j}$ is further to the right than $J_{i',j'}$ while moving upwards, and each $P^1_{\Fkt{h_3}{i,j}}$ is above the top path of $T'_1$, while avoiding the top strip of height $4\Fkt{f_P}{t}$.
	This is possible by the choice of our parametrisation.
	Next move each $J_{i,j}$ towards the right along its $P^1_{\Fkt{h_3}{i,j}}$ until the following requirements are met:
	\begin{itemize}
		\item the path $J_{1,1}$ ends on a horizontal cycle of $N_{L,1}$, which lies left of $T_1'$,
		\item the two paths $J_{1,2}$ and $J_{1,3}$ meet exactly the two horizontal cycles of $N_1'$ which contain the tail of $L_1$ and $R_1$ respectively, and
		\item the paths $J_{1,4},\dots,J_{t,t}$ end on horizontal cycles of $N_{R,1}$ such that $J_{i,j}$ ends on a cycle left of the cycle on which $J_{i',j'}$ ends, if $J_{i,j}$ passes $T_1'$ below $J_{i',j'}$.
	\end{itemize}
	Now, by growing the paths downwards, we may route $J_{1,2}$ and $J_{1,3}$ through the paths $L_1$ and $R_1$, thereby swapping their relative positions within the whole path family.
	We may then either grow the paths further downwards and then go right again, or use the switched horizontal cycles of $G_1$, to grow them back upwards within $U_{II}$.
	Either way, we are able to route the now neighbouring paths $J_{1,2}$ and $J_{1,4}$ (or $J_{2,1}$ in case $t=3$) through $T_2'$.
	Since we have $t^4$ such crosses in total, we may repeat this process until, after passing through $T_{t^4}'$, the paths are in the following order when moving downwards in our current embedding of $U_{II}$:
	\begin{align*}
		J_{1,1},J_{2,1},J_{3,1},\dots,J_{t,1},J_{1,2},\dots,J_{t,2},\dots,J_{1,t},\dots,J_{t,t}.
	\end{align*}
	We can now use $G_{t^4}$ to route all of these paths upwards and finally use the top strip of height $4\Fkt{f_P}{t}$, which we kept untouched so far, to route the paths back to the roots of our vertices.
	With the changed order we are no able to grow each path $J_{i,j}$ such that it indeed ends on the tail of $e_{j,i}$, while making sure that our paths stay pairwise vertex disjoint.
	Hence we have found a $\Bidirected{K_t}$-butterfly minor as desired and the proof of our claim is complete.\hfill$\blacksquare$

	From \cref{claim:crossesmeanKt} it follows that, in case all $t^4$ slices $\Split{N_i}$ are \textbf{not} $\Perimeter{\Split{N_i}}$-flat in $B$ with respect to $A$, then there exists a $K_{t,t}$-matching minor grasped by $W$ in $B$.
	In this case we are done and thus we may assume that there exists some $i\in[1,t^4]$ such that $\Split{N_i}$ is indeed $\Perimeter{\Split{N_i}}$-flat in $B$ with respect to $A$.
	Since $N_i$ is of width $r$, $\Split{N_i}$ contains a conformal matching $r$-wall $\widetilde{W}$ such that $\Perimeter{\widetilde{W}}=\Perimeter{\Split{N_i}}$ and thus our proof is complete.
	At last let us combine all assumptions on the $d_j$ to obtain the following bound on $\Fkt{\WallBound}{t,r}$:
	\begin{align*}
		\Fkt{\WallBound}{t,r}\leq\Brace{2^{140}t^{76}}^{2^{10}t^8+2^{12}}\Brace{r+2^{32+t^{30}}}.
	\end{align*} 
\end{proof}

Please note that many of the exponential parts of the bounds in $\ApexBound$ and $\WallBound$ are due to specific difficulties encountered in the digraphic setting.
Indeed, even the exponential lower bound on $f_W$ necessary for \cref{lemma:directedlongjumps} can probably be made polynomial by simply diving deeper into the increased freedom provided by the possibility of switching perfect matchings along directed cycles.
At this point it is not clear whether the exponential part of Phase II in the proof of \cref{thm:matchingflatwall} can also be made polynomial.
To do this, it is probably necessary to develop a technique for removing long jumps that does not rely on tilings.


\bibliographystyle{alphaurl}
\bibliography{literature}

\end{document}